\documentclass[12pt]{amsart}
\usepackage{amssymb,latexsym}

\numberwithin{equation}{section}

\newtheorem{theorem}{Theorem}[section]
\newtheorem{proposition}[theorem]{Proposition}
\newtheorem{conjecture}[theorem]{Conjecture}
\newtheorem{corollary}[theorem]{Corollary}
\newtheorem{lemma}[theorem]{Lemma}

\theoremstyle{definition}

\newtheorem{remark}[theorem]{Remark}
\newtheorem{example}[theorem]{Example}
\newtheorem{definition}[theorem]{Definition}

\def\aa{\mathbf{a}}
\def\bb{\mathbf{b}}
\def\cc{\mathbf{c}}
\def\dd{\mathbf{d}}
\def\ee{\mathbf{e}}

\def\gg{\mathbf{g}}

\def\xx{\mathbf{x}}
\def\yy{\mathbf{y}}
\def\zz{\mathbf{z}}
\def\pp{\mathbf{p}}

\def\ZZ{\mathbb{Z}}

\def\PP{\mathbb{P}}
\def\RR{\mathbb{R}}
\def\QQ{\mathbb{Q}}
\def\TT{\mathbb{T}}

\def\Acal{\mathcal{A}}
\def\Fcal{\mathcal{F}}
\def\Mcal{\mathcal{M}}
\def\Xcal{\mathcal{X}}

\def\Qsf{\QQ_{\,\rm sf}}

\def\Aprin{\Acal_\bullet}

\def\Phireal{\Phi^{\rm re}}
\def\Auniv{\Acal^{\rm univ}}

\def\verylongleftrightarrow{\longleftarrow\!\!\!\longrightarrow}

\def\Trop{\operatorname{Trop}}

\def\sgn{\operatorname{sgn}}

\newcommand{\twobyone}[2]{
\begin{bmatrix}#1\\#2
\end{bmatrix}
}

\newcommand{\overunder}[2]{
\!\begin{array}{c}
\scriptstyle{#1}\\[-.1in]
-\!\!\!-\!\!\!-\\[-.1in]
\scriptstyle{#2}
\end{array}
\!
}

\newcommand{\bmat}[4]{\left[\!\!\begin{array}{cc}
#1 & #2 \\ #3 & #4 \\
\end{array}\!\!\right]}

\renewcommand{\eqref}[1]{{\rm (\ref{#1})}}

\newcommand{\column}[2]
{\begin{bmatrix} #1\\[-3pt] \vdots \\[-2pt]
  #2\end{bmatrix} }

\newcommand{\smallcolumn}[2]
{\left[\begin{smallmatrix} #1\\ \cdot \\[-3pt] \cdot \\[-3pt] \cdot \\[-1pt]
  #2\end{smallmatrix}\right] }

\begin{document}


\title[Cluster algebras~IV]
{Cluster algebras~IV:
Coefficients
}

\author{Sergey Fomin}
\address{Department of Mathematics, University of Michigan,
Ann Arbor, MI 48109, USA} \email{fomin@umich.edu}

\author{Andrei Zelevinsky}
\address{\noindent Department of Mathematics, Northeastern University,
 Boston, MA 02115}
\email{andrei@neu.edu}

\subjclass[2000]{Primary
16S99, 
Secondary
05E15, 
22E46
}

\begin{abstract}
We study the dependence of a cluster algebra
on the choice of coefficients.
We write general formulas expressing the cluster variables in any
cluster algebra in terms of the initial data;
these formulas involve a family of polynomials associated with
a particular choice of ``principal'' coefficients.

We show that the exchange graph of a cluster algebra
with principal coefficients covers the exchange graph of any cluster
algebra with the same exchange matrix.

We investigate two families of parametrizations of cluster monomials by lattice
points, determined, respectively, by the denominators of their Laurent
expansions and by certain multi-gradings in cluster algebras with
principal coefficients.
The properties of these parametrizations, some proven and some
conjectural, suggest links to duality conjectures of V.~Fock and
A.~Goncharov.

The coefficient dynamics leads to a natural
generalization of Al.~Zamolodchikov's $Y$-systems.
We establish a Laurent phenomenon for such $Y$-systems,
previously known in finite type only, and sharpen the periodicity
result from an earlier paper.

For cluster algebras of finite type, we
identify a canonical ``universal'' choice of coefficients such that
an arbitrary cluster algebra can be obtained from
the universal one (of the same type) by an appropriate specialization of
coefficients.
\end{abstract}

\date{March 15, 2006; revised August 10, 2006.}

 \thanks{Research of S.~F.\ supported 
by the NSF grant DMS-0245385.
Research of A.~Z. supported by the NSF grants DMS-0200299 and
DMS-0500534, and by a Humboldt Research Award.}

\maketitle


\tableofcontents


\section{Introduction}
\label{sec:intro}

Since their introduction in \cite{ca1}, cluster algebras have
found applications in a diverse variety of settings which include
(in no particular order) total positivity, Lie theory, quiver
representations, Teichm\"uller theory, Poisson geometry, discrete
dynamical systems, tropical geometry, and algebraic combinatorics.
See, e.g., \cite{caldero-keller, fg-survey, fomin-reading, cdm}
and references therein.
This paper, the fourth in a series,
continues the study of the structural theory of cluster
algebras undertaken in its prequels \cite{ca1, ca2, ca3}.
As in those earlier papers, we try to keep the exposition reasonably
self-contained.

\medskip

A \emph{cluster algebra}~$\Acal$ of rank~$n$ is a 
subalgebra of an \emph{ambient field}~$\Fcal$ isomorphic to a
field of rational functions in $n$ variables.
Each cluster algebra comes equipped with a distinguished set
of generators called \emph{cluster variables}; this set is a union of
overlapping algebraically independent $n$-subsets of~$\Fcal$
called \emph{clusters}.
The clusters are related to each other by birational transformations
of the following kind:
for every cluster $\xx$ and every cluster variable $x \in \xx$,
there is another cluster $\xx' = \xx - \{x\} \cup \{x'\}$, with
the new cluster variable~$x'$ determined by an
\emph{exchange relation} of the form
\begin{equation}
\label{eq:exchange-not-detailed}
x\,x' = p^+ M^+ + p^- M^-  \,.
\end{equation}
Here $p^+$ and $p^-$ belong to a \emph{coefficient semifield}~$\PP$,
while $M^+$ and $M^-$ are two
monomials in the elements of $\xx - \{x\}$.
(See Definition~\ref{def:seed-mutation} for precise details.)
Each exchange relation involves two different kinds of data:
the \emph{exchange matrix}~$B$ encoding the exponents in $M^+$ and $M^-$,
and the two \emph{coefficients} $p^+$ and $p^-$.
In the previous papers of this series, a lot of attention
was given to exchange matrices and their dynamics.
The current paper brings into focus the dynamics of
coefficients.
Before describing the results in detail, we would like to
offer some
justification for the importance and timeliness of this study.

\medskip

$\bullet$
The coefficient dynamics is of interest in its own right.
It can be viewed as a far-reaching generalization of
Al.~Zamolodchikov's \emph{$Y$-systems}, introduced and studied
in~\cite{zamolodchikov} in the context of thermodynamic Bethe Ansatz.
Algebraic properties of $Y$-systems were investigated in~\cite{yga};  
in particular, we established their periodicity,
confirming a conjecture by Zamolodchikov (the type~$A$ case was
settled in~\cite{frenkel-szenes, gliozzi-tateo}).
Although our motivation for studying $Y$-systems
came from the coefficient dynamics in \cite[(5.4)--(5.5)]{ca1},
cluster algebras played no explicit role in \cite{yga}.
In the current paper, the study of $Y$-systems in the context of
cluster algebras allows us to extend and sharpen the results
of~\cite{yga}.

\medskip

$\bullet$
Systematically exploring the interplay between two types of
dynamics---that of cluster variables and that of coefficients---leads
to a better understanding of both phenomena.
The constructions we use to express this interplay are close in spirit
(albeit different in technical details) to the
notion of \emph{cluster ensemble} introduced and studied by
V.~Fock and A.~Goncharov~\cite{fg2} as a tool in higher
Teichm\"uller theory.
Our coefficient-based approach uncovers an unexpected ``common
source" of the two types of dynamics, expressing both the cluster variables and
the coefficients in terms of a new family of \emph{$F$-polynomials},
which generalize the Fibonacci polynomials of~\cite{yga}.
This approach also yields a new constructive way to
express the (conjectural) ``Langlands duality" between the two kinds
of dynamics, suggested in~\cite{fg2}.

\medskip

$\bullet$
Last but not least,
all examples of cluster algebras coming from
geometry of semisimple groups (see~\cite{ca2,ca3,cdm}) have
nontrivial coefficients.
(By contrast, most of
the recent developments in cluster algebra theory obtained in the
context of cluster categories and cluster tilted algebras deal
exclusively with the case of trivial coefficients.)
There are many examples of cluster algebras of geometric origin
which have the same combinatorics of exchange matrices but totally
different systems of coefficients.
This motivates the study of the dependence of a cluster algebra
structure on the choice of coefficients.
One surprising conclusion that we reach in this study is the
existence of ``universal formulas" for cluster variables;
these formulas reduce the case of arbitrary coefficients to a
very particular special case of \emph{principal coefficients}.

\medskip

The paper is organized as follows.
Requisite background
on cluster algebras
is presented in Section~\ref{sec:cluster-prelims}.
In keeping with the tradition established in~\cite{ca1,ca2,ca3},
we rephrase the definitions to better suit our current purposes.
In particular, we reconcile the setups in~\cite{ca1}
and~\cite{ca2} by emphasizing the difference between the \emph{labeled
seeds} attached to the vertices of an $n$-regular tree
(see Definitions~\ref{def:seed} and~\ref{def2.3:Tn}), and
the (unlabeled) seeds attached to the vertices of an \emph{exchange
graph} (see
Definitions \ref{def:seed-equivalence}--\ref{def:exchange-graph}).
(Informally,
a \emph{seed} is a cluster~$\xx$
together with the corresponding exchange matrix~$B$ and the collection of
coefficients~$p^{\pm}$ appearing in all exchanges
from~$\xx$.)

In Section~\ref{sec:separation}, we introduce some of the main
new concepts that play a central role in the paper,
most notably, the cluster algebras with principal coefficients
(Definition~\ref{def:principal-coeffs})
and, based on the latter, the $F$-polynomials (Definition~\ref{def:Aprin}).
Our first main result, Theorem~\ref{th:reduction-principal}
(sharpened in Corollary~\ref{cor:xjt=F/F}),
expresses any cluster variable in any cluster algebra in terms of the
(arbitrary) initial seed using the appropriate $F$-polynomials.
Thus, one can think of principal coefficients as a crucial special
case providing control over cluster algebras with arbitrary
coefficients.

The universal formulas of
Theorem~\ref{th:reduction-principal}/Corollary~\ref{cor:xjt=F/F}
have an \emph{a priori} unexpected
``separation of additions'' property:
each cluster variable is written as a ratio of two polynomial
expressions in the initial data, one of the them (the numerator)
employing only the ``ordinary addition'' in~$\Fcal$,
and another one (the denominator) involving only the ``auxiliary addition''
in the coefficient semifield~$\PP$.

These formulas lead to the following structural result obtained in
Section~\ref{sec:exchange-graphs}
(Theorem~\ref{th:exchange-graph-principal}): the exchange graph of
an algebra with principal coefficients serves as a cover for the
exchange graph of any cluster algebra with the same exchange matrix.
This is a step towards the much stronger
Conjecture~\ref{con:exchange-graph-independence}:
the exchange graph of a cluster algebra depends only on the exchange matrix at
any seed---but not on the choice of coefficients.

Properties of $F$-polynomials
are discussed in Section~\ref{sec:F-polynomials}.
Although the definition of $F$-polynomials is simple and elementary,
proving some of their properties resisted our efforts.
One example is Conjecture~\ref{con:conjectures-on-denoms},
an innocent-looking assertion that every
$F$-polynomial has constant term~$1$.

In Section~\ref{sec:principal-gradings}, we show that a cluster
algebra with principal coefficients possesses
a $\ZZ^n$-grading with respect to which all cluster
variables are homogeneous elements.
We then discuss general properties (partly proven and partly
conjectural) of the (multi-)degrees of cluster variables, which we
refer to as \emph{$\gg$-vectors}.
In particular, we show that these vectors can be expressed in
terms of (tropical evaluations of) $F$-polynomials. 

Section~\ref{sec:cluster-monomials} applies the
tools developed in earlier sections to the study of \emph{cluster
monomials}, that is, monomials in cluster variables all of which
belong to the same cluster.
The significance of cluster monomials stems from the expectation
that they belong to the suitably defined ``canonical basis'' of a
cluster algebra.
There is still a long way towards this result: even
the linear independence of cluster monomials remains open 
(see Conjecture~\ref{con:cl-mon-lin-indep}, reproducing
\cite[Conjecture~4.16]{cdm}).
As a possible tool in proving this conjecture,
we discuss two families of combinatorial
parametrizations of cluster monomials by integer vectors
in~$\ZZ^n$; they are somewhat reminiscent of
two families of parametrizations of canonical bases
(Lusztig's and string parametrizations) in the theory of quantum
groups.
One of these parametrizations of cluster monomials (by \emph{denominator vectors})
has already appeared in~\cite{ca2, cdm};
another one (by $\gg$-vectors) is new.
Note that in Section~\ref{sec:principal-gradings},
$\gg$-vectors are defined in restricted generality,
for the case of principal coefficients only; in
Section~\ref{sec:cluster-monomials},
we find a way to extend this definition to a much more general
class of coefficients.

In the rest of the paper
(Sections~\ref{sec:bipartite}--\ref{sec:universal-coeffs}), we
restrict our attention to some special classes of cluster algebras.
Sections~\ref{sec:bipartite}--\ref{sec:denom-bipartite} deal with the
\emph{bipartite cluster algebras}, a class
that includes all cluster algebras of finite type.
A bipartite cluster algebra has a seed $\Sigma$ such that, in
each of the $n$ exchange relations \eqref{eq:exchange-not-detailed}
emanating from~$\Sigma$,
one of the monomials $M^+$ and $M^-$ is equal to~$1$.
Such a seed $\Sigma$ is also called \emph{bipartite},
and is naturally included into a sequence of
bipartite seeds called the \emph{bipartite belt}.
In Section~\ref{sec:bipartite}, we show that the restriction of the
coefficient dynamics to the bipartite belt yields a
natural generalization of Zamolodchikov's $Y$-systems to the case
of arbitrary symmetrizable generalized Cartan matrices.
Two main results of this section
sharpen their counterparts obtained in~\cite{yga}:
Theorem~\ref{th:Laurent-Y-system}
(the ``Laurent phenomenon'')
asserts that the solutions of a
generalized $Y$-system can be expressed as Laurent polynomials in the initial
data, while Theorem~\ref{th:belt-periodic} improves upon the periodicity theorem
of~\cite{yga} by showing the necessity of the finite type assumption.

The proof of Theorem~\ref{th:belt-periodic} given in
Section~\ref{sec:denom-bipartite}
relies on the calculation of denominator vectors and $\gg$-vectors for
the cluster variables on the bipartite belt.
We also need some properties of the root system associated with a
symmetrizable generalized Cartan matrix.
These properties are discussed in Section~\ref{sec:prelim-roots},
which can be read independently of the rest of the paper.
We believe these results 
to be of independent interest.
In particular, Theorem~\ref{th:positive-roots-2} and
Corollary~\ref{cor:coxeter-elt-inf-type} 
strengthen the following result obtained in~\cite{blm}:
in the Weyl group associated to an indecomposable symmetrizable Cartan
matrix of infinite type,
the ``bipartite" Coxeter element has infinite order.

The last two sections \ref{sec:finite-type} and \ref{sec:universal-coeffs}
deal with the cluster algebras of finite type, i.e., those with
finitely many seeds.
The classification of these algebras obtained in \cite{ca2}
is parallel to the Cartan-Killing classification
of semisimple Lie algebras and finite root systems.
In Section~\ref{sec:finite-type}, we demonstrate that several conjectures
made in the preceding sections hold for the algebras of finite type.
In particular, we prove linear independence of cluster monomials
(Conjecture~\ref{con:cl-mon-lin-indep}) for the cluster algebras of finite type;
see Theorem~\ref{th:cluster-monomials-linear-indep}.
Other conjectures are proved by combining our results on bipartite
cluster algebras
with the structural results in~\cite{ca2}, such as the following
property (implicit in~\cite{ca2}):
in a cluster algebra of finite type, every cluster
variable appears in a seed lying on the bipartite belt.

In Theorem~\ref{th:fib=Fib}, we show that in finite type,
the $F$-polynomials coincide, up to a simple twist,
with the \emph{Fibonacci polynomials} of~\cite{yga}.

The concluding Section~\ref{sec:universal-coeffs} presents an
explicit construction of the ``universal" system of coefficients for
cluster algebras of each finite type.
This universality is understood in the following sense.
As a consequence of results in \cite{ca2}, for any two cluster
algebras $\Acal$ and $\bar \Acal$ of the same (finite)
Cartan-Killing type, there is a natural bijection between the sets
of cluster variables for $\Acal$ and $\bar \Acal$.
So the only feature that distinguishes these algebras from each
other is a choice of coefficients
(in turn predicated on the choice of the ground ring, the integral group
ring~$\ZZ \PP$ of the multiplicative group of the coefficient
semifield~$\PP$).
The main result of Section~\ref{sec:universal-coeffs}
(Theorem~\ref{th:universal}) asserts that, among all cluster algebras
of a given Cartan-Killing type, there is a \emph{universal} cluster algebra
$\Acal^{\rm univ}$ such that every other algebra of this type is
obtained from $\Acal^{\rm univ}$ by a (unique) ``change of base''
resulting from a homomorphism of multiplicative groups
of coefficient semifields.
The existence of $\Acal^{\rm univ}$ is by no means obvious, and it
is not clear how far this result can be extended beyond finite
type.

\medskip

The reader will notice a large number of conjectures put forth
in Sections~\ref{sec:exchange-graphs}--\ref{sec:cluster-monomials}.
We conclude the paper by an index of partial results that we obtained towards
these conjectures.
We note that a promising tool for attacking some of these conjectures
is provided by a geometric interpretation (found in~\cite{caldchap} and generalized
in~\cite{caldkell2}) of Laurent expansions of cluster variables
in terms of Grassmannians of quiver representations.


\section{Preliminaries on cluster algebras}
\label{sec:cluster-prelims}


To state our results, we 
recall the basic setup
introduced in~\cite{ca1,ca2,ca3}.
In keeping with the tradition established in the previous papers in
the series, we rephrase the definitions to better suit the current
purposes.

The definition of a cluster algebra~$\Acal$ starts with introducing
its ground ring.
Let $(\PP,\oplus, \cdot)$ be a \emph{semifield}, i.e.,
an abelian multiplicative group endowed with a binary operation of
\emph{(auxiliary) addition}~$\oplus$ which is commutative, associative, and
distributive with respect to the multiplication in~$\PP$.
The multiplicative group of~$\PP$ is torsion-free
\cite[Section~5]{ca1}, hence its group ring~$\ZZ\PP$---which will be
used as a \emph{ground ring} for~$\Acal$---is a domain.


The following two examples of semifields 
will be of particular importance to us.

\begin{definition}
[\emph{Universal semifield}]
\label{def:semifield-universal}
Let $\Qsf(u_1, \dots, u_\ell)$ denote the set of all rational
functions in $\ell$ independent variables $u_1, \dots, u_\ell$ which can be written as
\emph{subtrac\-tion-free} rational expressions in $u_1, \dots, u_\ell$.
For example, $u^2-u+1=\frac{u^3+1}{u+1}\in\Qsf(u)$.
The set $\Qsf(u_1, \dots, u_\ell)$ is a semifield with respect to
the usual operations of multiplication and addition.
This example is universal: any subtraction-free identity that
holds in $\Qsf(u_1, \dots, u_\ell)$ remains valid for any elements
$u_1,\dots,u_\ell$ in an arbitrary semifield
(see \cite[Lemma~2.1.6]{bfz-tp}).
\end{definition}

\begin{definition}
[\emph{Tropical semifield}]
\label{def:semifield-tropical}
Let $J$ be a finite set of labels,
and let $\Trop (u_j: j \in J)$ be an abelian group (written
multiplicatively) freely generated by the elements $u_j \, (j \in J)$.
We define the addition~$\oplus$ in $\Trop (u_j: j \in J)$ by
\begin{equation}
\label{eq:tropical-addition}
\prod_j u_j^{a_j} \oplus \prod_j u_j^{b_j} =
\prod_j u_j^{\min (a_j, b_j)} \,,
\end{equation}
and call $(\Trop (u_j: j \in J),\oplus,\cdot)$ a \emph{tropical
 semifield}.
To illustrate, $u_2\oplus u_1^2 u_2^{-1}=u_2^{-1}$ in $\Trop
 (u_1,u_2)$.
If $J$ is empty,
we obtain the trivial semifield consisting of a single element~$1$.

The group ring of $\Trop (u_j: j \in J)$ is the ring of Laurent
polynomials in the variables~$u_j\,$.
\end{definition}

The terminology in Definition~\ref{def:semifield-tropical}
is consistent with~\cite{ca1}.
It is different from (although closely related to)
the terminology in \cite[Section~2.1]{bfz-tp}, nowadays commonly used in
Tropical Geometry.

As an \emph{ambient field} for a cluster algebra~$\Acal$, we take a field $\Fcal$
isomorphic to the field of rational functions in $n$ independent
variables (here $n$ is the \emph{rank} of~$\Acal$),
with coefficients in~$\QQ \PP$.
Note that the definition of $\Fcal$ ignores the auxiliary addition
in~$\PP$. 

\begin{definition}
[\emph{Labeled seeds}]
\label{def:seed}
A (skew-symmetrizable) \emph{labeled $Y$-seed} in~$\PP$
is a pair $(\yy, B)$, where
\begin{itemize}
\item
$\yy = (y_1, \dots, y_n)$ is an $n$-tuple
of elements of $\PP$, and
\item
$B = (b_{ij})$ is an $n\!\times\! n$ integer matrix
which is \emph{skew-symmetrizable}.
\end{itemize}
That~is,
$d_i b_{ij} = - d_j b_{ji}$ for some positive integers $d_1, \dots, d_n$.
A \emph{labeled seed} in~$\Fcal$ is
a triple $(\xx, \yy, B)$, where
\begin{itemize}
\item
$(\yy, B)$ is a labeled $Y$-seed, and
\item
$\xx = (x_1, \dots, x_n)$ is an $n$-tuple of elements of~$\Fcal$
forming a \emph{free generating set}. 
\end{itemize}
\pagebreak[3]
That is, $x_1, \dots, x_n$
are algebraically independent over~$\QQ \PP$, and $\Fcal = \QQ \PP(x_1, \dots, x_n)$.
We refer to~$\xx$ as the (labeled)
\emph{cluster} of a labeled seed $(\xx, \yy, B)$,
to the tuple~$\yy$ as the \emph{coefficient tuple}, and to the
matrix~$B$ as the \emph{exchange matrix}.
\end{definition}

The (unlabeled) seeds as defined in \cite{ca2} and in
Definition~\ref{def:seed-equivalence} below, are obtained by
identifying labeled seeds that differ from each other by simultaneous permutations of
the components in $\xx$ and~$\yy$, and of the rows and columns of~$B$.
In this paper, we will mostly deal with labeled seeds;
meanwhile, we will sometimes refer to labeled seeds simply as seeds, when
there is no risk of confusion.

Throughout the paper, we use the notation
\begin{align*}
[x]_+ &= \max(x,0); \\
\sgn(x) &=
\begin{cases}
-1 & \text{if $x<0$;}\\
0  & \text{if $x=0$;}\\
 1 & \text{if $x>0$;}
\end{cases}\\
[1,n]&=\{1, \dots, n\}\,.
\end{align*}

\begin{definition}
[\emph{Seed mutations}]
\label{def:seed-mutation}
Let $(\xx, \yy, B)$ be a labeled seed in $\Fcal$, as in
Definition~\ref{def:seed} above,
and let $k \in [1,n]$.
The \emph{seed mutation} $\mu_k$ in direction~$k$ transforms
$(\xx, \yy, B)$ into the labeled seed
$\mu_k(\xx, \yy, B)=(\xx', \yy', B')$ defined as follows:
\begin{itemize}
\item
The entries of $B'=(b'_{ij})$ are given by
\begin{equation}
\label{eq:matrix-mutation}
b'_{ij} =
\begin{cases}
-b_{ij} & \text{if $i=k$ or $j=k$;} \\[.05in]
b_{ij} + \sgn(b_{ik}) \ [b_{ik}b_{kj}]_+
 & \text{otherwise.}
\end{cases}
\end{equation}
\item
The coefficient tuple $\yy'=(y_1',\dots,y_n')$ is given by
\begin{equation}
\label{eq:y-mutation}
y'_j =
\begin{cases}
y_k^{-1} & \text{if $j = k$};\\[.05in]
y_j y_k^{[b_{kj}]_+}
(y_k \oplus 1)^{- b_{kj}} & \text{if $j \neq k$}.
\end{cases}
\end{equation}
\item
The cluster $\xx'=(x_1',\dots,x_n')$ is given by
$x_j'=x_j$ for $j\neq k$,
whereas $x'_k \in \Fcal$ is determined
by the \emph{exchange relation}
\begin{equation}
\label{eq:exchange-rel-xx}
x'_k = \frac
{y_k \ \prod x_i^{[b_{ik}]_+}
+ \ \prod x_i^{[-b_{ik}]_+}}{(y_k \oplus 1) x_k} \, .
\end{equation}
\end{itemize}
We also say that the \emph{$Y$-seed mutation} in direction~$k$ transforms
$(\yy, B)$ into the labeled $Y$-seed $\mu_k(\yy, B)=(\yy', B')$ given by
\eqref{eq:matrix-mutation} and~\eqref{eq:y-mutation}.

It is easy to see that $B'$ is
skew-symmetrizable (with the same choice of $d_1, \dots, d_n$),
implying that $(\xx', \yy', B')$ is indeed a labeled seed.
Furthermore, the seed mutation $\mu_k$ is involutive,
that is, it transforms $(\xx', \yy', B')$ into the original labeled seed
$(\xx, \yy, B)$.
\end{definition}

A few comments are in order.


\begin{remark}
The transformation $\mu_k:B\mapsto B'$ defined by
\eqref{eq:matrix-mutation} is easily checked to coincide with
the \emph{matrix mutation} defined in~\cite[(4.3)]{ca1}).
Another equivalent way to define it is
\begin{equation}
\label{eq:alt-matrix-mutation}
b'_{ij} =
\begin{cases}
-b_{ij} & \text{if $i=k$ or $j=k$;} \\[.05in]
b_{ij} + [-b_{ik}]_+\,b_{kj} +b_{ik} [b_{kj}]_+
 & \text{otherwise.}
\end{cases}
\end{equation}
\end{remark}

\begin{remark}
The transformation $\yy\mapsto\yy'$
given by \eqref{eq:y-mutation} first appeared in
\cite[(5.5)]{ca1}, and then in a different context
in \cite[Lemma~1.3]{gsv1} and \cite[(6)]{fg2}
(cf.\ Proposition~\ref{pr:another-y-pattern} below).
\end{remark}

\begin{remark}
Originally \cite[Definition~5.3]{ca1} \cite[Section~1.2]{ca2},
we defined the coefficient tuple as a $2n$-tuple
$\pp = (p_1^\pm, \dots, p_n^\pm)$ of elements of $\PP$ satisfying the
\emph{normalization condition} $p_j^+ \oplus p_j^- = 1$ for all~$j$.
The equivalence of that setup with the current one was established
in \cite[(5.2)-(5.3)]{ca1}: by setting
\begin{equation}
\label{eq:y-through-p}
y_j = \frac{p_j^+}{p_j^-},
\end{equation}
the coefficients $p_j^\pm$ are recovered via
\begin{equation}
\label{eq:p-through-y}
p_j^+ = \frac{y_j}{y_j \oplus 1}, \quad
p_j^- = \frac{1}{y_j \oplus 1} \,,
\end{equation}
and the exchange relation \eqref{eq:exchange-rel-xx} takes
the usual form (cf.\ \cite[Definition~1.1]{ca2})
\begin{equation}
\label{eq:exch-rel-usual}
x'_k = x_k^{-1}\left(
p_k^+ \ \displaystyle\prod_{i=1}^n x_i^{[b_{ik}]_+}
+ p_k^- \ \displaystyle\prod_{i=1}^n x_i^{[-b_{ik}]_+}\right) \,.
\end{equation}
\end{remark}


\begin{definition}
[\emph{Regular $n$-ary tree}]
\label{def2.3:Tn}
As in \cite{ca1}, we consider the \emph{$n$-regular tree}~$\TT_n$
whose edges are labeled by the numbers $1, \dots, n$,
so that the $n$ edges emanating from each vertex receive
different labels.
We write $t \overunder{k}{} t'$ to indicate that vertices
$t,t'\in\TT_n$ are joined by an edge labeled by~$k$.
\end{definition}

\begin{definition}
[\emph{Patterns}]
\label{def:patterns}
A \emph{cluster pattern} (resp., \emph{$Y$-pattern}) is an assignment
of a labeled seed $\Sigma_t=(\xx_t, \yy_t, B_t)$
(resp., a labeled $Y$-seed $(\yy_t, B_t)$)
to every vertex $t \in \TT_n$, such that the seeds assigned to the
endpoints of any edge $t \overunder{k}{} t'$ are obtained from each
other by the seed mutation in direction~$k$.
The elements of $\Sigma_t$ are written as follows:
\begin{equation}
\label{eq:seed-labeling}
\xx_t = (x_{1;t}\,,\dots,x_{n;t})\,,\quad
\yy_t = (y_{1;t}\,,\dots,y_{n;t})\,,\quad
B_t = (b^t_{ij})\,.
\end{equation}
\end{definition}

Clearly, a cluster pattern (resp., $Y$-pattern) is uniquely determined
by each of its seeds 
(resp., $Y$-seeds), 
which can be chosen arbitrarily.

As will be explained in Section~\ref{sec:bipartite},
the notion of a $Y$-pattern is a far-reaching generalization of
Zamolodchikov's $Y$-systems~\cite{zamolodchikov, yga}.

\begin{example}
[\emph{Type~$A_2$; cf.\ \cite[Section~6]{ca1}}]
\label{example:A2-pattern}
Let $n=2$. Then the tree $\TT_2$ is an infinite chain.
We denote its vertices by
$\dots,t_{-1}\,,t_0\,,t_1\,,t_2\,,\dots$,
and label its edges as follows:
\begin{equation}
\label{eq:TT2}
\cdots
\overunder{2}{} t_{-1}
\overunder{1}{} t_0
\overunder{2}{} t_1
\overunder{1}{} t_2
\overunder{2}{} t_3
\overunder{1}{} \cdots \,.
\end{equation}
We denote the corresponding seeds by
$\Sigma_m=\Sigma_{t_m}=(\xx_m,\yy_m,B_m)$, for $m\in\ZZ$.
Let the initial seed $\Sigma_0$ be
\begin{equation}
\label{eq:A2-pattern-initial}
\xx_0=(x_1,x_2), \quad
\yy_0=(y_1,y_2), \quad
B_0=\begin{bmatrix}
 0 & 1\\
-1 & 0
\end{bmatrix}\,.
\end{equation}
We then recursively compute the seeds $\Sigma_1,\dots,\Sigma_5$ as
shown in Table~\ref{table:seeds-A2}.
Note that the labeled seed $\Sigma_5$ is obtained from $\Sigma_0$
by interchanging the indices~$1$ and $2$; the sequence then continues
by obvious periodicity (so that~$\Sigma_{10}$ becomes identical to~$\Sigma_0$,~etc.)
\end{example}

\begin{table}[ht]
\begin{equation*}
\begin{array}{|c|c|cc|cc|}
\hline
&&&&&\\[-4mm]
t & B_t & & \yy_t\hspace{20mm} & \hspace{20mm} \xx_t&  \\[1mm]
\hline
&&&&&\\[-3mm]
0 & \begin{bmatrix}0&1\\-1&0\end{bmatrix} &
y_1&y_2&x_1& x_2
\\[4.5mm]
\hline
&&&&&\\[-3mm]
1 & \begin{bmatrix}0&-1\\1&0\end{bmatrix} &
y_1(y_2\oplus 1) &
\dfrac{1}{y_2}
& x_1
& \dfrac{x_1y_2+1}{x_2(y_2\oplus 1)}
\\[4.5mm]
\hline
&&&&&\\[-3mm]
2 & \begin{bmatrix}0&1\\-1&0\end{bmatrix} &
\dfrac{1}{y_1(y_2\oplus 1)}&
\dfrac{y_1y_2\oplus y_1\oplus 1}{y_2}&
\dfrac{x_1y_1y_2+y_1+x_2}{(y_1y_2\oplus y_1\oplus 1)x_1x_2} &
\dfrac{x_1y_2+1}{x_2(y_2\oplus 1)}
\\[4.5mm]
\hline
&&&&&\\[-3mm]
3 & \begin{bmatrix}0&-1\\1&0\end{bmatrix} &
\dfrac{y_1\oplus 1}{y_1 y_2}&
\dfrac{y_2}{y_1y_2\oplus y_1\oplus 1}&
\dfrac{x_1y_1y_2+y_1+x_2}{(y_1y_2\oplus y_1\oplus 1)x_1x_2} &
\dfrac{y_1+x_2}{x_1(y_1\oplus 1)}
\\[4.5mm]
\hline
&&&&&\\[-3mm]
4 & \begin{bmatrix}0&1\\-1&0\end{bmatrix} &
\dfrac{y_1y_2}{y_1\oplus 1} &
\dfrac{1}{y_1}&
x_2 &
\dfrac{y_1+x_2}{x_1(y_1\oplus 1)}
\\[4.5mm]
\hline
&&&&&\\[-3mm]
5 & \begin{bmatrix}0&-1\\1&0\end{bmatrix} & y_2 & y_1 & x_2&x_1
\\[4.5mm]
\hline
\end{array}
\end{equation*}
\caption{Seeds in type~$A_2$}
\label{table:seeds-A2}
\end{table}

\pagebreak[3]

Now everything is in place for defining cluster algebras.

\begin{definition}
[\emph{Cluster algebra}]
\label{def:cluster-algebra}
Given a cluster pattern, we denote
\begin{equation}
\label{eq:cluster-variables}
\Xcal
= \bigcup_{t \in \TT_n} \xx_t
= \{ x_{i,t}\,:\, t \in \TT_n\,,\ 1\leq i\leq n \} \ ,
\end{equation}
the union of clusters of all the seeds in the pattern.
We refer to the elements $x_{i,t}\in \Xcal$ as \emph{cluster variables}.
The 
\emph{cluster algebra} $\Acal$ associated with a
given cluster pattern is the $\ZZ \PP$-subalgebra of the ambient field $\Fcal$
generated by all cluster variables: $\Acal = \ZZ \PP[\Xcal]$.
We denote $\Acal = \Acal(\xx, \yy, B)$, where
$(\xx,\yy,B)=(\xx_t,\yy_t,B_t)$
is any labeled seed in the underlying cluster pattern.
\end{definition}


\begin{definition}
[\emph{Geometric type}]
\label{def:geom-type}
A cluster algebra (or a $Y$-pattern, or a cluster pattern)
is of \emph{geometric type} if the coefficient
semifield $\PP$ is a tropical semifield of
Definition~\ref{def:semifield-tropical}.
\end{definition}

For patterns of geometric type, it is convenient to denote the generators of $\PP$ by
$x_{n+1},\dots,x_m$ (for some integer $m\geq n$), so that
$\PP=\Trop(x_{n+1},\dots,x_m)$.
As the coefficients $y_{1;t},\dots,y_{n;t}$ at each seed
$\Sigma_t=(\xx_t,\yy_t,B_t)$
are Laurent monomials in
$x_{n+1},\dots,x_m$, we define the integers $b^t_{ij}$,
for $j\in [1,n]$ and $n<i\leq m$, by
\begin{equation}
\label{eq:yt-geometric}
y_{j;t}=\prod_{i=n+1}^m x_i^{b^t_{ij}}\,.
\end{equation}
That is, we include the exchange matrix $B_t$ as a submatrix into a
larger $m\times n$ matrix
\begin{equation}
\label{eq:tilde-B}
\tilde B_t=(b^t_{ij}) \qquad (1\leq i\leq m, \ \ 1\leq j\leq n)
\end{equation}
whose matrix elements $b^t_{ij}$ with $i>n$ are used to encode the
coefficients~$y_j = y_{j;t}\,$.
It then follows from Definition~\ref{def:semifield-tropical} that the
formulas \eqref{eq:p-through-y} for a seed $\Sigma_t$ become
\[
p_j^+ = \prod_{i=n+1}^m x_i^{[b^t_{ij}]_+}\,, \quad
p_j^- = \prod_{i=n+1}^m x_i^{[-b^t_{ij}]_+}\,,
\]
and the exchange relations \eqref{eq:exch-rel-usual} become
\begin{equation}
\label{eq:exch-rel-geom}
x'_k = x_k^{-1}\left(
\displaystyle\prod_{i=1}^m x_i^{[b^t_{ik}]_+}
+ \displaystyle\prod_{i=1}^m x_i^{[-b^t_{ik}]_+}\right) ,
\end{equation}
where we denote $x_i=x_{i;t}$ for $i\in [1,n]$.
Furthermore, the $Y$-seed mutation formulas \eqref{eq:y-mutation}
take the form \eqref{eq:matrix-mutation}, with $i>n$.
(Cf.\ \cite[Proposition~5.8]{ca1}.)
Thus, extending the notion of matrix mutation to rectangular matrices,
we have $\tilde B_{t'} = \mu_k(\tilde B_t)$ whenever
$t \overunder{k}{} t'$.
Summing up, the entire $Y$-pattern (hence, up to isomorphism, the
corresponding cluster algebra of geometric type)
is uniquely determined by a single rectangular matrix~$\tilde B_t$.

\section{Separation formulas}
\label{sec:separation}

We now introduce a special class of coefficient
patterns that will play a central role in this paper.

\begin{definition}
[\emph{Principal coefficients}]
\label{def:principal-coeffs}
We say that a $Y$-pattern $t \mapsto (\yy_t,B_t)$ on $\TT_n$
(or a cluster pattern $t \mapsto (\xx_t, \yy_t, B_t)$,
or the corresponding cluster algebra~$\Acal$)  has
\emph{principal coefficients at a vertex~$t_0$} if
$\PP= \Trop(y_1, \dots, y_n)$ and $\yy_{t_0}= (y_1, \dots, y_n)$. \linebreak[3]
In this case, we denote $\Acal=\Aprin(B_{t_0})$.
\end{definition}

\begin{remark}
\label{rem:principal-tildeB}
Definition~\ref{def:principal-coeffs} can be rephrased
as follows: a cluster algebra~$\Acal$ has principal coefficients at a
vertex~$t_0$ if $\Acal$ is of geometric type,
and is associated with the matrix $\tilde B_{t_0}$ of order $2n \times n$
whose principal (i.e., top $n\times n$) part is $B_{t_0}$,
and whose complementary (i.e., bottom) part is
the $n \times n$ identity matrix (cf.\ \cite[Corollary~5.9]{ca1}).
\end{remark}



In Theorem~\ref{th:reduction-principal} below we give
a formula that expresses a cluster variable in an arbitrary
cluster algebra with a given initial exchange matrix $B^0 =B_{t_0}$
in terms of similar formulas for the corresponding cluster algebra
with principal coefficients.
In order to write these formulas,
we need to introduce some notation.

\begin{definition}
[\emph{The functions~$X_{\ell,t}$ and~$F_{\ell,t}$}]
\label{def:Aprin}
Let~$\Acal=\Aprin(B^0)$ be the cluster \linebreak[3]
algebra with principal coefficients at a
vertex~$t_0$, defined by the initial seed
$\Sigma_{t_0}=(\xx_{t_0}\,,\yy_{t_0}\,,B_{t_0})$ with
\begin{equation}
\label{eq:initial-seed}
\xx_{t_0} = (x_1, \dots, x_n), \quad
\yy_{t_0} = (y_1, \dots, y_n), \quad
B_{t_0} = B^0 = (b^0_{ij})\,.
\end{equation}
Thus $\PP=\Trop(y_1,\dots,y_n)$, and all coefficients in all exchange
relations  \eqref{eq:exchange-rel-xx} are monomials in
$y_1,\dots,y_n\,$.
By iterating these exchange relations, we
can express every cluster variable $x_{\ell;t}$ as a (unique)
rational function in $x_1, \dots, x_n, y_1, \dots, y_n$
given by a subtraction-free rational expression; we
denote this rational function by
\begin{equation}
\label{eq:X-sf}
X_{\ell;t} 
= X_{\ell;t}^{B^0;t_0}\in
\Qsf(x_1, \dots, x_n; y_1, \dots, y_n).
\end{equation}
Let $F_{\ell;t} =
F_{\ell;t}^{B^0;t_0}\in
\Qsf(y_1, \dots, y_n)$
denote the rational function obtained from $X_{\ell;t}$ by specializing
all the~$x_j$ to~$1$:
\begin{equation}
\label{eq:F-def}
F_{\ell;t}(y_1, \dots, y_n) = X_{\ell;t}(1, \dots, 1; y_1, \dots, y_n).
\end{equation}
\end{definition}

For instance, we have
$X_{\ell;t_0} = x_\ell$ and $F_{\ell;t_0} = 1$
for all~$\ell$; and if $t_1 \overunder{k}{} t_0\,$, then
\[
X_{k;t_1}^{B^0;t_0} = \frac
{y_k \ \prod x_i^{[b^0_{ik}]_+}
+ \ \prod x_i^{[-b^0_{ik}]_+}}{x_k}, \quad
F_{k;t_1}^{B^0;t_0} = y_k + 1\,.
\]

\begin{example}
[\emph{Principal coefficients, type~$A_2$}]
\label{example:A2-pattern-principal}
Let $n=2$ and
$B^0=\left[
\begin{smallmatrix}
0 & 1\\
-1&0
\end{smallmatrix}
\right]$.
That is, we consider the special case of Example~\ref{example:A2-pattern}
in which the coefficient semifield is $\PP=\Trop(y_1,y_2)$.
In this case, the formulas in Table~\ref{table:seeds-A2} simplify
considerably. The results are shown in the first four columns of
Table~\ref{table:seeds-geom-A2}, where we use the notation
$\xx_t=(X_{1;t}\,,X_{2;t})$, in agreement with
Definition~\ref{def:Aprin}.
Setting $x_1=x_2=1$ yields the polynomials~$F_{i;t}\,$,
shown in the fifth column of the table.

\begin{table}[ht]
\begin{equation*}
\begin{array}{|c|c|cc|cc|cc|}
\hline
&&&&&&&\\[-4mm]
t & \tilde B_t & \hspace{7mm}\yy_t&  & X_{1;t} & X_{2;t}&
F_{1;t} & F_{2;t}
 \\[1mm]
\hline
&&&&&&&\\[-3mm]
0 &
\left[\begin{smallmatrix}0&1\\-1&0\\1&0\\0&1\end{smallmatrix}\right] &
y_1&y_2&
x_1&x_2 &
1 & 1
\\[4.5mm]
\hline
&&&&&&&\\[-3mm]
1 & \left[\begin{smallmatrix}0&-1\\1&0\\1&0\\0&-1\end{smallmatrix}\right] &
y_1 &\dfrac{1}{y_2}&
x_1 & \dfrac{x_1y_2+1}{x_2}&
1& y_2+1
\\[4.5mm]
\hline
&&&&&&&\\[-3mm]
2 & \left[\begin{smallmatrix}0&1\\-1&0\\-1&0\\0&-1\end{smallmatrix}\right] &
\dfrac{1}{y_1}&
\dfrac{1}{y_2}&
\dfrac{x_1y_1y_2+y_1+x_2}{x_1x_2}&
\dfrac{x_1y_2+1}{x_2}&
y_1y_2+y_1+1 &
y_2+1
\\[4.5mm]
\hline
&&&&&&&\\[-3mm]
3 & \left[\begin{smallmatrix}0&-1\\1&0\\-1&0\\-1&1\end{smallmatrix}\right] &
\dfrac{1}{y_1 y_2}&
y_2 &
\dfrac{x_1y_1y_2+y_1+x_2}{x_1x_2}&
\dfrac{y_1+x_2}{x_1}&
y_1y_2+y_1+1 &
y_1+1
\\[4.5mm]
\hline
&&&&&&&\\[-3mm]
4 & \left[\begin{smallmatrix}0&1\\-1&0\\1&-1\\1&0\end{smallmatrix}\right] &
y_1y_2 &
\dfrac{1}{y_1}&
x_2 &
\dfrac{y_1+x_2}{x_1}&
1 &
y_1+1
\\[4.5mm]
\hline
&&&&&&&\\[-3mm]
5 &
\left[\begin{smallmatrix}0&-1\\1&0\\0&1\\1&0\end{smallmatrix}\right] &
y_2 & y_1 & x_2&x_1 & 1 & 1
\\[4.5mm]
\hline
\end{array}
\end{equation*}
\smallskip

\caption{Type~$A_2$, principal coefficients}
\label{table:seeds-geom-A2}
\end{table}
\end{example}

\pagebreak[3]

As we shall now demonstrate, each rational function
$X_{\ell;t}$ is in fact a Laurent
polynomial with integer coefficients, while each $F_{\ell;t}$ is
actually a polynomial in~$y_1,\dots,y_n$.

\begin{theorem}[\emph{Laurent phenomenon} {\cite[Theorem~3.1]{ca1}}]
\label{th:Laurent-phenom}
The cluster algebra~$\Acal$ associated
with a seed $(\xx, \yy, B)$ is contained in the
Laurent polynomial ring $\ZZ \PP[\xx^{\pm 1}]$,
i.e., every element of $\Acal$ is a
Laurent polynomial over $\ZZ \PP$ in the cluster
variables from~$\xx$.
\end{theorem}

In view of \cite[Proposition~11.2]{ca2}, for a cluster algebra
with principal coefficients, Theorem~\ref{th:Laurent-phenom} can be
sharpened as follows.

\begin{proposition}
\label{pr:Laurent-sharpened}
Let~$\Acal=\Aprin(B^0)$ be the cluster algebra with principal coefficients
at a vertex $t_0$ with the initial seed given by \eqref{eq:initial-seed}.
Then~$\Acal \subset \ZZ [x_1^{\pm 1}, \dots, x_n^{\pm 1}; y_1, \dots,
  y_n]$.
That is, every element of $\Acal$ is a
Laurent polynomial in $x_1, \dots, x_n$ whose coefficients
are integer polynomials in $y_1, \dots, y_n$.
Thus
\begin{equation}
\label{eq:X-poly}
X_{\ell;t} \in \ZZ [x_1^{\pm 1}, \dots, x_n^{\pm 1}; y_1, \dots, y_n],
\quad F_{\ell;t} \in \ZZ [y_1, \dots, y_n] \ .
\end{equation}
\end{proposition}

The polynomials in \eqref{eq:X-poly} are conjectured to have
positive coefficients.
We note that not until the positivity of coefficients of the (Laurent) polynomial
expansion for some $X_{j;t}$ (or $F_{j;t}$) has been established, can
such an expansion be used for evaluation in an arbitrary semifield.
On the other hand,  one can always use an
appropriate subtraction-free expression.

As explained in Section~\ref{sec:finite-type},
in the special case of cluster algebras of finite type, the
polynomials~$F_{\ell;t}$ are closely related to the ``Fibonacci
polynomials'' of~\cite{yga}.


We will need one more piece of notation in order to state our next
theorem.
If $F$ is a subtraction-free rational expression over $\QQ$
in several variables, $\PP$ a semifield,
and $u_1,\dots,u_\ell$ some elements of~$\PP$,
then we denote by $F|_\PP(u_1,\dots,u_\ell)$ the evaluation of $F$ at
$u_1,\dots,u_\ell\,$.
This is well defined in view of the universality property of the
semifield $\Qsf(u_1, \dots, u_\ell)$ alluded to in
Definition~\ref{def:semifield-universal}.
For example, if
$F(u_1,u_2)= 
u_1^2 - u_1 u_2 + u_2^2
\in\Qsf(u_1,u_2)$,
and $\PP=\Trop(y_1,y_2)$, then
$F|_\PP(y_1,y_2)=\frac{y_1^3\oplus y_2^3}{y_1\oplus y_2}=1$.

\begin{theorem}
\label{th:reduction-principal}
Let $\Acal$ be a cluster algebra over an arbitrary semifield $\PP$,
with a seed at an initial vertex $t_0$ given by \eqref{eq:initial-seed}.
Then the cluster variables in~$\Acal$ can be expressed as follows:
\begin{equation}
\label{eq:xjt-reduction-principal}
x_{\ell;t} = \frac{X_{\ell;t}^{B^0;t_0}|_\Fcal (x_1, \dots, x_n;y_1, \dots, y_n)}
{F_{\ell;t}^{B^0;t_0}|_\PP (y_1, \dots, y_n)} \, .
\end{equation}
\end{theorem}

Formula \eqref{eq:xjt-reduction-principal} exhibits the ``separation
of additions'' phenomenon: the numerator in the right-hand side
is totally independent of the auxiliary addition~$\oplus$,
while the denominator does not involve the (ordinary) addition
in~$\Fcal$.

\begin{example}
To obtain the expressions for the cluster variables in
the last column of Table~\ref{table:seeds-A2},
one needs to replace each $+$ sign in the fifth
column of Table~\ref{table:seeds-geom-A2} by the
$\oplus$ sign
and then
divide each expression in the fourth column of
Table~\ref{table:seeds-geom-A2} by its counterpart in the fifth
column.
\end{example}

The proof of Theorem~\ref{th:reduction-principal},
to be given later in this section, requires some preparation.
%
Our starting point is the following generalization of
\cite[Lemma~1.3]{gsv1} and \hbox{\cite[(6)]{fg2}.}

\pagebreak[3]

\begin{proposition}
\label{pr:another-y-pattern}
Let $t \mapsto (\xx_t, \yy_t, B_t)$ be a cluster pattern in~$\Fcal$,
with
\begin{equation}
\xx_t = (x_{1;t}, \dots, x_{n;t}),\quad
\yy_t = (y_{1;t}, \dots, y_{n;t}),\quad
B_t = (b^t_{ij})
\end{equation}
(as in \eqref{eq:seed-labeling}).
For $t \in \TT_n$, let
$\hat \yy_t = (\hat y_{1;t}, \dots, \hat y_{n;t})$ be the $n$-tuple of
elements in $\Fcal$ given by
\begin{equation}
\label{eq:yhat}
\hat y_{j;t} = y_{j;t} \prod_i x_{i;t}^{b^t_{ij}} \,.
\end{equation}
Then $t \mapsto (\hat \yy_t, B_t)$ is a $Y$-pattern
in $\Fcal$.
In other words,
for $t \overunder{k}{} t'$, we have (cf.~\eqref{eq:y-mutation}):
\begin{equation}
\label{eq:y-mutation-detailed}
\hat y_{j;t'} =
\begin{cases}
\hat y_{k;t}^{-1}  & \text{if $j = k$};\\[.05in]
\hat y_{j;t}  \hat y_{k;t}^{[b^t_{kj}]_+}
(\hat y_{k;t} + 1)^{- b^t_{kj}} & \text{if $j \neq k$.}
\end{cases}
\end{equation}
\end{proposition}

\begin{proof}
We denote
\begin{equation}
\label{eq:xx'yy'bb'}
x_i=x_{i,t}\,,\quad
x'_i=x_{i,t'}\,,\quad
y_i=y_{i,t}\,,\quad
y'_i=y_{i,t'}\,, \quad
b_{ij}=b^t_{ij}\,, \quad
b'_{ij}=b^{t'}_{ij}
\end{equation}
(so as to match the notation of Definition~\ref{def:seed-mutation})
and
\begin{equation}
\label{eq:hat-yy'}
\hat y_j=\hat y_{j,t}\,,\quad
\hat y'_j=\hat y_{j,t'}\,.
\end{equation}
With this notation, we get
\[
\hat y'_k = y'_k \prod_{i \neq k} x_i^{b'_{ik}}
= y_k^{-1} \prod_{i \neq k} x_i^{-b_{ik}} =
\hat y_k^{-1} \,,
\]
proving the first case in \eqref{eq:y-mutation-detailed}.
To check the second case, note that
the exchange relation \eqref{eq:exchange-rel-xx}
can be rewritten as follows:
\begin{equation}
\label{eq:exchange-rel-rewrite}
 x'_k = \frac{\hat y_k + 1}{y_k \oplus 1} \, x_k^{-1}
\prod_i x_i^{[-b_{ik}]_+} \, .
\end{equation}
Using this, for $j \neq k$, we get
\[
\hat y'_j = y'_j (x'_k)^{-b_{kj}}
\prod_{i \neq k} x_i^{b'_{ij}} =
y_j y_k^{[b_{kj}]_+} (\hat y_k + 1)^{-b_{kj}}
x_k^{b_{kj}} \prod_{i \neq k}
x_i^{b'_{ij} - [-b_{ik}]_+ b_{kj}} \ .
\]
Applying \eqref{eq:alt-matrix-mutation}, we get
\[
\hat y'_j  = (y_j \prod_{i \neq j} x_i^{b_{ij}})
(y_k \prod_{i \neq k} x_i^{b_{ik}})^{[b_{kj}]_+}
(\hat y_k + 1)^{-b_{kj}} = \hat y_j \hat y_k^{[b_{kj}]_+}
(\hat y_k + 1)^{-b_{kj}} \,,
\]
as required.
\end{proof}

To draw some conclusions from
Proposition~\ref{pr:another-y-pattern},
we need one more piece of notation.

\begin{definition}
[\emph{Rational functions $Y_{\ell,t}$}]
\label{def:Y-functions}
We denote by
\[
Y_{j;t} 
= Y_{j;t}^{B^0;t_0} \in
 \Qsf(y_1, \dots, y_n) \,.
\]
the components of the $Y$-pattern with values in the universal
semifield $\Qsf(y_1, \dots, y_n)$ (see
Definition~\ref{def:semifield-universal}), and the initial $Y$-seed
$((y_1, \dots, y_n), B^0)$ at~$t_0$.

\end{definition}

\begin{example}
For a $Y$-pattern of type~$A_2$ (see
Example~\ref{example:A2-pattern}),
the functions $Y_{j;t}^{B^0;t_0}$ appear in the third column of
Table~\ref{table:seeds-A2} (replace $\oplus$ by $+$ throughout).
\end{example}

By the universality property of the semifield $\Qsf(y_1, \dots, y_n)$,
for any $Y$-pattern $(\yy_t,B_t)$ with values in an arbitrary
semifield~$\PP$ and the initial $Y$-seed
$((y_1, \dots, y_n), B^0)$ at~$t_0$, the elements
$y_{j;t} \in \PP$ are given by
\begin{equation}
\label{eq:yY-general}
y_{j;t} =  Y_{j;t}^{B^0;t_0}|_\PP(y_1, \dots, y_n).
\end{equation}
In particular, using the notation of Definition~\ref{def:Y-functions},
Proposition~\ref{pr:another-y-pattern} can be restated as
\begin{equation}
\label{eq:y-and-haty-2}
\hat y_{j;t} =  Y_{j;t}^{B^0;t_0}|_\Fcal(\hat y_1, \dots, \hat y_n),
\end{equation}
where we abbreviate $\hat y_j=\hat y_{j;t_0}\,$.
Another important specialization of \eqref{eq:yY-general} is for
the patterns of geometric type: using the notation in
\eqref{eq:yt-geometric}--\eqref{eq:tilde-B}, in this case we
have
\begin{equation}
\label{eq:yY-geometric}
Y_{j;t}|_{\Trop(x_{n+1},\dots,x_m)}(y_1, \dots, y_n) =
\prod_{i=n+1}^m x_i^{b^t_{ij}}\,,
\end{equation}
where $y_j = \prod_{i=n+1}^m x_i^{b^0_{ij}}$.

Comparing Definitions~\ref{def:Aprin} and~\ref{def:Y-functions}, we
note that the rational functions $Y_{j;t}$ 
are defined by means of ``general coefficients'' while the
definition of the polynomials~$F_{i;t}$ involves principal coefficients.
Still, both $Y_{j;t}$ and $F_{i;t}$ lie in the same semifield
$\Qsf(y_1, \dots, y_n)$.
In fact, these two families are closely related, as we will now see.


\begin{proposition}
\label{pr:Y-F-2}
Let $t \overunder{\ell}{} t'$ and $B_t=(b_{ij})$. Then
\begin{equation}
\label{eq:Y-F-2}
F_{\ell;t'} = \frac{Y_{\ell;t} + 1}
{(Y_{\ell;t} + 1)|_{\Trop(y_1, \dots, y_n)}} \cdot
F_{\ell;t}^{-1} \prod_{i=1}^n F_{i;t}^{[-b_{i\ell}]_+} \ .
\end{equation}
\end{proposition}

\begin{proof}
The equality \eqref{eq:Y-F-2}
is a special case of~\eqref{eq:exchange-rel-rewrite}
for the cluster algebra with principal
coefficients 
under the specialization $x_1\!=\!\cdots\!=\!x_n\!=\!1$.
\end{proof}


\begin{proof}[Proof of Theorem~\ref{th:reduction-principal}]
We prove \eqref{eq:xjt-reduction-principal} by induction on the
distance between~$t$ and~$t_0$ in the tree~$\TT_n$.
The equality in question is trivial for $t = t_0$, both sides of
\eqref{eq:xjt-reduction-principal} being equal to~$x_i$.
So it suffices to show the following: if $t \overunder{\ell}{} t'$,
and all $x_{i;t}$ satisfy \eqref{eq:xjt-reduction-principal} then
the same is true for $x_{\ell;t'}$.

Applying \eqref{eq:exchange-rel-rewrite}, we obtain
\begin{equation}
\label{eq:x-X}
x_{\ell;t'} = \frac{(Y_{\ell;t} + 1)|_\Fcal(\hat y_1, \dots, \hat y_n)}
{(Y_{\ell;t} + 1)|_\PP(y_1, \dots, y_n)} \cdot
x_{\ell;t}^{-1} \prod_{i=1}^n x_{i;t}^{[-b_{i\ell}]_+} \ ,
\end{equation}
where the elements $\hat y_j=\hat y_{j;t_0} \in \Fcal$ are given by \eqref{eq:yhat},
and $B = B_t$.
Applying \eqref{eq:x-X} to the cluster algebra with principal
coefficients at~$t_0$, we get
\begin{equation}
\label{eq:x-X-principal}
X_{\ell;t'} = \frac{(Y_{\ell;t} + 1)(\hat y_1, \dots, \hat y_n)}
{(Y_{\ell;t} + 1)|_{\Trop(y_1, \dots, y_n)}(y_1, \dots, y_n)} \cdot
X_{\ell;t}^{-1} \prod_{i=1}^n X_{i;t}^{[-b_{i\ell}]_+}
\end{equation}
(note that further specializing all~$x_i$ to~$1$ leads to
\eqref{eq:Y-F-2}).
To see that $x_{\ell;t'}$ satisfies
\eqref{eq:xjt-reduction-principal}, it suffices to substitute
for each $x_{i;t}$ in the right-hand side of \eqref{eq:x-X} its expression given by
\eqref{eq:xjt-reduction-principal}, and then to use
\eqref{eq:x-X-principal} in conjunction with
the equality \eqref{eq:Y-F-2} evaluated in~$\PP$.
\end{proof}

The following proposition establishes another relation between
the rational functions $Y_{j;t}$ and the polynomials $F_{i;t}$.

\begin{proposition}
\label{pr:Y-F-1}
For any $t \in \TT_n$ and $j \in [1,n]$, we have
\begin{equation}
\label{eq:Y-F-1}
Y_{j;t} = Y_{j;t}|_{\Trop(y_1, \dots, y_n)}
\prod_{i=1}^n F_{i;t}^{b_{ij}} \ ,
\end{equation}
where $(b_{ij}) = (b_{ij}^t) = B_t$ is the exchange matrix at~$t$.
\end{proposition}

\begin{proof}
Apply \eqref{eq:y-and-haty-2} to the cluster algebra with principal
coefficients 
and specialize $x_1\!=\!\cdots\!=\!x_n\!=\!1$
(implying $\hat y_i\!=\!y_i$).
Then  \eqref{eq:Y-F-1} becomes a special case of~\eqref{eq:yhat}.
\end{proof}

\begin{remark}
\label{rem:Ytrop-principal}
We note that the monomial $Y_{j;t}|_{\Trop(y_1, \dots, y_n)}$
appearing in~\eqref{eq:Y-F-1} is nothing but the value of $y_{j;t}$
in the cluster algebra with principal coefficients at~$t_0\,$.
For our running example of type~$A_2\,$,
the monomials $Y_{j;t}|_{\Trop(y_1,\dots,y_n)}$ appear in the third
column of Table~\ref{table:seeds-geom-A2}.
\end{remark}

Proposition~\ref{pr:Y-F-1} has the following useful corollary.

\begin{proposition}
\label{pr:laurent-Y-sink}
Suppose that $t\in\TT_n$ and $j \in [1,n]$ are such that
$b_{ij}^t\geq 0$ for all~$i$.
Then $Y_{j;t}$ is a Laurent polynomial in $y_1,\dots,y_n\,$.
\end{proposition}

\begin{proof}
Immediate from \eqref{eq:Y-F-1}.
\end{proof}

\section{Exchange graphs}
\label{sec:exchange-graphs}

We start by recalling some basic definitions from \cite{ca1,ca2}.

\begin{definition}
[\emph{Seeds}]
\label{def:seed-equivalence}
We say that two labeled seeds $\Sigma=(\xx, \yy, B)$
and $\Sigma'=(\xx', \yy', B')$
define the same \emph{seed}
if $\Sigma'$ is obtained from $\Sigma$ by simultaneous relabeling of the
$n$-tuples $\xx$ and $\yy$ and the corresponding relabeling of
the rows and columns of~$B$.
In other words, seeds are equivalence classes of labeled seeds, where
\[
\Sigma=(\xx, \yy, B),\quad \xx=(x_1,\dots,x_n),\quad
\yy=(y_1,\dots,y_n),\quad B=(b_{ij})
\]
and
\[
\Sigma'=(\xx', \yy', B'),\quad \xx'=(x'_1,\dots,x'_n),\quad
\yy'=(y'_1,\dots,y'_n),\quad B'=(b'_{ij})
\]
are equivalent (denoted $\Sigma\sim\Sigma'$)
if there exists a
permutation~$\sigma$ of indices~$1, \dots, n$ such that
\[
x'_i = x_{\sigma(i)}, \quad y'_j = y_{\sigma(j)}, \quad
b'_{i,j} = b_{\sigma(i), \sigma(j)}
\]
for all~$i$ and~$j$.
We denote by $[\Sigma]$ the seed represented by a labeled seed
$\Sigma$.
\end{definition}

For instance, in Example~\ref{example:A2-pattern}, the labeled
seeds $\Sigma_0$ and $\Sigma_5$ define the same seed, i.e.,
$[\Sigma_0] = [\Sigma_5]$.

The (unlabeled) \emph{cluster} of a seed 
is an (unordered)
$n$-element subset of~$\Fcal$ which is a free generating set for~$\Fcal$
(cf. Definition~\ref{def:seed}).
In view of Definition~\ref{def:seed-mutation},
for every element $z$ of the cluster of a seed
$[\Sigma]$, there is a well-defined seed
$\mu_z([\Sigma]) = [\mu_k(\Sigma)]$, where $z = x_k$ in the
labeling provided by~$\Sigma$; we say that the seed
$\mu_z([\Sigma])$ is obtained from $[\Sigma]$ by mutation in
direction~$z$.

The combinatorics of seed mutations
is captured by the exchange graph of a cluster algebra.

\begin{definition}
[\emph{Exchange graphs}]
\label{def:exchange-graph}
The \emph{exchange graph} of a cluster pattern (and of the
corresponding cluster algebra) is the $n$-regular (finite or infinite)
connected graph whose vertices are the seeds of the cluster pattern and whose
edges connect the seeds related by a single mutation.

To illustrate,   the exchange graph in
Example~\ref{example:A2-pattern} is a $5$-cycle.

The exchange graph can be obtained as a quotient of the tree~$\TT_n$
modulo the equivalence relation on vertices defined by setting
$t \sim t'$ whenever $\Sigma_t\sim\Sigma_{t'}$
(equivalently, $[\Sigma_t] = [\Sigma_{t'}]$).

Each exchange graph is endowed with a canonical (discrete) \emph{connection},
in the sense of \cite{bgh, guillemin-zara}.
That is, for any two adjacent vertices $[\Sigma]$ and $[\Sigma']$ in
the exchange graph, the $n-1$ edges incident to $[\Sigma]$ but not
to $[\Sigma']$ are in a natural bijection with the edges
incident to $[\Sigma']$ but not to $[\Sigma]$.
Specifically, an edge corresponding to the exchange of a cluster
variable $z$ from $[\Sigma]$ is matched to the  edge corresponding to
the exchange of $z$ from~$[\Sigma']$.
\end{definition}

Clearly, the exchange graph of a cluster algebra depends only on
its underlying $Y$-pattern, which is in turn determined by a
$Y$-seed at an arbitrary vertex~$t_0$ in~$\TT_n$.
The following stronger statement was conjectured in
\cite[Conjecture~4.14(1)]{cdm}.

\begin{conjecture}
\label{con:exchange-graph-independence}
The exchange graph of a cluster algebra~$\Acal(\xx, \yy, B)$
(and the canonical connection on this graph)
depends only on the matrix~$B$.
\end{conjecture}

\begin{definition}
[\emph{Finite type}]
\label{def:finite-type}
A cluster algebra $\Acal$ 
is of \emph{finite type} if its exchange graph is finite, that is, $\Acal$
has finitely many distinct seeds. 
\end{definition}

Cluster algebras of finite type were classified in \cite{ca2}:
they correspond to finite root systems.
As shown in \cite{ca2}, the property of a seed $(\xx, \yy, B)$ to
define a cluster algebra of finite type depends only on the matrix~$B$.

For cluster algebras of finite type,
Conjecture~\ref{con:exchange-graph-independence} was proved
in~\cite[Theorem~1.13]{ca2}.
It is also known to hold for $n\leq 2$ \cite[Example~7.6]{ca1}.

As an application of the results in Section~\ref{sec:separation},
we obtain a partial result towards
Conjecture~\ref{con:exchange-graph-independence}.

\begin{definition}
[\emph{Coverings of exchange graphs}]
\label{def:covering}
Let $\Acal'$ and $\Acal$ be two cluster algebras
defined by cluster patterns $t\mapsto \Sigma_t=(\xx_t,\yy_t,B_t)$
and $t\mapsto \Sigma'_t=(\xx'_t,\yy'_t,B_t)$
with the same exchange matrices~$B_t\,$.
We say that the exchange graph of~$\Acal'$ \emph{covers}
the exchange graph of~$\Acal$ if
the canonical projection of $\TT_n$ onto the exchange
graph of~$\Acal$ factors through the similar projection
for~$\Acal'$; in other words,
$\Sigma'_{t_1}\sim\Sigma'_{t_2}$ implies
$\Sigma_{t_1}\sim\Sigma_{t_2}$.
\end{definition}

\begin{theorem}
\label{th:exchange-graph-principal}
The exchange graph of an arbitrary cluster algebra $\Acal$
is covered by the exchange graph
of the algebra~$\Aprin(t_0)$ which has the same exchange matrices~$B_t$ ($t\in\TT_n$)
as $\Acal$ and has principal coefficients at some vertex~$t_0$.
\end{theorem}

In other words, among all cluster algebras with a given exchange
matrix~$B^0 = B_{t_0}\,$, the one with principal coefficients has the
``largest'' exchange graph, i.e.,
one that covers all exchange graphs in this class.
It is immediate from Theorem~\ref{th:exchange-graph-principal} that
the exchange graph of $\Aprin(t_0)$
is independent of the choice of~$t_0$.

It follows from Definition~\ref{def:covering} that the ``smallest'' exchange
graph---that is, the one covered by all others---arises in the case of
\emph{trivial} coefficients (all equal to~$1$),
i.e., in the case of a cluster algebra over the one-element semifield~$\PP=\{1\}$.


\pagebreak[3]

\begin{proof}[Proof of Theorem~\ref{th:exchange-graph-principal}]
Suppose that in $\Aprin(t_0)$,
the labeled seeds at two vertices~$t_1$ and~$t_2$
are equivalent. 
We need to show that $\Sigma_{t_1}\sim\Sigma_{t_2}$
in an arbitrary cluster algebra~$\Acal$
with the same exchange matrix~$B^0$ at~$t_0$.
By definition, we have
$X_{i;t_2}^{B^0;t_0} = X_{\sigma(i);t_1}^{B^0;t_0}$,
$Y_{j;t_2}^{B^0;t_0}|_{\Trop(y_1, \dots, y_n)} =
 Y_{\sigma(j);t_1}^{B^0;t_0}|_{\Trop(y_1, \dots, y_n)}$, and
$b_{ij}^{t_2} = b_{\sigma(i), \sigma(j)}^{t_1}$ for some permutation
$\sigma$ of $[1,n]$.
This implies in particular that $F_{i;t_2}^{B^0;t_0} = F_{\sigma(i);t_1}^{B^0;t_0}$.
We need to show that in~$\Acal$, we have $x_{i;t_2} = x_{\sigma(i);t_1}$
and $y_{j;t_2} = y_{\sigma(j);t_1}$ for all~$i$ and~$j$.
The former equality follows at once from
\eqref{eq:xjt-reduction-principal}; the latter one follows from
the equality $Y_{j;t_2}^{B^0;t_0} =
Y_{\sigma(j);t_1}^{B^0;t_0}$, which is a consequence of
\eqref{eq:Y-F-1}.
\end{proof}

The following conjecture asserts that the exchange graph of
$\Aprin(t_0)$ has a purely combinatorial description in terms of
matrix mutations.

\begin{conjecture}
\label{con:principal-exchange-graph-combinatorial}
Let $(\tilde B_t)_{t \in \TT_n}$ be a family of $2n \times n$ extended
exchange matrices associated with the algebra $\Aprin(t_0)$ (see
Remark~\ref{rem:principal-tildeB}).
Then the labeled seeds at two vertices $t$ and $t'$ define
the same seed of $\Aprin(t_0)$ if and only if $\tilde B_{t'}$ is
obtained from $\tilde B_t$ by a simultaneous permutation of rows
and columns of the principal part~$B_t$. \linebreak[3]
In other words, in the cluster algebra~$\Aprin(t_0)$,
the cluster of a seed $(\xx, \yy, B)$ is uniquely
determined by the corresponding $Y$-seed $(\yy, B)$.
\end{conjecture}

Conjecture~\ref{con:principal-exchange-graph-combinatorial} is
not yet checked even for algebras of finite type.
Note that it implies the following combinatorial characterization
of finite type exchange matrices.

We say that a matrix $\tilde B$ has \emph{finite mutation type} if its mutation
equivalence class is finite, i.e., only finitely many
matrices can be obtained from $\tilde B$ by repeated matrix
mutations.

\begin{conjecture}
\label{con:finite-type-via-2nbyn-matrices}
An exchange matrix $B_{t_0}$ gives rise to cluster algebras of
finite type if and only if the corresponding $2n \times n$ matrix
$\tilde B_{t_0}$ defined as in Remark~\ref{rem:principal-tildeB}
is of finite mutation type.
\end{conjecture}

It is the ``if" part of Conjecture~\ref{con:finite-type-via-2nbyn-matrices}
that still remains open.
At the moment we can only prove the following weaker statement.

\begin{proposition}
\label{pr:finite-eq-finite-mutation}
An exchange matrix $B$ gives rise to cluster algebras of
finite type if and only if any integer matrix $\tilde B$ having
$B$ as its principal part is of finite mutation type.
\end{proposition}

Again, only the ``if" part needs a proof.
It will be given in Section~\ref{sec:denom-bipartite}.

\section{Properties of $F$-polynomials}
\label{sec:F-polynomials}

In this section we discuss some properties (partly conjectural) of
the polynomials $F_{\ell;t} = F_{\ell;t}^{B^0;t_0}\in
\ZZ[y_1, \dots, y_n]$ given by \eqref{eq:F-def}.
We start by writing down, in a closed form, the recurrence relations
for these polynomials.

\begin{proposition}
\label{pr:F-recurrence}
Let $t \mapsto \tilde B_t=(b^t_{ij}) \, (t\in\TT_n)$ be the family of
$2n\times n$ 
matrices associated with the cluster algebra $\Aprin(B_{t_0})$.
That is, $\tilde B_{t_0}$ is the matrix in
Remark~\ref{rem:principal-tildeB}, and
$\tilde B_{t'} = \mu_k(\tilde B_t)$ whenever
$t \overunder{k}{} t'$.
Then the polynomials
$F_{\ell;t} = F_{\ell;t}^{B_{t_0};t_0}(y_1,\dots,y_n)$
are uniquely determined by the initial conditions
\begin{equation}
\label{eq:F-initial}
F_{\ell;t_0} = 1 \quad (\ell = 1,
\dots, n) \ ,
\end{equation}
together with the recurrence relations
\begin{align}
\label{eq:F-mut1}
F_{\ell;t'}&= F_{\ell;t} \quad \text{for $\ell\neq
  k$;}\\
\label{eq:F-mut2}
F_{k;t'} &= \frac{\prod_{j=1}^n y_j^{[b_{n+j,k}^t]_+}
\prod_{i=1}^n F_{i;t}^{[b_{ik}^t]_+}
+
\prod_{j=1}^n y_j^{[-b_{n+j,k}^t]_+}
\prod_{i=1}^n F_{i;t}^{[-b_{ik}^t]_+}}{F_{k;t}}\,,
\end{align}
for every edge $t \overunder{k}{} t'$ such that $t$ lies on
the (unique) path from $t_0$ to $t'$ in~$\TT_n\,$.
\end{proposition}

\begin{proof}
The relations \eqref{eq:F-mut1}--\eqref{eq:F-mut2} are immediate consequences
of the exchange relations in the cluster algebra with principal
coefficients at~$t_0$, which are a special case of
\eqref{eq:exch-rel-geom}:
\begin{equation}
\label{eq:X-mut}
X_{k;t}\, X_{k;t'}
=
\prod_{j=1}^n y_j^{[b_{n+j,k}^t]_+}
\prod_{i=1}^n X_{i;t}^{[b_{ik}^t]_+}
+
\prod_{j=1}^n y_j^{[-b_{n+j,k}^t]_+}
\prod_{i=1}^n X_{i;t}^{[-b_{ik}^t]_+} \, .
\end{equation}
To deduce \eqref{eq:F-mut2} from \eqref{eq:X-mut}, specialize all initial cluster
variables to~$1$.
\end{proof}

We next draw some consequences from
\eqref{eq:xjt-reduction-principal}.

\begin{proposition}
\label{pr:F-primitive}
Each of the polynomials $F_{\ell;t}(y_1, \dots, y_n)$ is not
divisible by any~$y_j$.
\end{proposition}

\begin{proof}
Applying \eqref{eq:xjt-reduction-principal} to the algebra with
principal coefficients at~$t_0$, we conclude that the denominator on
the right-hand side becomes equal to~$1$, that is,
\begin{equation}
\label{eq:F-primitive}
F_{\ell;t}^{B^0;t_0}|_{\Trop(y_1, \dots, y_n)} = 1\, ,
\end{equation}
which is precisely our statement.
\end{proof}

Our next result relates the polynomials  $F_{\ell;t}^{B^0;t_0}$
and $F_{\ell;t}^{-B^0;t_0}$.

\begin{proposition}
\label{pr:FB-F-B}
The polynomials  $F_{\ell;t}^{B^0;t_0}$
and $F_{\ell;t}^{-B^0;t_0}$ are related by
\begin{equation}
\label{eq:FB-F-B}
F_{\ell;t}^{B^0;t_0}(y_1,\dots,y_n)
= \frac{F_{\ell;t}^{-B^0;t_0}(y_1^{-1},\dots,y_n^{-1})}
{F_{\ell;t}^{-B^0;t_0}|_{\Trop(y_1, \dots, y_n)}(y_1^{-1},\dots,y_n^{-1})} \, .
\end{equation}
\end{proposition}

\begin{proof}
Using the notation in Definition~\ref{def:Aprin},
let us introduce the cluster algebra~$\Acal'$
which shares with $\Acal=\Aprin(B^0)$ the same
ambient field, the same coefficient semifield
$\PP = \Trop(y_1, \dots, y_n)$ and the same initial cluster
$\xx_{t_0} = (x_1, \dots, x_n)$, but has the initial $Y$-seed
$(\yy'_{t_0}\,,B'_{t_0})$
given by
$\yy'_{t_0} = (y_1^{-1}, \dots, y_n^{-1}), \,\,
B'_{t_0} = -B^0$.
Comparing the definitions, we see that all the cluster variables
$x_{\ell;t}$ in~$\Acal'$ are the same as in~$\Acal$.
Applying \eqref{eq:xjt-reduction-principal}, we get
\begin{equation}
\label{eq:XB-X-B}
X_{\ell;t}^{B^0;t_0}(x_1, \dots, x_n; y_1,\dots,y_n)
= \frac{X_{\ell;t}^{-B^0;t_0}(x_1, \dots, x_n; y_1^{-1},\dots,y_n^{-1})}
{F_{\ell;t}^{-B^0;t_0}|_{\Trop(y_1, \dots, y_n)}(y_1^{-1},\dots,y_n^{-1})} \, ,
\end{equation}
which implies \eqref{eq:FB-F-B} by specializing all~$x_i$ to~$1$.
\end{proof}

We now present a tantalizing conjecture which is a strengthening of
Proposition~\ref{pr:F-primitive}.

\begin{conjecture}
\label{con:conjectures-on-denoms}
Each polynomial $F_{\ell;t}^{B^0;t_0}(y_1, \dots, y_n)$ has constant term~$1$.
\end{conjecture}

In view of Proposition~\ref{pr:FB-F-B}, Conjecture~\ref{con:conjectures-on-denoms}
is equivalent to the following.

\begin{conjecture}
\label{con:conjectures-on-denoms-2}
Each polynomial $F_{\ell;t}^{B^0;t_0}(y_1, \dots, y_n)$ has a unique
monomial of maximal degree.
Furthermore, this monomial has coefficient~$1$, and it is
divisible by all the other occurring monomials.
\end{conjecture}

Some special cases of Conjecture~\ref{con:conjectures-on-denoms}
are proved in subsequent sections.
Here we give some equivalent reformulations.

\begin{proposition}
\label{pr:equivalent-conjectures-on-denoms}
For a given initial seed,
the following are equivalent:
\begin{enumerate}
\item[(i)]
Each polynomial $F_{\ell;t}(y_1, \dots, y_n)$ has constant term~$1$.
\item[(ii)]
Each tropical evaluation $Y_{\ell;t}|_{\Trop(y_1, \dots, y_n)}$
is a monomial in $y_1, \dots,y_n$ which has either all exponents
nonnegative or all exponents nonpositive.
\item[(iii)]
In the corresponding cluster algebra with principal coefficients,
in each exchange relation
\eqref{eq:exchange-rel-xx}/\eqref{eq:exch-rel-usual},  exactly
one of the coefficients $p_k^\pm$ is
equal to~$1$.
\end{enumerate}
\end{proposition}

\begin{proof}
In the notation of Proposition~\ref{pr:F-recurrence},
we have
\begin{equation}
\label{eq:Y-tropical-principal}
Y_{\ell;t}|_{\Trop(y_1, \dots, y_n)} =
\prod_{j=1}^n y_j^{b^t_{n+j,\ell}}\, ;
\end{equation}
this is a special case of \eqref{eq:yY-geometric}.
Thus (ii) can be rephrased as follows:
\begin{enumerate}
\item[(ii${}'$)]
For any~$\ell$ and~$t$, the integer vector
$\cc_{\ell;t} = (b^t_{n+1,\ell}, \dots, b^t_{2n,\ell})$
has either all components
nonnegative or all components nonpositive.
\end{enumerate}
In view of \eqref{eq:X-mut}, the coefficients in the
exchange relations are of the form
\[
p^\pm  = \prod_{j=1}^n y_j^{[\pm b^t_{n+j,\ell}]_+} \,.
\]
Hence (iii)$\Rightarrow$(ii${}'$).
To prove the converse implication (ii${}'$)$\Rightarrow$(iii),
it remains to show that (ii${}'$) implies that
for any~$\ell$ and~$t$, we have $\cc_{\ell;t} \neq 0$ .
Recall that $\tilde B_{t'} = \mu_k(\tilde B_t)$ whenever
$t \overunder{k}{} t'$.
Under the assumption (ii${}'$), the matrix mutation rule \eqref{eq:matrix-mutation}
can be rewritten as follows:
\begin{equation}
\label{eq:extended-matrix-mutation}
\cc_{\ell;t'} =
\begin{cases}
-\cc_{\ell;t} & \text{if $\ell=k$;} \\[.05in]
\cc_{\ell;t} + [b_{k\ell}]_+ \ \cc_{k;t}
 & \text{if $b^t_{n+i,k} \geq 0$ for all $i$;}\\[.05in]
\cc_{\ell;t} + [-b_{k\ell}]_+ \ \cc_{k;t}
 & \text{if $b^t_{n+i,k} \leq 0$ for all $i$.}
\end{cases}
\end{equation}
It follows that, as long as (ii${}'$) holds, matrix mutations act
by invertible linear transformations on the bottom part
of~$\tilde B_t\,$. The latter must therefore have full rank, and so
has no zero columns.

It remains to demonstrate that (i)$\Leftrightarrow$(iii).
The implication (i)$\Rightarrow$(iii) is immediate
from~\eqref{eq:F-mut2}: if all participating $F$-polynomials
take value~$1$ at $y_1 = \cdots = y_n = 0$ then exactly
one of the coefficients on the right must be equal to~$1$.
The converse implication (iii)$\Rightarrow$(i) also follows from
\eqref{eq:F-mut2} by induction on the distance from $t$ to~$t_0$
in~$\TT_n\,$.
\end{proof}

\section{Principal $\ZZ^n$-gradings}
\label{sec:principal-gradings}

In this section we introduce and study a natural
$\ZZ^n$-grading in a cluster algebra with principal coefficients.
As an application, we obtain, in Corollary~\ref{cor:xjt=F/F} below, a refinement of
Theorem~\ref{th:reduction-principal}.

\begin{proposition}
\label{pr:homogeneity}
Every Laurent polynomial $X_{\ell;t}^{B^0;t_0}$ is homogeneous with respect to
the $\ZZ^n$-grading in
$\ZZ [x_1^{\pm 1}, \dots, x_n^{\pm 1}; y_1, \dots, y_n]$ given by
\begin{equation}
\label{eq:degrees-xy}
{\rm deg}(x_i) = \ee_i, \quad {\rm deg}(y_j) = - \bb_j^0 \ ,
\end{equation}
where $\ee_1, \dots, \ee_n$ are the standard  basis (column) vectors in
$\ZZ^n$, and $\bb_j^0 = \sum_i b^0_{ij}\,\ee_i$ is the $j$th
column of~$B^0$.
\end{proposition}

\begin{proof}
Let $\hat y_1\,,\dots,\hat y_n$ be given by
\begin{equation}
\label{eq:yhat-0}
\hat y_j = y_j \prod_i x_i^{b^0_{ij}}
\end{equation}
(cf.~\eqref{eq:yhat}).
In view of \eqref{eq:degrees-xy}, we have
\begin{equation}
\label{eq:degrees-yhat}
{\rm deg}(\hat y_j) = 0\,.
\end{equation}
Now let us apply \eqref{eq:x-X-principal}.
Using induction on the distance from~$t_0$ in~$\TT_n$, it suffices
to show that the element $X_{\ell;t'}$ on the left-hand side is
homogeneous provided so are all the elements $X_{i;t}$ on the
right-hand side.
In other words, it remains to show the homogeneity of
$$\frac{(Y_{\ell;t} + 1)(\hat y_1, \dots, \hat y_n)}
{(Y_{\ell;t} + 1)|_{\Trop(y_1, \dots, y_n)}(y_1, \dots, y_n)}\ .$$
The numerator $(Y_{\ell;t} + 1)(\hat y_1, \dots, \hat y_n)$ is
homogeneous of degree~$0$ by \eqref{eq:degrees-yhat}.
The denominator is a Laurent monomial in
$y_1, \dots, y_n$, so it is homogeneous as well.
\end{proof}

Proposition~\ref{pr:homogeneity} can be restated as follows
(cf.\ Definitions~\ref{def:principal-coeffs} and~\ref{def:Aprin}).

\begin{corollary}
\label{cor:principal-grading}
Under the $\ZZ^n$-grading given by \eqref{eq:degrees-xy},
the cluster algebra $\Aprin(B_{t_0})$
is a $\ZZ^n$-graded
subalgebra of $\ZZ [x_1^{\pm 1}, \dots, x_n^{\pm 1}; y_1, \dots,
  y_n]$.
All cluster variables in $\Aprin(B_{t_0})$ are homogeneous elements.
\end{corollary}

We will use the notation
\begin{equation}
\label{eq:degree-vector}
\gg_{\ell;t} = \gg_{\ell;t}^{B^0;t_0} 
= \column{g_1}{g_n}
= \deg(X_{\ell;t}^{B^0;t_0}) \in \ZZ^n
\end{equation}
(see \eqref{eq:degrees-xy}) for the (multi-) degrees of the cluster
variables in~$\Aprin(B_{t_0})$, and refer to the integer vectors thus obtained
as \emph{$\gg$-vectors}.

We can now refine the ``separation of additions'' result in
Theorem~\ref{th:reduction-principal}, as follows.

\begin{corollary}
\label{cor:xjt=F/F}
Cluster variables in an arbitrary cluster
algebra
$\Acal=\Acal(\xx_0,\yy_0,B^0)$ can be expressed in terms of the initial seed
\eqref{eq:initial-seed} by the formula
\begin{equation}
\label{eq:xjt=F/F}
x_{\ell;t} = \frac{F_{\ell;t}^{B^0;t_0}|_\Fcal(\hat y_1, \dots, \hat y_n)}
{F_{\ell;t}^{B^0;t_0}|_\PP (y_1, \dots, y_n)} \, x_1^{g_1} \cdots
x_n^{g_n} \,,
\end{equation}
where we use the notation \eqref{eq:yhat-0} and~\eqref{eq:degree-vector}.
\end{corollary}

\begin{proof}
Proposition~\ref{pr:homogeneity} can be restated as follows:
\begin{equation}
\label{eq:X-rescaling}
X_{\ell;t}^{B^0;t_0}(\gamma_1 x_1,\dots,\gamma_n x_n;
\dots,\prod_k \gamma_k^{-b^0_{kj}}y_j,\dots)
=\Bigl(\prod_k \gamma_k^{g_k}\Bigr) X_{\ell;t}^{B^0;t_0}(x_1,\dots,x_n;y_1,\dots,y_n)
\end{equation}
for any rational functions~$\gamma_1,\dots,\gamma_n$.
Applying \eqref{eq:X-rescaling} with
$\gamma_k=x_k^{-1}\,$, we obtain
\begin{equation}
\label{eq:degree-rescaling}
X_{\ell;t}^{B^0;t_0}(x_1,\dots,x_n;y_1,\dots,y_n)
= x_1^{g_1} \cdots x_n^{g_n} F_{\ell;t}^{B^0;t_0}(\hat y_1, \dots, \hat y_n)\,.
\end{equation}
Substituting \eqref{eq:degree-rescaling} into
\eqref{eq:xjt-reduction-principal} yields \eqref{eq:xjt=F/F}.
\end{proof}

\begin{example}
In the special case of type~$A_2$ (Examples~\ref{example:A2-pattern}
and~\ref{example:A2-pattern-principal}),
the $\gg$-vectors $\gg_{\ell;t}=\gg_{\ell;t}^{B^0;t_0}$ can be read off the fourth column of
Table~\ref{table:seeds-geom-A2},
by taking the monomials not involving the $y$ variables and recording
their exponents.
The results are shown in the second column of Table~\ref{table:typeA2-degrees}, where we
use the shorthand
\[
x^{\gg_{\ell;t}}=\prod_i x_i^{g_i}\,.
\]
Substituting
$\hat y_1=y_1x_2^{-1},\,\, \hat y_2 = y_2x_1$
into the fifth column of Table~\ref{table:seeds-geom-A2},
we obtain the Laurent polynomials $F_{\ell;t}(\hat y_1,\hat y_2)$
shown in the third column of Table~\ref{table:typeA2-degrees}.
The fourth column of this table lists the polynomials
$F_{\ell;t}|_\PP(y_1,y_2)$.
Thus, Table~\ref{table:typeA2-degrees} provides all ingredients for
the right-hand side of~\eqref{eq:xjt=F/F} for the cluster algebra
under consideration.
The resulting formulas for the cluster variables
match the ones in the fourth column of Table~\ref{table:seeds-A2}.
\end{example}

\begin{table}[ht]
\begin{equation*}
\begin{array}{|c|cc|cc|cc|}
\hline
&&&&&&\\[-4mm]
t & x^{\gg_{1;t}} & x^{\gg_{2;t}} &
F_{1;t}(\hat y_1,\hat y_2) & F_{2;t}(\hat y_1,\hat y_2)&
F_{1;t}|_\PP & F_{2;t}|_\PP
 \\[1mm]
\hline
&&&&&&\\[-3mm]
0 &
x_1&x_2 &
1 & 1&
1 & 1
\\[4.5mm]
\hline
&&&&&&\\[-3mm]
1 &
x_1&
\dfrac{1}{x_2}&
1 &
x_1y_2+1 &
1 &
y_2\oplus 1
\\[4.5mm]
\hline
&&&&&&\\[-3mm]
2 &
\dfrac{1}{x_1}&
\dfrac{1}{x_2}&
\dfrac{x_1y_1y_2+y_1+x_2}{x_2}&
x_1y_2+1 &
y_1y_2\oplus y_1\oplus 1 &
y_2\oplus 1
\\[4.5mm]
\hline
&&&&&&\\[-3mm]
3 &
\dfrac{1}{x_1}&
\dfrac{x_2}{x_1}&
\dfrac{x_1y_1y_2+y_1+x_2}{x_2}&
\dfrac{y_1+x_2}{x_2} &
y_1y_2\oplus y_1\oplus 1 &
y_1\oplus 1
\\[4.5mm]
\hline
&&&&&&\\[-3mm]
4 &
x_2 &
\dfrac{x_2}{x_1}&
1&
\dfrac{y_1+x_2}{x_2} &
1 &
y_1\oplus 1
\\[4.5mm]
\hline
&&&&&&\\[-3mm]
5 & x_2 & x_1 & 1 & 1 & 1 & 1
\\[4.5mm]
\hline
\end{array}
\end{equation*}
\caption{$\gg$-vectors and Laurent polynomials $F_{\ell;t}(\hat y_1,\hat
  y_2)$ in type~$A_2$}
\label{table:typeA2-degrees}
\end{table}

\begin{remark}
[\emph{The $+/\oplus$-symmetry}]  
The formula \eqref{eq:xjt=F/F} exhibits certain symmetry between
the ordinary addition in~$\Fcal$ and the auxiliary addition in~$\PP$.
To make this symmetry more transparent, we note that this formula
(as well as many other formulas in this paper) extends to the
following more general setup.
Let $M$ be an abelian multiplicative group, with
two distinguished subgroups $\PP$ and $\hat \PP$, each of which is
endowed with a semifield structure, which we write as $(\PP, \oplus)$ and
$(\hat \PP, +)$, respectively.
(In our current setup, think of~$M$ as the multiplicative group
of~$\Fcal$, with $\PP$ the coefficient semifield, and
$\hat \PP$ the subsemifield in~$\Fcal$ generated by all the elements
$\hat y_{j;t}$ given by \eqref{eq:yhat}.)
We can now define a (labeled) \emph{$(M,\PP, \hat \PP)$-seed} as follows:
\begin{equation}
\label{eq:MPhatP-seed}
\Sigma = ((x_1, \dots, x_n; y_1, \dots, y_n; \hat y_1, \dots,
\hat y_n), B^0) \quad (x_j \in M, \,\, y_j \in \PP,
\,\, \hat y_j \in \hat \PP) \ ,
\end{equation}
where $B^0$ is a skew-symmetrizable integer $n \times n$ matrix
(see Definition~\ref{def:seed}), and the elements $x_j, y_j$ and
$\hat y_j$ satisfy \eqref{eq:yhat-0}.
Then, for any $k = 1, \dots, n$, the $(M,\PP, \hat \PP)$-seed
mutation in direction~$k$ transforms $\Sigma$ into
$$\Sigma' = \mu_k(\Sigma) = ((x'_1, \dots, x'_n; y'_1, \dots, y'_n;
{\hat y}'_1, \dots, {\hat y}'_n), B^1)$$
given as follows:
\begin{itemize}
\item
$B^1$ is obtained from $B^0$ by the matrix mutation in
direction~$k$ (see \eqref{eq:matrix-mutation}).
\item
The tuple $(y'_1, \dots, y'_n)$ (resp.,
$({\hat y}'_1, \dots, {\hat y}'_n)$) is obtained from
$(y_1, \dots, y_n)$ (resp.,
$(\hat {y}_1, \dots, \hat {y}_n)$) by the $Y$-seed mutation
\eqref{eq:y-mutation} in the semifield~$\PP$ (resp.,~$\hat\PP$).
\item
$x'_i = x_i$ for $i \neq k$, and $x'_k$ is given by
\eqref{eq:exchange-rel-rewrite} with $b_{ik} = b_{ik}^0$.
\end{itemize}
The fact that $\Sigma'$ is indeed a $(M,\PP, \hat \PP)$-seed
is shown by calculations in the proof of
Proposition~\ref{pr:another-y-pattern}.

Following Definition~\ref{def:patterns}, we can define
a $(M,\PP, \hat \PP)$-pattern as an assignment
of an $(M,\PP, \hat \PP)$-seed $\Sigma_t$
to every vertex $t \in \TT_n$, such that the seeds assigned to the
endpoints of any edge $t \overunder{k}{} t'$ are obtained from each
other by the mutation in direction~$k$.
If the seed $\Sigma_{t_0}$ is given by \eqref{eq:MPhatP-seed},
and $\Sigma_t$ is written as
\begin{equation}
\label{eq:MPhatP-seed-t}
\Sigma_t =((x_{1;t},\dots, x_{n;t}; y_{1;t},\dots, y_{n;t};
\hat y_{1;t},\dots, \hat y_{n;t}), B_t)\,,
\end{equation}
then the elements of $\Sigma_t$ are given  by the expressions
\eqref{eq:yY-general}, \eqref{eq:y-and-haty-2} and
\eqref{eq:xjt=F/F}, where in the last two formulas the
evaluation in~$\Fcal$ is replaced by that in~$\hat \PP$.
These expressions imply the following \emph{$+/\oplus$-duality}:
every $(M,\PP, \hat \PP)$-pattern gives rise to a
$(M, \hat \PP, \PP)$-pattern obtained by replacing each $(M,\PP, \hat \PP)$-seed
\eqref{eq:MPhatP-seed-t} with the $(M, \hat \PP, \PP)$-seed
$$((x_{1;t}^{-1},\dots, x_{n;t}^{-1}; \hat y_{1;t},\dots, \hat y_{n;t};
y_{1;t},\dots, y_{n;t}), B_t) \ .$$
\end{remark}

\medskip

Returning to cluster algebras, we note that,
in practical terms, using the formula \eqref{eq:xjt=F/F} to compute
some cluster variable~$x_{\ell;t}$ in terms of an initial cluster at
a vertex~$t_0$ requires 
calculating the polynomial
$F_{\ell;t}^{B^0;t_0}(y_1,\dots,y_n)$ and the $\gg$-vector $\gg_{\ell;t}^{B^0;t_0}$.
The polynomial $F_{\ell;t}^{B^0;t_0}(y_1,\dots,y_n)$ can be
calculated recursively by Proposition~\ref{pr:F-recurrence}.
As to computation of the $\gg$-vector, it can be accomplished using
recurrences \eqref{eq:deg-mut1}--\eqref{eq:deg-mut2}
below.

\begin{proposition}
\label{pr:gg-recurrence}
In the notation of Proposition~\ref{pr:F-recurrence},
the $\gg$-vectors $\gg_{\ell;t} = \gg_{\ell;t}^{B^0;t_0}$
are uniquely determined by the initial conditions
\begin{equation}
\label{eq:gg-initial}
\gg_{\ell;t} = \ee_\ell \quad (\ell = 1,
\dots, n)
\end{equation}
together with the recurrence relations
\begin{align}
\label{eq:deg-mut1}
\gg_{\ell;t'}&=\gg_{\ell;t} \quad \text{for $\ell\neq
  k$;}\\
\label{eq:deg-mut2}
\gg_{k;t'} &= -\gg_{k;t} +\sum_{i=1}^n [b^t_{ik}]_+ \gg_{i;t}
-\sum_{j=1}^n  [b^t_{n+j,k}]_+ \bb^{0}_j
 \,,
\end{align}
where the (unique) path from $t_0$ to $t'$ in $\TT_n$ ends
with the edge $t \overunder{k}{} t'$, and $\bb^{0}_j$ stands for
the $j$th column of $B^0$.
\end{proposition}

\begin{proof}
The relations \eqref{eq:deg-mut1}--\eqref{eq:deg-mut2} follow
from \eqref{eq:X-mut} by equating the
degree of the product $X_{k;t}\, X_{k;t'}$ with that of the first
term on the right.
\end{proof}

Note that \eqref{eq:deg-mut2} can be replaced by the relation
\begin{equation}
\label{eq:deg-mut3}
\gg_{k;t'} = -\gg_{k;t} +\sum_{i=1}^n [-b^t_{ik}]_+ \gg_{i;t}
-\sum_{j=1}^n  [-b^t_{n+j,k}]_+ \bb^{0}_j\, ,
\end{equation}
obtained by taking the second term on the right in \eqref{eq:X-mut}
instead of the first one.
The fact that \eqref{eq:deg-mut2} and \eqref{eq:deg-mut3} agree is
equivalent to the identity
\begin{equation}
\label{eq:deg-identity}
\sum_{i=1}^n b^t_{ik} \gg_{i;t} = \sum_{j=1}^n  b^t_{n+j,k} \bb^{0}_j\, ,
\end{equation}
which is simply another way of saying that the element
$Y_{k;t}(\hat y_1, \dots, \hat y_n)$ is homogeneous of degree~$0$
(see \eqref{eq:degrees-yhat} and \eqref{eq:y-and-haty-2}).

\begin{example}
\label{example:A2-pattern-degrees}
The results of $\gg$-vector calculations based on
\eqref{eq:deg-mut1}--\eqref{eq:deg-mut2}
for a cluster algebra of type~$A_2$ (Examples~\ref{example:A2-pattern}
and~\ref{example:A2-pattern-principal}) are shown in
the third column of Table~\ref{table:degrees-recursive}.
\end{example}

\begin{table}[ht]
\begin{equation*}
\begin{array}{|c|c|cc|cc|}
\hline
&&&&&\\[-4mm]
t & \tilde B_t & \gg_{1;t} & \gg_{2;t} & \dd_{1;t} & \dd_{2;t}
 \\[1mm]
\hline
&&&&&\\[-3mm]
0 &
\left[\begin{smallmatrix}0&1\\-1&0\\1&0\\0&1\end{smallmatrix}\right] &
\twobyone{1}{0} & \twobyone{0}{1}  &
\twobyone{-1}{0} & \twobyone{0}{-1}
\\[4.5mm]
\hline
&&&&&\\[-3mm]
1 & \left[\begin{smallmatrix}0&-1\\1&0\\1&0\\0&-1\end{smallmatrix}\right] &
\twobyone{1}{0} & \twobyone{0}{-1} &
\twobyone{-1}{0} & \twobyone{0}{1}
\\[4.5mm]
\hline
&&&&&\\[-3mm]
2 & \left[\begin{smallmatrix}0&1\\-1&0\\-1&0\\0&-1\end{smallmatrix}\right] &
\twobyone{-1}{0} & \twobyone{0}{-1}  &
\twobyone{1}{1} & \twobyone{0}{1}
\\[4.5mm]
\hline
&&&&&\\[-3mm]
3 & \left[\begin{smallmatrix}0&-1\\1&0\\-1&0\\-1&1\end{smallmatrix}\right] &
\twobyone{-1}{0} & \twobyone{-1}{1}  &
\twobyone{1}{1} & \twobyone{1}{0}
\\[4.5mm]
\hline
&&&&&\\[-3mm]
4 & \left[\begin{smallmatrix}0&1\\-1&0\\1&-1\\1&0\end{smallmatrix}\right] &
\twobyone{0}{1} & \twobyone{-1}{1}  &
\twobyone{0}{-1} & \twobyone{1}{0}
\\[4.5mm]
\hline
&&&&&\\[-3mm]
5 &
\left[\begin{smallmatrix}0&-1\\1&0\\0&1\\1&0\end{smallmatrix}\right] &
\twobyone{0}{1} & \twobyone{1}{0} &
\twobyone{0}{-1} & \twobyone{-1}{0}
\\[4.5mm]
\hline
\end{array}
\end{equation*}
\caption{
$\gg$-vectors and denominators in type~$A_2$}
\label{table:degrees-recursive}
\end{table}

We now present an alternative recursive description of the
$\gg$-vectors $\gg_{\ell;t}^{t_0}$ that is based on fixing~$\ell$
and~$t$ and allowing the initial vertex~$t_0$ to vary.


\begin{proposition}
\label{pr:degree-transition}
Let $t_0 \overunder{k}{} t_1$, and $B^1 = \mu_k(B^0)$.
Then the $\gg$-vectors
\[
\gg_{\ell;t}^{B^0;t_0} = \column{g_1}{g_n} \quad \text{and}\quad
\gg_{\ell;t}^{B^1;t_1} = \column{g'_1}{g'_n}
\]
are related as
follows:
\begin{equation}
\label{eq:g-transition}
g_i =
\begin{cases}
-g'_k  & \text{if $i = k$};\\[.05in]
g'_i + [-b^0_{ik}]_+ g'_k + b^0_{ik} h'_k
 & \text{if $i \neq k$},
\end{cases}
\end{equation}
where the integer $h'_k$ is defined by
\begin{equation}
\label{eq:hk-def}
\text{$u^{h'_k}=F^{B^1;t_1}_{\ell,t}|_{\Trop(u)}
(u^{[b^0_{k1}]_+},\dots,u^{-1},\dots,u^{[b^0_{kn}]_+})$
($u^{-1}$ in the $k$th position).}
\end{equation}
\end{proposition}

\begin{proof}
We start by preparing some needed notation.
We will work in the cluster algebra~$\Acal=\Aprin(B^0)$
with principal coefficients at a vertex~$t_0$, and adopt the
notation in Definition~\ref{def:Aprin}.
Thus, the seed $\Sigma_{t_0}=(\xx_{t_0}\,,\yy_{t_0}\,,B_{t_0})$
is given by \eqref{eq:initial-seed}, and
the coefficient semifield is $\PP=\Trop(y_1,\dots,y_n)$.
Let $\Sigma_{t_1}=(\xx_{t_1}\,,\yy_{t_1}\,,B_{t_1})$ be the
adjacent seed at~$t_1$ obtained from $\Sigma_{t_0}$ by the
mutation in direction~$k$ given by
\eqref{eq:matrix-mutation} - \eqref{eq:exchange-rel-xx}.
In more detail, $\Sigma_{t_1}$ is given as follows:
\begin{itemize}
\item
Its exchange matrix is  $B_{t_1} = B^1 = \mu_k(B^0)$.
\item
Its coefficient system is $\yy_{t_1} = (y'_1, \dots, y'_n)$, where
the $y'_j$ are given by \eqref{eq:y-mutation} evaluated in~$\PP$:
\begin{equation}
\label{eq:y-mutation-principal}
y'_j =
\begin{cases}
y_k^{-1} & \text{if $j = k$};\\[.05in]
y_j y_k^{[b^0_{kj}]_+}
 & \text{if $j \neq k$}.
\end{cases}
\end{equation}
\item
Its cluster is
$\xx_{t_1} = (x'_1, \dots, x'_n)$, where $x'_i = x_i$ for
$i \neq k$, and
\begin{equation}
\label{eq:exchange-rel-xx-principal}
x'_k = x_k^{-1} (y_k  \prod x_i^{[b_{ik}]_+}
+ \ \prod x_i^{[-b_{ik}]_+})  \ .
\end{equation}
\end{itemize}
Finally, let $(\hat y_1, \dots, \hat y_n)$ be given by
\eqref{eq:yhat-0}, and let $(\hat y'_1, \dots, \hat y'_n)$
be obtained from them by \eqref{eq:y-mutation}, where the addition
is the one in the ambient field $\Fcal$.

To prove \eqref{eq:g-transition}, we evaluate the cluster variable
$x_{\ell;t} \in \Acal$ in two different ways  by applying
the formula \eqref{eq:xjt=F/F}
with respect to the seeds at~$t_0$ and at~$t_1$.
Equating the two expressions yields the identity
\begin{equation}
\label{eq:xjt-t0-t1}
\frac{F_{\ell;t}^{B^0;t_0}|_\Fcal(\hat y_1, \dots, \hat y_n)}
{F_{\ell;t}^{B^0;t_0}|_\PP (y_1, \dots, y_n)} \, \prod_i x_i^{g_i}
= \frac{F_{\ell;t}^{B^1;t_1}|_\Fcal(\hat y'_1, \dots, \hat y'_n)}
{F_{\ell;t}^{B^1;t_1}|_\PP (y'_1, \dots, y'_n)} \, \prod_i (x'_i)^{g'_i}
 \, .
\end{equation}
Note that the denominator on the left of \eqref{eq:xjt-t0-t1}
is equal to~$1$ by \eqref{eq:F-primitive}.
By the same token, it is easy to see that the denominator on the
right of \eqref{eq:xjt-t0-t1} is equal to~$y_k^{h'_k}$
(see \eqref{eq:hk-def}).
Therefore \eqref{eq:xjt-t0-t1} simplifies to
\begin{equation}
\label{eq:xjt-t0-t1-2}
F_{\ell;t}^{B^0;t_0}|_\Fcal(\hat y_1, \dots, \hat y_n) \, \prod_i x_i^{g_i}
= F_{\ell;t}^{B^1;t_1}|_\Fcal(\hat y'_1, \dots, \hat y'_n)
 \,y_k^{-h'_k} \prod_i (x'_i)^{g'_i}
 \, .
\end{equation}

We now specialize \eqref{eq:xjt-t0-t1-2} by setting $x_i = 1$ for all~$i$.
By \eqref{eq:exchange-rel-xx-principal}, $x'_k$ specializes
to~$y_k + 1$, while by \eqref{eq:yhat-0}, each~$\hat y_j$ specializes to~$y_j$.
As for the elements $\hat y'_j$, they specialize to the tuple of
elements of $\Fcal$ obtained from the $y_j$ by \eqref{eq:y-mutation}, with the addition
again understood as the one in~$\Fcal$.
To distinguish the latter elements from the ones in
\eqref{eq:y-mutation-principal}, we denote the specialized $\hat y'_j$
by $\overline y'_j$; we can thus assume that the elements
$\overline y'_j$ belong to the universal semifield
$\Qsf(y_1, \dots, y_n)$ and are given by \eqref{eq:y-mutation}
there.
In this notation, \eqref{eq:xjt-t0-t1-2}
specializes to
\begin{equation}
\label{eq:xjt-t0-t1-spec}
F_{\ell;t}^{B^0;t_0}(y_1, \dots, y_n)
= F_{\ell;t}^{B^1;t_1}(\overline y'_1, \dots, \overline y'_n)
 \,y_k^{-h'_k} (y_k + 1)^{g'_k}
 \, .
\end{equation}

Using the universality of $\Qsf(y_1, \dots, y_n)$, we conclude
that \eqref{eq:xjt-t0-t1-spec} implies the following:
\begin{equation}
\label{eq:xjt-t0-t1-spec-2}
F_{\ell;t}^{B^0;t_0}|_\Fcal(\hat y_1, \dots, \hat y_n)
= F_{\ell;t}^{B^1;t_1}|_\Fcal(\hat y'_1, \dots, \hat y'_n)
 \,\hat y_k^{-h'_k} (\hat y_k + 1)^{g'_k}
 \, .
\end{equation}
Substituting \eqref{eq:xjt-t0-t1-spec-2} into \eqref{eq:xjt-t0-t1-2}
and performing cancellations, we obtain
\begin{equation}
\label{eq:xjt-t0-t1-3}
 \prod_i x_i^{g_i} = (\hat y_k/y_k)^{h'_k} (\hat y_k + 1)^{-g'_k} \prod_i (x'_i)^{g'_i}
 \, .
\end{equation}
It follows from \eqref{eq:yhat-0} and
\eqref{eq:exchange-rel-xx-principal} that
$$\hat y_k/y_k = \prod_i x_i^{b^0_{ik}}, \quad
\hat y_k + 1 = \prod_i x_i^{-[-b^0_{ik}]_+} x_k x'_k \ .$$
Substituting these expressions into \eqref{eq:xjt-t0-t1-3}, we see
that the powers of $x'_k$ on the right cancel out, so both sides
are monomials in $x_1, \dots, x_n$.
The desired formula \eqref{eq:g-transition}
is obtained by equating the exponents of each~$x_i$ on both sides.
\end{proof}

As a byproduct of the above proof, we can express the $\gg$-vector
$\gg_{\ell;t}^{B^0;t_0}$ 
in terms of $F$-polynomials.

\begin{proposition}
\label{pr:g-thru-F}
In the notation of Proposition~\ref{pr:degree-transition}, the
component~$g_k$ of the
$\gg$-vector $\gg_{\ell;t}^{B^0;t_0}$ is
given by
\begin{equation}
\label{eq:gk-thru-F}
\text{$u^{g_k}=\frac{F^{B^0;t_0}_{\ell,t}|_{\Trop(u)}
(u^{[-b^0_{k1}]_+},\dots,u^{-1},\dots,u^{[-b^0_{kn}]_+})}
{F^{B^1;t_1}_{\ell,t}|_{\Trop(u)}
(u^{[b^0_{k1}]_+},\dots,u^{-1},\dots,u^{[b^0_{kn}]_+})}$}
\end{equation}
($u^{-1}$ in the $k$th position).
\end{proposition}

\begin{proof}
Define the integer~$h_k$ by
\begin{equation}
\label{eq:hk-def-2}
\text{$u^{h_k}=F^{B^0;t_0}_{\ell,t}|_{\Trop(u)}
(u^{[-b^0_{k1}]_+},\dots,u^{-1},\dots,u^{[-b^0_{kn}]_+})$
($u^{-1}$ in the $k$th position).}
\end{equation}
Then the desired equality \eqref{eq:gk-thru-F}
takes the form
\begin{equation}
\label{eq:gk=hk-hk'}
g_k=h_k-h'_k \ .
\end{equation}
To prove \eqref{eq:gk=hk-hk'}, we exploit the symmetry
between~$t_0$ and $t_1$.
To do this, we use the formula obtained from
\eqref{eq:xjt-t0-t1-spec} by interchanging $t_0$ and $t_1$.
Comparing \eqref{eq:hk-def-2} with
\eqref{eq:hk-def} (and using the first case in
the matrix mutation rule \eqref{eq:matrix-mutation}),
we see that interchanging $t_0$ and $t_1$ replaces $h'_k$ with
$h_k$, and $g'_k$ with $g_k$; thus, the resulting counterpart of
\eqref{eq:xjt-t0-t1-spec} takes the form
\begin{equation}
\label{eq:xjt-t0-t1-spec-3}
F_{\ell;t}^{B^1;t_1}(\overline y'_1, \dots, \overline y'_n)
= F_{\ell;t}^{B^0;t_0}(y_1, \dots, y_n)
 \,(\overline y'_k)^{-h_k} (\overline y'_k + 1)^{g_k}
 \, .
\end{equation}
Using the fact that $\overline y'_k = y_k^{-1}$ (by the first case
in \eqref{eq:y-mutation}), we can rewrite \eqref{eq:xjt-t0-t1-spec-3}
as follows:
\begin{equation}
\label{eq:xjt-t0-t1-spec-4}
F_{\ell;t}^{B^0;t_0}(y_1, \dots, y_n) =
F_{\ell;t}^{B^1;t_1}(\overline y'_1, \dots, \overline y'_n)
 \,y_k^{g_k - h_k} (y_k + 1)^{-g_k}
 \, .
\end{equation}
Comparing \eqref{eq:xjt-t0-t1-spec-4} with \eqref{eq:xjt-t0-t1-spec}
yields the desired equality \eqref{eq:gk=hk-hk'}.
\end{proof}



We conclude this section with several conjectural properties of the $\gg$-vectors.
The first one sharpens \eqref{eq:gk=hk-hk'}.

\begin{conjecture}
\label{con:hk-gk}
In the  notation of Proposition~\ref{pr:degree-transition} and
\eqref{eq:hk-def-2}, we have
\begin{equation}
\label{eq:hk-gk-conjecture}
h'_k = -[g_k]_+, \quad h_k = -[-g_k]_+ = \min(0, g_k) \, .
\end{equation}
\end{conjecture}

The next conjecture somewhat surprisingly expresses the $\gg$-vector
$\gg_{\ell;t}^{B^0;t_0}$ in terms of
the polynomial $F_{\ell;t}^{B^0;t_0}$ alone.

\begin{conjecture}
\label{con:g-thru-same-F}
Suppose that $B^0, t_0, \ell$ and $t$ are such that the polynomial
$F_{\ell;t}^{B^0;t_0}$ is not identically equal to~$1$.
Then the $\gg$-vector
$\gg_{\ell;t}^{B^0;t_0} = \smallcolumn{g_1}{g_n}$ is given by
\begin{equation}
\label{eq:gk-thru-same-F}
u_1^{g_1} \cdots u_n^{g_n} =
\frac{F^{B^0;t_0}_{\ell,t}|_{\Trop(u_1, \dots, u_n)}
(u_1^{-1},\dots,u_n^{-1})}
{F^{B^0;t_0}_{\ell,t}|_{\Trop(u_1, \dots, u_n)}
\Bigl(\prod_i u_i^{b^0_{i1}}, \dots, \prod_i u_i^{b^0_{in}}\Bigr)} \, .
\end{equation}
\end{conjecture}

To state our last conjecture, the following terminology comes
handy.

\begin{definition}
\label{def:coherent-vectors}
We say that a collection of vectors in~$\ZZ^n$ (or~$\RR^n)$ are
\emph{sign-coherent} (to each other) if, for any $i \in [1,n]$, the $i$th coordinates of all
these vectors are either all nonnegative, or all nonpositive.
\end{definition}

\begin{conjecture}
\label{con:signs-gi}
For any given $B^0, t_0$, and $t$ as above,
the vectors $\gg_{1;t}^{B^0;t_0}, \dots, \gg_{n;t}^{B^0;t_0}$ are sign-coherent.
\end{conjecture}

The significance of these conjectures will become clearer in
the next sections.



\section{Cluster monomials and their parameterizations}
\label{sec:cluster-monomials}


\begin{definition}[{\cite[Definition~4.15]{cdm}}]
\label{def:cluster-monomial}
A \emph{cluster monomial} in a cluster algebra~$\Acal$
is a monomial in cluster variables all of which
belong to the same cluster.
Thus a cluster monomial at a vertex $t \in \TT_n$ is a product
of the form
\begin{equation}
\label{eq:cluster-monomial}
x_{\aa;t} = \prod_\ell x_{\ell;t}^{a_\ell} \in\Fcal \,,
\quad \text{with}
\quad
\aa = \column{a_1}{a_n} \in \ZZ_{\geq 0}^n \,.
\end{equation}
\end{definition}

Cluster monomials are destined to play an important role in the
emerging structural theory of cluster algebras and its applications.
We refer the reader to \cite[Section~4.3]{cdm} for a discussion of
the properties of cluster monomials, both proven and conjectural.
Based on the results for rank~$2$ obtained in \cite{sherzel},
we expect every cluster monomial to belong to the ``canonical basis"
that has yet to be defined for an arbitrary cluster algebra.
Even so, the following basic property has not been established in
general.

\begin{conjecture}[{\cite[Conjecture~4.16]{cdm}}]
\label{con:cl-mon-lin-indep}
Cluster monomials in any cluster algebra are linearly independent over
the ground ring.
\end{conjecture}

For cluster algebras of finite type, this conjecture is proved in
Theorem~\ref{th:cluster-monomials-linear-indep}.

In this section, we discuss two families of combinatorial
parameterizations of cluster monomials by integer vectors
in~$\ZZ^n$. (They are somewhat reminiscent of Lusztig's and string
parameterizations of canonical bases in the theory of quantum groups.)
One of these parametrizations (by \emph{denominator vectors})
has already appeared in~\cite{ca2, cdm};
another one (by \emph{$\gg$-vectors}) is new.
(Note that in Section~\ref{sec:principal-gradings},
$\gg$-vectors were defined in restricted generality,
for the case of principal coefficients only.)

We start by recalling the notion of denominator vectors
(see, e.g., \cite[Section~4.3]{cdm}).
Let $\xx = (x_1, \dots, x_n)$ be a labeled cluster
from some seed in~$\Acal$.
By Theorem~\ref{th:Laurent-phenom}, every nonzero element
$z \in \Acal$ can be uniquely written as
\begin{equation}
\label{eq:Laurent-normal-form}
z = \frac{N(x_1, \dots, x_n)}{x_1^{d_1} \cdots x_n^{d_n}} \, ,
\end{equation}
where $N(x_1, \dots, x_n)$ is a polynomial with coefficients in~$\ZZ \PP$
which is not divisible by any cluster variable~$x_i\in\xx$.
We denote
\begin{equation}
\label{eq:denominator-vector}
\dd(z) = \dd_\xx(z) = \column{d_1}{d_n} \in \ZZ^n,
\end{equation}
and call the integer vector~$\dd(z)$ the \emph{denominator vector}
of~$z$ with respect to the cluster~$\xx$.
Note that the map $z \mapsto \dd(z)$ has the
following multiplicative property:
\begin{equation}
\label{eq:delta-multiplicative}
\dd(z_1 z_2) = \dd(z_1) + \dd(z_2) \,.
\end{equation}

Suppose that~$\Acal$ is associated with a (labeled) cluster
pattern $t \mapsto (\xx_t, \yy_t, B_t)$ on~$\TT_n$, and that~$\xx$
is the cluster of the initial seed \eqref{eq:initial-seed} at~$t_0$.
We will use the notation
\begin{equation}
\label{eq:denom-cluster-variable}
\dd_{\ell;t}^{B^0;t_0} = \dd_\xx(x_{\ell;t})
\end{equation}
for the denominator vectors of cluster variables.
This notation is unambiguous since it is easy to see that
$\dd_\xx(x_{\ell;t})$ is independent of the choice of the
coefficient system $\yy = \yy_{t_0}$; in fact, the mutation
mechanism for generating cluster variables makes it clear that the
whole Newton polytope of $x_{\ell;t}$ (as a Laurent polynomial in
$x_1, \dots, x_n$) is independent of this choice.

More explicitly. the vectors $\dd_{\ell;t} = \dd_{\ell;t}^{B^0;t_0}$ are uniquely
determined by the initial conditions
\begin{equation}
\label{eq:denom-cluster-variable-initial}
\dd_{\ell;t_0}^{B^0;t_0} = -\ee_\ell
\end{equation}
(as before, $\ee_1, \dots, \ee_n$ are the standard basis vectors
in~$\ZZ^n$) together with the recurrence relations implied by
the exchange relation~\eqref{eq:exchange-rel-xx}:
\begin{equation}
\label{eq:exchange-denominator}
\dd_{\ell;t'} =
\begin{cases}
\dd_{\ell;t} & \text{if $\ell\neq k$;}\\
- \dd_{k;t} +
\max\Bigl(\displaystyle\sum_i [b_{ik}^t]_+ \dd_{i;t},
\sum_i [-b_{ik}^t]_+ \dd_{i;t}\Bigr)
& \text{if $\ell=k$}
\end{cases}
\end{equation}
for $t \overunder{k}{} t'$;
here and in what follows, the operations $\max$, $\min$,
$[a]_+$, etc., on vectors are performed component-wise.

\begin{example}
The denominator vectors $\dd_{\ell;t}^{B^0;t_0}$ for a cluster algebra of type~$A_2$
(cf.\ Examples~\ref{example:A2-pattern},
\ref{example:A2-pattern-principal},
and~\ref{example:A2-pattern-degrees}) are shown in
the fourth column of Table~\ref{table:degrees-recursive}.
\end{example}

As above, let~$\Acal$ be a cluster algebra with the initial
exchange matrix $B_{t_0} = B^0$.
If a cluster variable $x_{\ell;t}$ belongs to
the initial cluster, that is, $x_{\ell;t} = x_{i;t_0}$ for some~$i$,
then 
its denominator is given
by~\eqref{eq:denom-cluster-variable-initial}.
Otherwise, we expect
its denominator vector to have the following properties.

\begin{conjecture}
\label{con:denominator-1}
Suppose $x_{\ell;t}$ does not belong to the initial cluster.
Let $\dd_{\ell;t}^{B^0;t_0} \!=\! \smallcolumn{d_1}{d_n}$ \linebreak[3]
be the denominator vector of~$x_{\ell;t}\,$.
Then:
\begin{enumerate}
\item
All components $d_i$ are nonnegative.
\item
We have $d_i = 0$ if and only if
there is a cluster containing both
$x_{\ell;t}$ and~$x_{i;t_0}\,$.
\item
Each component $d_i$ depends only on $x_{\ell;t}$ and $x_{i;t_0}$, not on
the vertices $t_0$ and~$t$.
\end{enumerate}
\end{conjecture}

We note that Conjecture~\ref{con:denominator-1} implies the following
counterpart of Conjecture~\ref{con:signs-gi}
(cf.\ Definition~\ref{def:coherent-vectors}).

\begin{conjecture}
\label{con:sign-coherence-denom}
For any~$t$ and~$t_0\,$, the vectors
$\dd_{1;t}^{B^0;t_0}, \dots, \dd_{n;t}^{B^0;t_0}$ are sign-coherent.
\end{conjecture}


Using \eqref{eq:delta-multiplicative}, we can extend the notation
\eqref{eq:denom-cluster-variable} from cluster variables to
cluster monomials (cf.\ \eqref{eq:cluster-monomial})
by setting
\begin{equation}
\label{eq:denominator-monomial}
\dd_{\aa;t}^{B^0;t_0} = \dd_\xx(x_{\aa;t})
= \sum_{\ell = 1}^n a_\ell \ \dd_{\ell;t}^{B^0;t_0} \ .
\end{equation}
In particular, in view of \eqref{eq:denom-cluster-variable-initial}, we have
\begin{equation}
\label{eq:denom-cluster-monomial-initial}
\dd_{\aa;t_0}^{B^0;t_0} = -\aa \qquad (\aa \in \ZZ^n_{\geq 0}) \ .
\end{equation}

\begin{conjecture}
[{\rm \cite[Conjecture~4.17]{cdm}}]
\label{con:denominator-2}
Different cluster monomials have different denominators with
respect to a given initial seed~$\Sigma_{t_0}\,$.
In particular, for any fixed $t \in \TT_n$, the correspondence
$\aa \mapsto \dd_{\aa;t}^{B^0;t_0}$ is an embedding
$\ZZ^n_{\geq 0} \to \ZZ^n$, hence the vectors
$\dd_{1;t}^{B^0;t_0}, \dots, \dd_{n;t}^{B^0;t_0}$
form a basis of the vector space~$\QQ^n$.
\end{conjecture}

\begin{remark}
\label{rem:denom-even}
If the initial seed $\Sigma_{t_0}$ is \emph{acyclic} (see~\cite{ca3}),
we expect the vectors
$\dd_{1;t}^{B^0;t_0}, \dots, \dd_{n;t}^{B^0;t_0}$
to form a $\ZZ$-basis of the lattice~$\ZZ^n$.
We note that this property does not hold in general.
The simplest counterexample is in finite type~$A_3\,$, for a
non-acyclic seed with the exchange matrix
\[
B^0 = B_{t_0}=\left[\begin{matrix}
0 & 1 & -1\\
-1 & 0 & 1\\
1 & -1& 0
\end{matrix}\right]\,.
\]
Let~$t$ be the vertex connected to~$t_0$ by the path
$$t_0 \overunder{1}{} t_1 \overunder{2}{} t_2 \overunder{3}{} t_3
\overunder{1}{} t \, .$$
Applying the recurrence \eqref{eq:exchange-denominator},
we compute the vectors $\dd_{\ell;t}^{B^0;t_0}$:
\[
\dd_{1;t}^{B^0;t_0} = \begin{bmatrix}0\\ 1\\ 1\end{bmatrix},\quad
\dd_{2;t}^{B^0;t_0} = \begin{bmatrix}1\\ 1\\ 0\end{bmatrix},\quad
\dd_{3;t}^{B^0;t_0} = \begin{bmatrix}1\\ 0\\ 1\end{bmatrix} \,.
\]
These vectors generate a proper sublattice
of index~$2$ in~$\ZZ^3$.
\end{remark}

In view of Remark~\ref{rem:denom-even}, it is unlikely that
there exists a simple formula relating the
denominator vectors $\dd_{\ell;t}^{B^0;t_0}$ and
$\dd_{\ell;t}^{B^1;t_1}$
(for $t_0 \overunder{k}{} t_1$ and $\mu_k(B^0) = B^1$).
This feature is expected to be rectified by passing from the
denominator vectors to the $\gg$-vector parametrization.
Introducing the latter will require a little preparation.

\medskip

Let~$\Acal$ be a cluster algebra of geometric type associated to an
$m \times n$ matrix $\tilde B_{t_0} = \tilde B^0$
(see~\eqref{eq:yt-geometric}--\eqref{eq:tilde-B}).
We assume that~$\Acal$ satisfies the following property:
\begin{equation}
\label{eq:full-rank}
\text{The matrix $\tilde B^0$ has full rank~$n$.}
\end{equation}
As shown in \cite[Lemma~1.2]{gsv1} and \cite[Lemma~3.2]{ca3}, the same is then true for all
matrices~$\tilde B_t$.
For instance, all the algebras with principal coefficients satisfy \eqref{eq:full-rank}.
The property \eqref{eq:full-rank} easily implies that, for every~$t$, the elements
\begin{equation}
\label{eq:yhat-geometric}
\hat y_{j;t} = y_{j;t} \prod_{i=1}^n x_{i;t}^{b^t_{ij}} =
\prod_{i=1}^m x_{i;t}^{b^t_{ij}} \quad (j \in [1,n])
\end{equation}
(cf.~\eqref{eq:yhat}) are algebraically independent in~$\Fcal$;
here the elements $x_{i;t} = x_i$ for $n < i \leq m$ are the
generators of the coefficient semifield~$\PP = \Trop(x_{n+1}, \dots, x_m)$.

Let $\Mcal$ denote the set of all nonzero
elements $z \in \Acal$ that can be written in the form
\begin{equation}
\label{eq:yhat-presentation}
z = R(\hat y_{1;t}, \dots, \hat y_{n;t}) \
\prod_{i=1}^m x_{i;t}^{a_i}
\end{equation}
for some vertex $t \in \TT_n$,
where~$R$ is a rational function in~$n$ variables with
coefficients in~$\QQ$, and $a_i \in \ZZ$ for $i \in [1,m]$.
In view of the formulas \eqref{eq:y-mutation-detailed} and
\eqref{eq:exchange-rel-rewrite}, if~$z$ is of the form
\eqref{eq:yhat-presentation} for some~$t$, then it can be written
in the same form for \emph{any} choice of~$t$.
Therefore (or by Corollary~\ref{cor:xjt=F/F}),
each cluster variable---hence each cluster monomial---belongs
to~$\Mcal$.
Clearly, $\Mcal$ is a subgroup of
the multiplicative group of the ambient field~$\Fcal$.

We call a rational
function~$R(u_1, \dots, u_n) \in \QQ(u_1, \dots, u_n)$
\emph{primitive} if it can be written as a ratio of two 
polynomials not divisible by any~$u_i$.
In view of \eqref{eq:yhat-geometric}, any $z \in \Mcal$ can
be written in the form \eqref{eq:yhat-presentation} with~$R$ primitive.
We claim that such a presentation is unique.

\begin{proposition}
\label{pr:yhat-presentation-unique}
If~$\Acal$ satisfies \eqref{eq:full-rank}, then every $z \in \Mcal$
has a unique presentation in the form \eqref{eq:yhat-presentation} with~$R$ primitive.
\end{proposition}

\begin{proof}
It is enough to show that if
\begin{equation}
\label{eq:two-presentations}
P(\hat y_{1;t}, \dots, \hat y_{n;t}) =
Q(\hat y_{1;t}, \dots, \hat y_{n;t}) \
\prod_{i=1}^m x_{i;t}^{a_i}
\end{equation}
for some primitive polynomials~$P$ and~$Q$ and some integers
$a_1, \dots, a_m$, then $a_i = 0$ for all~$i$, and so $P = Q$.
Since each $\hat y_{j;t}$ is a Laurent monomial in~$m$
independent variables $x_{1;t}, \dots, x_{m;t}$, the relation
\eqref{eq:two-presentations} implies that
$\prod_{i=1}^m x_{i;t}^{a_i}$ can be expressed as a Laurent
monomial in 
$\hat y_{1;t}, \dots, \hat y_{n;t}$.
Since $P$ and $Q$ are primitive, their ratio can be a Laurent monomial
only if they are equal to each other, and we are done.
\end{proof}

Now everything is in place for defining the $\gg$-vector
parametrizations (of~$\Mcal$).

\begin{definition}
\label{def:g-vectors}
For any $z \in \Mcal$ and any $t \in \TT_n$, the
\emph{$\gg$-vector} of~$z$ with respect to~$t$ is the vector $\gg_t(z) \in \ZZ^n$
defined as follows: if~$z$ is expressed in the form
\eqref{eq:yhat-presentation} with~$R$ primitive, then we set
$\gg_t(z) = \smallcolumn{a_1}{a_n}$.
(This is well defined by
Proposition~\ref{pr:yhat-presentation-unique}.)
\end{definition}

Definition~\ref{def:g-vectors} implies at once the following
analogue of \eqref{eq:delta-multiplicative}:
\begin{equation}
\label{eq:g-multiplicative}
\gg_t(z_1 z_2) = \gg_t(z_1) + \gg_t(z_2) \,.
\end{equation}

By Corollary~\ref{cor:xjt=F/F} and Proposition~\ref{pr:F-primitive},
we have
\begin{equation}
\label{eq:g-cluster-variable}
\gg_{t_0}(x_{\ell;t}) = \gg_{\ell;t}^{B^0;t_0},
\end{equation}
as defined by \eqref{eq:degree-vector};
so our terminology is consistent with that of
Section~\ref{sec:principal-gradings}.

In analogy with \eqref{eq:denominator-monomial}, we
use \eqref{eq:g-multiplicative} to extend the notation
\eqref{eq:g-cluster-variable} from cluster variables to
cluster monomials by setting
\begin{equation}
\label{eq:g-monomial}
\gg_{\aa;t}^{B^0;t_0} = \gg_{t_0}(x_{\aa;t})
= \sum_{\ell = 1}^n a_\ell \ \gg_{\ell;t}^{B^0;t_0} \ .
\end{equation}
In particular, 
we have (cf.~\eqref{eq:denom-cluster-monomial-initial})
\begin{equation}
\label{eq:g-cluster-monomial-initial}
\gg_{\aa;t_0}^{B^0;t_0} = \aa \quad (\aa \in \ZZ^n_{\geq 0}) \ .
\end{equation}

We suggest the following counterpart of
Conjecture~\ref{con:denominator-2}.

\begin{conjecture}
\label{con:g-vector-2}
{\ }
\begin{enumerate}
\item
Different cluster monomials have different $\gg$-vectors with
respect to a given initial seed.
\item
For every~$t \in \TT_n$, the vectors
$\gg_{1;t}^{B^0;t_0}, \dots, \gg_{n;t}^{B^0;t_0}$
form a $\ZZ$-basis of the lattice~$\ZZ^n$.
\end{enumerate}
\end{conjecture}

\begin{remark}
\label{rem:g-injectivity-implies-lin-indep}
It is not hard to see that Part~1 of Conjecture~\ref{con:g-vector-2}
together with Conjecture~\ref{con:conjectures-on-denoms} imply
Conjecture~\ref{con:cl-mon-lin-indep} (in a cluster algebra
satisfying \eqref{eq:full-rank}).
Indeed, assuming Conjecture~\ref{con:conjectures-on-denoms} and
applying Corollary~\ref{cor:xjt=F/F}, we see that the Laurent
expansion of any cluster monomial~$z$ in terms of the elements
$x_i = x_{i;t_0} \,\, (i \in [1,m])$ contains a Laurent
monomial $x_1^{a_1} \cdots x_m^{a_m}$ such that
$\aa=\smallcolumn{a_1}{a_n} = \gg_{t_0}(z)$, and the exponent vector of
any other Laurent monomial in the same expansion is obtained from
$\aa$ by adding some nonnegative integer linear
combination of the columns of~$\tilde B^0$.
Since the columns of~$\tilde B^0$ are assumed to be linearly
independent, the monomial $x_1^{a_1} \cdots x_m^{a_m}$ is the leading term
of~$z$ with respect to an appropriate term order.
Then Part~1 of Conjecture~\ref{con:g-vector-2} implies that
different cluster monomials have different leading terms, and so
are linearly independent.
\end{remark}

The property in Part 2 of Conjecture~\ref{con:g-vector-2}
is stronger than the corresponding (conjectural) property of the denominator
vectors (cf.~Remark~\ref{rem:denom-even}).
Another (conjectural) advantage of the $\gg$-vectors is the
following nice transition rule.

\begin{conjecture}
\label{con:g-transition}
Let $t_0 \overunder{k}{} t_1$ be two adjacent vertices in~$\TT_n$,
and let $B^1 = \mu_k(B^0)$.
Then, for any $t \in \TT_n$ and $\aa \in \ZZ_{\geq 0}^n$,
the $\gg$-vectors $\gg_{\aa;t}^{B^0;t_0} = (g_1, \dots, g_n)$
and $\gg_{\aa;t}^{B^1;t_1} = (g'_1, \dots, g'_n)$ are related as
follows:
\begin{equation}
\label{eq:Langlands-dual-trop}
g'_j =
\begin{cases}
-g_k  & \text{if $j = k$};\\[.05in]
g_j + [b^0_{jk}]_+ g_k
  - b^0_{jk} \min(g_k,0)
 & \text{if $j \neq k$}.
\end{cases}
\end{equation}
\end{conjecture}

\begin{remark}
\label{rem:conjectures-interplay-1}
Conjecture~\ref{con:g-transition} is easily seen to be a
consequence of Conjectures~\ref{con:hk-gk} and
\ref{con:signs-gi}.
Indeed, combining \eqref{eq:g-transition} with
\eqref{eq:hk-gk-conjecture} (and interchanging~$t_0$ and~$t_1$),
we conclude that \eqref{eq:Langlands-dual-trop} holds for the
$\gg$-vectors of cluster variables; the fact that \eqref{eq:Langlands-dual-trop}
extends to cluster monomials follows from Conjecture~\ref{con:signs-gi}.
\end{remark}

\pagebreak[3]

\begin{remark}
\label{rem:conjectures-interplay-2}
If $\gg_1, \dots, \gg_n$ are sign-coherent
vectors forming a $\ZZ$-basis in~$\ZZ^n$, then the transformation
\eqref{eq:Langlands-dual-trop} sends them to  a $\ZZ$-basis in $\ZZ^n$.
Thus, the validity of Part~2 of Conjecture~\ref{con:g-vector-2} for
\emph{all} initial seeds would follow from its validity for any single
choice of~$t_0$ in conjunction with Conjectures~\ref{con:signs-gi}
and~\ref{con:g-transition}.
\end{remark}

\begin{remark}
\label{rem:Fock-Goncharov}
Let $\ZZ_{\min, +}$ denote the set of integers supplied with the following
semifield structure: the multiplication in $\ZZ_{\min, +}$ is the
usual addition of integers, while the ``addition" is given by
$a \oplus b = \min(a,b)$.
Comparing \eqref{eq:Langlands-dual-trop} with~\eqref{eq:y-mutation},
we see that Conjecture~\ref{con:g-transition}
can be restated as follows: for every cluster monomial~$z$, the
correspondence $t \mapsto (\gg_t(z), B_t^T)$ (here $B_t^T$ denotes
the transpose of~$B_t$) is a $Y$-pattern (see
Definition~\ref{def:patterns}) with values in $\ZZ_{\min, +}\,$.
The fact that the exchange matrices $B_t$ get replaced by
their transposes brings to mind the \emph{Langlands duality}.
According to Part 1 of Conjecture~\ref{con:g-vector-2}, each
cluster monomial gives rise to a certain ``Langlands dual
tropical $Y$-pattern.''
This conjecture is consistent with the series of conjectures
by V.~Fock and A.~Goncharov in  \cite[Section~5]{fg2}
on a ``canonical'' basis in a cluster algebra.
In fact, our construction of $\gg$-vectors was motivated by a
desire to understand \cite[Conjecture~5.1]{fg2} and make the
parameterizations introduced there more explicit.
\end{remark}

We conclude this section by a proposition and a conjecture
relating denominator vectors of cluster variables to the
corresponding $\gg$-vectors.

\pagebreak[3]

\begin{proposition}
\label{pr:denom-gg}
Assume that condition~\eqref{eq:full-rank} holds.
Then the denominator vector
$\dd_{\ell;t}^{B^0;t_0} = \smallcolumn{d_1}{d_n}$
and the $\gg$-vector \linebreak[3]
$\gg_{\ell;t}^{B^0;t_0} = \smallcolumn{g_1}{g_n}$
of the same cluster
variable with respect to the same initial seed satisfy the
following relation:
\begin{equation}
\label{eq:dk+gk-thru-F}
u_1^{-d_1-g_1} \cdots u_n^{-d_n-g_n} =
F^{B^0;t_0}_{\ell,t}|_{\Trop(u_1, \dots, u_n)}
\Bigl(\prod_i u_i^{b^0_{i1}}, \dots, \prod_i u_i^{b^0_{in}}\Bigr) \,.
\end{equation}
\end{proposition}

\begin{proof}
Remembering the definition \eqref{eq:Laurent-normal-form} of the
denominator vector, we see that \eqref{eq:dk+gk-thru-F} is a
consequence of \eqref{eq:xjt=F/F}.
\end{proof}

In view of Proposition~\ref{pr:denom-gg},
Conjecture~\ref{con:g-thru-same-F} can be equivalently
restated as follows.

\begin{conjecture}
\label{con:denominators-thru-F}
Suppose that $B^0, t_0, \ell$ and $t$ are such that the
cluster variable $x_{\ell;t}$ does not belong to the cluster at~$t_0$.
Then the denominator vector
$\dd_{\ell;t}^{B^0;t_0} = \smallcolumn{d_1}{d_n}$ is given by
\begin{equation}
\label{eq:denom-thru-F}
u_1^{-d_1} \cdots u_n^{-d_n} =
F^{B^0;t_0}_{\ell,t}|_{\Trop(u_1, \dots, u_n)}
(u_1^{-1},\dots,u_n^{-1}) \, .
\end{equation}
\end{conjecture}

\begin{remark}
\label{rem:denom-leading-in-F}
Comparing Conjecture~\ref{con:denominators-thru-F}
with Conjecture~\ref{con:conjectures-on-denoms-2}, we see that
their conjunction can be stated as follows:
the vector $[\dd_{\ell;t}^{B^0;t_0}]_+$ appears as the vector
of exponents in the unique monomial
in the polynomial $F_{\ell;t}^{B^0;t_0}(y_1, \dots, y_n)$
that is divisible by all other occurring monomials.
\end{remark}

\pagebreak[3]

\section{Bipartite belt and $Y$-systems}
\label{sec:bipartite}


\begin{definition}
\label{def:bipartite-seed}
We call a (labeled) seed $(\xx, \yy, B)$ (and its exchange matrix
$B=(b_{ij})$) \emph{bipartite} if there is a function
$\varepsilon: [1,n] \to \{1, -1\}$ such that, for all $i$ and~$j$,
\begin{equation}
\label{eq:bipartite-B}
b_{ij} > 0 \Longrightarrow
\begin{cases}
\varepsilon(i) = 1 \,, \\
\varepsilon(j) = -1\,.
\end{cases}
\end{equation}
\end{definition}

For the rest of this paper, we are working in a cluster algebra
$\Acal=\Acal(\xx_0, \yy_0, B)$ of rank~$n$ with a bipartite initial exchange
matrix $B=(b_{ij})$.

The \emph{Cartan counterpart} of~$B$ \cite[(1.6)]{ca2}
is the (generalized) \emph{Cartan matrix} $A = A(B)=(a_{ij})$ defined by
\begin{equation}
\label{eq:assoc-cartan}
a_{ij} =
\begin{cases}
2 & \text{if $i=j$;} \\ 
- |b_{ij}| & \text{if $i\neq j$.}
\end{cases}
\end{equation}
The \emph{Coxeter graph} of~$A$ has the vertices $1, \dots, n$,
with~$i$ and~$j$ joined by an edge whenever $a_{ij} < 0$.
A Cartan matrix~$A$ is the Cartan counterpart of
a  bipartite matrix $B$ if and only if
the Coxeter graph of $A$ is bipartite,
which explains the terminology.


Recall that $\mu_k$ denotes the mutation in direction~$k$;
we use the same notation for seeds, $Y$-seeds, and exchange matrices.
Note that $\mu_k \mu_\ell = \mu_\ell \mu_k$ whenever
$b_{k \ell} = 0$.
Therefore the following two ``composite"
mutations are well defined on bipartite seeds:
\begin{equation}
\label{eq:mu-pm}
\mu_+ = \prod_{\varepsilon(k)= 1} \mu_k \,,\qquad
\mu_- = \prod_{\varepsilon(k)=-1} \mu_k \,.
\end{equation}
Furthermore, each of $\mu_+$ and $\mu_-$ is involutive and
transforms bipartite seeds to bipartite ones; in fact, $\mu_\pm (B)=-B$.
In the terminology of \cite[Section~9.1]{ca2},
these mutation transformations are ``shape-preserving.''

\begin{definition}
[\emph{Bipartite belt}]
\label{def:bipartite-belt}
The initial bipartite seed $\Sigma_0=(\xx_0, \yy_0, B)$
is naturally included into a \emph{bipartite belt} consisting of the seeds
\[
\Sigma_m=(\xx_{m}, \yy_{m}, (-1)^m B) \qquad (m \in \ZZ)
\]
defined by setting, for each $r > 0$,
\begin{align}
\label{eq:big-circle-seeds+}
\Sigma_r=(\xx_{r}, \yy_{r}, (-1)^r B) &=
\underbrace{\mu_\pm \cdots \mu_- \mu_+ \mu_-}_{r \text{~factors}}
(\Sigma_0 
),\\
\label{eq:big-circle-seeds-}
\Sigma_{-r}=(\xx_{- r}, \yy_{- r}, (-1)^r B) &=
\underbrace{\mu_\mp \cdots \mu_+ \mu_- \mu_+}_{r \text{~factors}}
(\Sigma_0 
).
\end{align}
We write the cluster $\xx_m$ and the coefficient tuple $\yy_m$ as
\begin{equation}
\label{eq:xm-ym}
\xx_m = (x_{1;m}, \dots, x_{n;m}), \quad \yy_m = (y_{1;m}, \dots, y_{n;m})\ .
\end{equation}
\end{definition}

The mutation rules of Definition~\ref{def:seed-mutation} imply that,
for $i,j\in [1,n]$ and $m\in\ZZ$,
\begin{align}
\label{eq:x-parity}
x_{i;m+1} = x_{i;m}
&\quad \text{if}\quad \varepsilon(i) = (-1)^m \,,\\ 
\label{eq:y-parity}
y_{j;m+1} = y_{j;m}^{-1}
&\quad \text{if}\quad \varepsilon(j) = (-1)^{m-1} \,.
\end{align}
Hence in studying the elements $x_{i;m}$ and $y_{j;m}$,
we may concentrate on the families
\begin{align}
\label{eq:xim}
&\{
x_{i;m} \,:\, \varepsilon(i) = (-1)^m
\}\,,
\\
\label{eq:yjm-Y-system}
&\{
y_{j;m} \,:\, \varepsilon(j) = (-1)^{m-1}
\}\,.
\end{align}

The following exchange relations are easy consequences of
\eqref{eq:y-mutation}--\eqref{eq:exchange-rel-xx}.
First, for  $\varepsilon(i) = (-1)^m$, we have
\begin{equation}
\label{eq:ym-exchange}
y_{i;m-1}\, y_{i;m+1} = \prod_{\varepsilon(j) = -\varepsilon(i)}
(y_{j;m} \oplus 1)^{-a_{ji}}
\,.
\end{equation}
Second, for $\varepsilon(j) = (-1)^{m-1}$, we have
\begin{equation}
\label{eq:xm-exchange}
x_{j;m-1}\, x_{j;m+1} = \frac{y_{j;m} \prod_{\varepsilon(i) = -\varepsilon(j)}
x_{i;m}^{-a_{ij}} + 1}{y_{j;m} \,\oplus\,1}
\,.
\end{equation}

\begin{example}
[\emph{Type~$A_2$}]
\label{example:A2-xy-bipartite}
Let
\begin{equation}
\label{eq:A2-epsilon}
B=\begin{bmatrix}
0 & 1\\
-1 & 0
\end{bmatrix}
\,,\qquad
A=A(B)=\begin{bmatrix}
2 & -1\\
-1 & 2
\end{bmatrix}
\,,\qquad
\varepsilon(1)=1
\,,\quad
\varepsilon(2)=-1.
\end{equation}
(This choice of signs is consistent with
\cite[Section~3.5]{yga}.)
Then $\mu_+=\mu_1\,$, $\mu_-=\mu_2\,$, and
\eqref{eq:big-circle-seeds+}--\eqref{eq:big-circle-seeds-} become
\begin{align*}
\cdots\cdots\cdots\cdots&\cdots\cdots\cdots\cdots\cdots\\
(\xx_{-2},\yy_{-2},B)&=\mu_2(\xx_{-1},\yy_{-1},-B),\\
(\xx_{-1},\yy_{-1},-B)&=\mu_1(\xx_0,\yy_0,B),\\
(\xx_1,\yy_1,-B)&=\mu_2(\xx_0,\yy_0,B),\\
(\xx_2,\yy_2,B)&=\mu_1(\xx_1,\yy_1,-B),\\
\cdots\cdots\cdots\cdots&\cdots\cdots\cdots\cdots\cdots
\end{align*}
in agreement with the notation
used in Example~\ref{example:A2-pattern}.
The families \eqref{eq:xim} and \eqref{eq:yjm-Y-system}
take the form
$\{x_{i;m} \,:\, i\not\equiv m\bmod 2\}$
and
$\{y_{j;m} \,:\, j\equiv m\bmod 2\}$,
respectively.
These elements are
shown in Table~\ref{table:bipartite-A2},
obtained by erasing half of the entries from
Table~\ref{table:seeds-A2} and substituting
$y_1=u_1^{-1}$ and $y_2=u_2\,$, so that the $y_{j;m}$
are expressed in terms of
\[
u_1=y_{1;-1}=y_1^{-1}
\quad \text{and} \quad
u_2=y_{2;0}=y_2
\]
(cf.\ \cite[Example~1.3]{yga}), whereas
the $x_{i;m}$ are expressed in terms of $u_1$, $u_2$,
\[
x_1=x_{1;0}
\quad \text{and} \quad
x_2=x_{2;-1}=x_{2;0}\,.
\]
\end{example}

\begin{table}[ht]
\begin{equation*}
\begin{array}{|c|cc|cc|cc|}
\hline
&&&&&&\\[-4mm]
m & y_{1;m} & y_{2;m} & x_{1;m} & x_{2;m} & \alpha(1;m) & \alpha(2;m) \\[1mm]
\hline
&&&&&&\\[-3mm]
-1 &
u_1
&
&
& x_2 & & -\alpha_2
\\[2.5mm]
\hline
&&&&&&\\[-3mm]
0 &
&u_2&x_1& &-\alpha_1 &
\\[2.5mm]
\hline
&&&&&&\\[-3mm]
1 &
\dfrac{u_2\oplus 1}{u_1} &
& & \dfrac{x_1u_2+1}{x_2(u_2\oplus 1)}
& & \alpha_2
\\[4.5mm]
\hline
&&&&&&\\[-3mm]
2 &
&
\dfrac{u_1\oplus u_2\oplus 1}{u_1 u_2}&
\dfrac{x_1u_2+x_2u_1+1}{(u_1\oplus u_2\oplus 1)x_1x_2} &
& \alpha_1+\alpha_2 &
\\[4.5mm]
\hline
&&&&&&\\[-3mm]
3 &
\dfrac{u_1\oplus 1}{u_2}&
&
& \dfrac{x_2 u_1+1}{x_1(u_1\oplus 1)}
& & \alpha_1
\\[4.5mm]
\hline
&&&&&&\\[-3mm]
4 &
& u_1 &
x_2 &
&-\alpha_2 &
\\[2.5mm]
\hline
&&&&&&\\[-3mm]
5 & u_2 &
&
&x_1 & & 
\\[2.5mm]
\hline
\end{array}
\end{equation*}
\caption{Bipartite dynamics in type~$A_2$}
\label{table:bipartite-A2}
\end{table}


Formula \eqref{eq:ym-exchange} can be viewed as a recurrence relation
for the family \eqref{eq:yjm-Y-system} of elements $y_{j;m}\,$.
This recurrence is a natural generalization of Zamolodchikov's
\hbox{$Y$-systems} \cite{yga,zamolodchikov},
so we
refer to it as
a (generalized) \emph{$Y$-system} (with values in a semifield~$\PP$)
associated with a (generalized) Cartan matrix~$A$.

There are two natural choices for the initial data in a $Y$-system.
The first one utilized in \cite{yga} takes an arbitrary $n$-tuple
$(u_1, \dots, u_n) \in \PP^n$ given by
\begin{equation}
\label{eq:yjm-Y-system-initial}
u_i =
\begin{cases}
y_{i;-1} & \text{if $\varepsilon(i) = +1$;} \\[.05in]
y_{i;0} & \text{if $\varepsilon(i) = -1$.}
\end{cases}
\end{equation}
An alternative choice is to use as the initial settings
an arbitrary $n$-tuple
$\yy_0=(y_1,\dots,y_n)\in\PP^n$.
In view of \eqref{eq:y-parity}, the two choices are related by
\begin{equation}
\label{eq:yjm-Y-system-initial-y0}
u_i =
\begin{cases}
y_i^{-1} & \text{if $\varepsilon(i) = +1$;} \\[.05in]
y_i & \text{if $\varepsilon(i) = -1$,}
\end{cases}
\end{equation}
so passing from one to another does not present any problems.

While dealing with $Y$-systems, we may (and will) work in 
the universal semifield $\Qsf(y_1, \dots, y_n)$ of
Definition~\ref{def:semifield-universal}.  
Using $\yy_0=(y_1,\dots,y_n)$ as the initial data and
remembering the notation introduced in
Definition~\ref{def:Y-functions}, we can write the 
solution of a $Y$-system as
\begin{equation}
\label{eq:Y-solution}
y_{j;m} = Y_{j;t_m}^{B;t_0}(y_1, \dots, y_n) \, ,
\end{equation}
where $t_0 \in \TT_n$ (resp., $t_m \in \TT_n$)
is the vertex to which the labeled seed
$\Sigma_0$ (resp.,~$\Sigma_m$) is attached.
(For a general semifield~$\PP$, one interprets \eqref{eq:Y-solution} as the
corresponding specialization of the rational functions
$Y_{j;t_m}^{B;t_0}$.) 

\begin{remark}
\label{rem:A-transpose}
In a traditional $Y$-system \cite{yga,zamolodchikov},
$A$~is a Cartan matrix of finite type,
the initial data have form~\eqref{eq:yjm-Y-system-initial},
and $\PP=\Qsf(u_1,\dots,u_n)$.
Note that there is a notational difference between the $Y$-system
relation \eqref{eq:ym-exchange} and its
counterpart \cite[(1.1)]{yga}: these two relations feature the exponents $-a_{ji}$
and $-a_{ij}$, respectively.
Accordingly, while translating the results in \cite{yga} into our
current notation (or vice versa), the Cartan matrix~$A$ should be
transposed.
\end{remark}

\pagebreak[3]

\begin{theorem}
[{Laurent phenomenon for generalized $Y$-systems}]
\label{th:Laurent-Y-system}
In the $Y$-system \eqref{eq:ym-exchange} associated with an
arbitrary symmetrizable generalized Cartan matrix~$A$, every element
$y_{j;m}$ with $\varepsilon(j) = (-1)^{m-1}$ viewed as a rational
function in the initial data $u_1, \dots, u_n$ (or $y_1, \dots,
y_n$) is in fact an integer Laurent polynomial.
\end{theorem}


Theorem~\ref{th:Laurent-Y-system} generalizes the finite type
result \cite[Theorem~1.5]{yga};
note however that we do not obtain the positivity property
established in \emph{loc.~cit.}

\begin{proof}
The theorem follows from 
\eqref{eq:Y-solution} and
Proposition~\ref{pr:laurent-Y-sink}
once we realize that for $Y_{j;t} = y_{j;m}$ with $\varepsilon(j) = (-1)^{m-1}$,
all the exponents $b_{ij}$ in \eqref{eq:Y-F-1} are nonnegative.
\end{proof}

\begin{example}
[\emph{$Y$-system of rank~$2$; cf.\ \cite[(6.9)]{ca1}}]
For $n=2$, a 
(generalized) Cartan matrix is of the form
\[
A=\begin{bmatrix}
2 & -b\\
-c & 2
\end{bmatrix}\,,
\]
where the integers $b$ and~$c$ are either both positive or both equal
to~$0$.
Setting~$\varepsilon(1)=1$ and $\varepsilon(2)=-1$,
and assuming that $\oplus$ is the ordinary addition,
the $Y$-system 
becomes:
\begin{align*}
y_{2;m+1}=
\dfrac{(y_{1;m} + 1)^b}{y_{2;m-1}} & \quad \text{if $m$ is
  odd;}\\[.1in]
y_{1;m+1}=
\dfrac{(y_{2;m} + 1)^c}{y_{1;m-1}} & \quad \text{if $m$ is even.}
\end{align*}
Using the initial data $u_1 = y_{1;-1}, \, u_2 = y_{2;0}$
(cf.~\eqref{eq:yjm-Y-system-initial}),
we obtain:
\[
y_{1;1}=\dfrac{(u_2 + 1)^c}{u_1}\,,\ \
y_{2;2}=\dfrac{\bigl({\frac{(u_2 + 1)^c}{u_1} +
    1\bigr)}^b}{u_2}\,,\ \
y_{1;3}=\biggl(\dfrac{\bigl({\frac{(u_2 + 1)^c}{u_1} +
    1\bigr)}^b}{u_2} + 1\biggr)^c \dfrac{u_1}{(u_2 + 1)^c}\,,
\ \text{etc.}
\]
To illustrate Theorem~\ref{th:Laurent-Y-system}: to see that
$y_{1;3}$ is indeed a Laurent polynomial in $u_1$ and $u_2$, we can rewrite it as
$$y_{1;3} = \dfrac{1}{u_1^{bc-1} u_2^c}
\biggl(u_1^b + \dfrac{\bigl((u_2+1)^c + u_1\bigr)^b - u_1^b}{u_2+1}
\biggr)^c \, .$$
\end{example}

\begin{example}
[\emph{$Y$-system of type~$A_{2r-1}^{(1)}$}]
Let $n$ be an even integer $\geq 4$, and let~$A$ be the Cartan
matrix of affine type $A_{n-1}^{(1)}$.
So the Coxeter graph of~$A$ is a $n$-cycle, and $a_{ij} = -1$ for any two adjacent
vertices $i$ and~$j$.
The corresponding $Y$-system can be described in simple terms as follows.
Place initial variables $u_1, \dots, u_n$ around a circle.
Alternately color them red and blue.
Replace each red variable~$u$ by $\frac{(v+1)(w+1)}{u}$, where
$v$ and $w$ are the two neighbors of~$u$.
Then do the same for all blue variables, then again for all the red ones, etc.
Theorem~\ref{th:Laurent-Y-system} asserts that all rational functions
obtained this way are integer Laurent polynomials in $u_1, \dots, u_n$.
In particular, if the initial values are all specialized
to~$1$, then all the numbers generated in this $Y$-system become integers.
We leave it to the reader to check that these integers are squares of Fibonacci
numbers.
\end{example}


Without losing ground, we can assume that the matrix~$A$
is \emph{indecomposable}, that is,
its Coxeter graph is connected.
In most situations,
the general case reduces to this one in an obvious way.

Al.~B.~Zamolodchikov conjectured~\cite{zamolodchikov} and we
proved~\cite[Theorems~1.1--1.2]{yga} that if~$A$ is an indecomposable
Cartan matrix of finite type, then the corresponding \hbox{$Y$-system}
exhibits the following periodicity property:
\begin{equation}
\label{eq:y-period}
y_{j;m+2(h+2)}=y_{j;m} \quad (m \in \ZZ, \,\, j \in [1,n]) \, ,
\end{equation}
where $h$ is the \emph{Coxeter number}.
(In the special case of type~$A_n$, this was proved earlier
in~\cite{frenkel-szenes, gliozzi-tateo}.)
Our next result sharpens \eqref{eq:y-period} as follows.

\begin{theorem}
\label{th:belt-periodic}
Suppose an exchange matrix~$B$ is bipartite, and its Cartan
counterpart $A = A(B)$ is indecomposable.
\begin{enumerate}
\item
If~$A$ is of finite type, then the corresponding
bipartite belt (see Definition~\ref{def:bipartite-belt}) has
the following periodicity property: the labeled seeds $\Sigma_m$
and $\Sigma_{m+ 2(h+2)}$ are equal to each other for all~$m \in \ZZ$.
\item
If~$A$ is of infinite type, then all
the elements $x_{i;m}$ in \eqref{eq:xim}
are distinct, as are
all the elements $y_{j;m}$ in~\eqref{eq:yjm-Y-system} viewed as
Laurent polynomials in the initial data.
\end{enumerate}
\end{theorem}

To rephrase Theorem~\ref{th:belt-periodic}:
the bipartite belt (or the associated $Y$-system)
is periodic if and only if the Cartan matrix~$A$ is of finite type;
if it is periodic, then $2(h+2)$ is a period;
if it is not, then all cluster variables on the bipartite belt
are distinct, as are all Laurent polynomials in the associated
$Y$-system.

Theorem~\ref{th:belt-periodic} is proved in Section~\ref{sec:denom-bipartite}.

\section{Bipartite dynamics on roots}
\label{sec:prelim-roots}

As before, let $A=(a_{ij})$ be a (symmetrizable generalized) $n \times n$ Cartan
matrix with a bipartite Coxeter graph.
We identify the lattice $\ZZ^n$ with the root lattice~$Q$
associated with~$A$.
Specifically, the standard basis vectors in~$\ZZ^n$ are denoted by
$\alpha_1, \dots, \alpha_n$ and identified with the simple roots.
Let $W \subset GL(Q)$ be the corresponding Weyl group; it is generated by
the \emph{simple reflections} $s_1, \dots, s_n$ which act on the
simple roots by
\begin{equation}
\label{eq:si-action}
s_i (\alpha_j) = \alpha_j - a_{ij} \alpha_i.
\end{equation}
Following \cite{yga}, we define the elements $t_\pm \in W$ by setting
\begin{equation}
\label{eq:Tpm}
t_+ = \prod_{\varepsilon(k)=1} s_k\,,\qquad
t_- = \prod_{\varepsilon(k)=-1} s_k\,.
\end{equation}
As in \eqref{eq:mu-pm},
these two elements are well defined because the factors in each
of the two products mutually commute.
By the same reason, each of them is an involution: $t_+^2 = t_-^2 = 1$.
To illustrate \eqref{eq:si-action} and~\eqref{eq:Tpm}, the actions of
$t_+$ and $t_-$ on the simple roots are given by
\begin{equation}
\label{eq:t+-action}
t_+ (\alpha_j) =
\begin{cases}
- \alpha_j & \text{if $\varepsilon(j) = +1$;} \\[.05in]
\alpha_j - \sum_{i \neq j} a_{ij} \alpha_i
 & \text{if $\varepsilon(j) = -1$;}
\end{cases}
\end{equation}
\begin{equation}
\label{eq:t--action}
t_- (\alpha_j) =
\begin{cases}
\alpha_j - \sum_{i \neq j} a_{ij} \alpha_i
 & \text{if $\varepsilon(j) = +1$;} \\[.05in]
- \alpha_j & \text{if $\varepsilon(j) = -1$.}
\end{cases}
\end{equation}

\begin{definition}
\label{def:alpha-im}
We define the vectors $\alpha(i;m) \in Q$,
for all $i \in [1,n]$ and all $m \in \ZZ$ such that $\varepsilon(i) = (-1)^m$,
by setting, for all $r \geq 0$:
\begin{align}
\label{eq:dems-m-positive}
 \alpha(i;r) &=
\underbrace{t_- t_+ \cdots t_{\varepsilon(i)}}_{r
  \text{~factors}}(-\alpha_i) \qquad\, \text{for $\varepsilon(i) = (-1)^r$;}\\
\label{eq:dems-m-negative}
 \alpha(j;-r-1) &=
\underbrace{t_+ t_- \cdots t_{\varepsilon(j)}}_{r
\text{~factors}}(-\alpha_j) \qquad \text{for $\varepsilon(j) = (-1)^{r-1}$.}
\end{align}
In particular, we have
$\alpha(i;m) = - \alpha_i$ for $m \in \{0, -1\}$.
(Note that this notation is slightly different from the one in
\cite{yga}.)
\end{definition}

\begin{example}
[{\cite[Example~2.3]{yga}}]
\label{example:A2-t+t-}
In type $A_2$ (cf.\
Example~\ref{example:A2-xy-bipartite}), the elements $t_+=s_1$
and $t_-=s_2$ act as follows:
\begin{equation}
\label{eq:A2-t+t-}
\begin{array}{ccc}
-\alpha_1 & \stackrel{\textstyle t_+}{\verylongleftrightarrow}
~\alpha_1~ \stackrel{\textstyle t_-}{\verylongleftrightarrow}
~\alpha_1\!+\!\alpha_2~
\stackrel{\textstyle t_+}{\verylongleftrightarrow} ~\alpha_2~
\stackrel{\textstyle t_-}{\verylongleftrightarrow}
& -\alpha_2\,.
\end{array}
\end{equation}
The corresponding vectors $\alpha(i;m)$ with $m\in [-1,4]$ are shown
in Table~\ref{table:bipartite-A2} .
\end{example}

Clearly, all the vectors $\alpha(i;m)$ are (real) roots; in
particular, for each of them, all nonzero components are of the
same sign.
To determine these signs, we assume that the Cartan
matrix~$A$ is indecomposable; the general case easily reduces to
this one.

\pagebreak[3]

\begin{proposition}
\label{pr:positive-roots-1}
Let~$A$ be an indecomposable Cartan matrix of finite type (i.e.,
one of the Cartan-Killing types $A_n, \dots, G_2$).
Let $h$ be the corresponding Coxeter number, that is, the order of
the Coxeter element $t_+ t_-$ in~$W$.
Then $\alpha(i;m)$ is a positive root for any
$m \in [-h-1,h] - \{-1, 0\}$ and any $i$ with $\varepsilon(i) = (-1)^m$.
All these positive roots are distinct, and every
positive root associated to~$A$ appears among them.
Furthermore, we have:
\begin{align}
\label{eq:roots-periodicity-1}
\alpha(i;-h-2) &= -\alpha_{i^*} \quad \text{for $\varepsilon(i) =
  (-1)^h$,}\\
\label{eq:roots-periodicity-2}
\alpha(j;h+1)  &= -\alpha_{j^*} \quad \text{for $\varepsilon(j) =
  (-1)^{h-1}$,}
\end{align}
where $i \mapsto i^*$ is the involution induced by the longest
element $w_0 \!\in\! W$: $w_0(\alpha_i)\!=\! - \alpha_{i^*}$.
\end{proposition}

\begin{proof}
This is an immediate consequence of a classical
result of R.~Steinberg \cite{steinberg} (cf.~\cite[Lemma~2.1,
Proposition~2.5]{yga}).
\end{proof}

\begin{theorem}
\label{th:positive-roots-2}
Let~$A$ be an indecomposable Cartan matrix of infinite type
whose Coxeter graph is bipartite.
Then the vectors $\alpha(i;m)$ (see Definition~\ref{def:alpha-im})
for $m \in \ZZ - \{-1, 0\}$ and $\varepsilon(i) = (-1)^m$
are distinct positive roots.
\end{theorem}


\begin{proof}
Consider the matrix $A' = 2 - A$.
Since $A'$ is symmetrizable, all of its eigenvalues are real.
Let $\rho$ denote the maximal eigenvalue of~$A'$.
As in~\cite{blm}, we observe that
$\rho \geq 2$.
Indeed, a Cartan matrix $A$ is of finite type if and only if its
symmetrization is positive definite,
i.e., all eigenvalues of~$A$ are positive or, equivalently,
all eigenvalues of $A'$ are less than~$2$.
(As shown in \cite{blm}, $A$ is of affine type if and only if
$\rho = 2$, but we will not need this.)

We next note that all entries of $A'$ are nonnegative.
Moreover, since $A$ is indecomposable,
some power of $A'$ has all entries positive.
Hence the Perron-Frobenius theorem applies to~$A'$,
implying in particular that
an eigenvector for the eigenvalue~$\rho$ of the transpose of $A'$ can
be chosen to have all components positive.
Explicitly, there are positive real numbers
$c_1, \dots, c_n$ such that
\begin{equation}
\label{eq:eigenvector}
- \sum_{i \neq j} a_{ij} c_i = \rho c_j
\quad
\text{for any~$j \in [1,n]$.}
\end{equation}

Let $\RR^2$ be a $2$-dimensional real vector space with a fixed  basis
$\{\ee_+, \ee_-\}$ and two linear involutions $s_+$ and $s_-$ acting
in this basis by the matrices
\begin{equation}
\label{s-pm}
s_+ = \bmat{-1}{\rho}{0}{1}, \quad
s_- = \bmat{1}{0}{\rho}{-1} \ .
\end{equation}
Let $U: Q \to \RR^2$ be the linear map defined by
$U(\alpha_i) = c_i \,\ee_{\varepsilon(i)}$ for all~$i$, where
$(c_1, \dots, c_n)\in\RR_{> 0}^n$
satisfies~\eqref{eq:eigenvector}.

\begin{lemma}
\label{lem:Ut=sU}
\ $U t_{\pm} = s_{\pm} U$.
\end{lemma}

\begin{proof}
The two cases are completely analogous, so we only treat one.
This is a simple calculation that uses \eqref{eq:t+-action},
\eqref{eq:eigenvector}, and~\eqref{s-pm} (in this order):
\begin{align*}
Ut_+(\alpha_j)
&=
\begin{cases}
-c_j\,\ee_{+} & \text{if $\varepsilon(j) = +1$;} \\[.05in]
c_j\,\ee_{-} - (\sum_{i \neq j} a_{ij} c_i)\,\ee_{+}
 & \text{if $\varepsilon(j) = -1$;}
\end{cases} \\
&=
\begin{cases}
-c_j\,\ee_{+} & \text{if $\varepsilon(j) = +1$;} \\[.05in]
c_j\,\ee_{-} +\rho\, c_j\,\ee_{+}
 & \text{if $\varepsilon(j) = -1$;}
\end{cases} \\
&=c_j\,s_+(\ee_{\varepsilon(j)})\\
&=s_+\,U(\alpha_j)\,,
\end{align*}
as desired.
\end{proof}

Note that~$U$ sends positive (resp., negative) roots to positive (resp.,
negative) linear combinations of $\ee_+$ and~$\ee_-\,$.
To show that each $\alpha(i;m)$ with $m \in \ZZ - \{-1, 0\}$
and $\varepsilon(i) = (-1)^m$ is a positive root, it therefore suffices
to prove that $U(\alpha(i;m))$ is a positive linear combination
of $\ee_+$ and~$\ee_-\,$.
By Lemma~\ref{lem:Ut=sU}, 
\[
U(\alpha(i;m))
=U\, t_\pm t_\mp \cdots t_{\varepsilon(i)} (-\alpha_i)
=c_i\,s_\pm s_\mp\cdots s_{-\varepsilon(i)} (\ee_{\varepsilon(i)}).
\]
Thus,
we need to show the following elementary statement:
if $s_+$ and $s_-$ are given by \eqref{s-pm} with $\rho \geq 2$, then
any vector of the form $s_\pm s_\mp \cdots s_\varepsilon (\ee_{-\varepsilon})$ is
a positive linear combination of $\ee_+$ and $\ee_-$.
This statement is hardly new, and easy to prove.
Direct inspection shows that the vectors in
question are of the form $(U_m(\rho), U_{m+1}(\rho))$ and
$(U_{m+1}(\rho), U_{m}(\rho))$ for $m>0$,
where the sequence $(U_m(\rho))$ is defined by the linear recurrence
\[
U_0(\rho)=0, \quad
U_1(\rho)=1, \quad
U_{m+1}(\rho)=\rho\,U_m(\rho)-U_{m-1}(\rho)
\quad (m > 0) \ .
\]
Assuming that $U_m(\rho) > U_{m-1}(\rho) \geq 0$ for some $m > 0$
and using the fact that $\rho \geq 2$, we conclude that
$U_{m+1}(\rho)-U_m(\rho)\geq U_m(\rho)-U_{m-1}(\rho) > 0$, implying
by induction on $m$ that
$U_{m}(\rho)>0$ for all $m>0$, as desired.
%
%

We have shown that all the roots in \eqref{eq:dems-m-positive}
and \eqref{eq:dems-m-negative} with $r > 0$, i.e., the roots
$\alpha(i;r)$ and $\alpha(j;-r-1)$ with $r > 0$, are positive.
To finish the proof of Theorem~\ref{th:positive-roots-2}, it
remains to show that all these roots are distinct.
Thus, we need to show that every equality of the form
\begin{equation}
\label{eq:aim=ajn}
t_\pm t_\mp \cdots t_{\varepsilon(i)}(-\alpha_i) =
t_\pm t_\mp \cdots t_{\varepsilon(j)}(-\alpha_j)
\end{equation}
implies that $i = j$ and there are the same number of factors on
both sides.
Let~$w$ (resp.,~$w'$) denote the product
$t_\pm t_\mp \cdots t_{\varepsilon(i)}$ (resp.,
$t_\pm t_\mp \cdots t_{\varepsilon(j)}$) on the left (resp., on the
right) side in \eqref{eq:aim=ajn}.
Rewriting \eqref{eq:aim=ajn} as
$$- \alpha_i = w^{-1} w' (-\alpha_j),$$
we note that, unless the products $w$ and $w'$ are identical, the
root $w^{-1} w' (-\alpha_j)$ is positive as shown above, and so
cannot be equal to $-\alpha_i$, finishing the proof.
\end{proof}

Theorem~\ref{th:positive-roots-2} can be rephrased as follows.

\begin{corollary}
\label{cor:coxeter-elt-inf-type}
Let~$A$ be an $n\times n$ indecomposable
Cartan matrix of infinite type whose Coxeter graph is bipartite.
Then $\ell((t_+ t_-)^r) = rn$ for all $r > 0$, where $\ell(w)$ is
the length of an element $w$ in the Coxeter group~$W$.
Equivalently, the concatenation of any number of copies of
$\prod_{\varepsilon(k)=1} s_k \prod_{\varepsilon(k)=-1} s_k$
is a reduced decomposition in~$W$.
\end{corollary}

\begin{proof}
The equivalence of Corollary~\ref{cor:coxeter-elt-inf-type} and
Theorem~\ref{th:positive-roots-2} is clear from the following
well-known characterization of
reduced decompositions: a decomposition $s_{i_1} \cdots s_{i_\ell}$
is reduced if and only if all the roots $s_{i_1} \cdots
s_{i_{k-1}}(\alpha_{i_k})$ are positive for $k = 1, \dots, \ell$.
\end{proof}

\begin{remark}
\label{rem:positive-roots-2}
Corollary~\ref{cor:coxeter-elt-inf-type}
sharpens the following result in~\cite{blm}:
if~$A$ is an indecomposable Cartan matrix of infinite type
with a bipartite Coxeter graph,
then the element $t_+ t_- \in W$ is of infinite order.
(In fact~\cite{howlett}, any Coxeter element in any infinite
Coxeter group has infinite order.)
The above proof of Theorem~\ref{th:positive-roots-2}
is based on similar ideas as the proof in \cite{blm}
but seems to be simpler
although we prove a stronger statement.
\end{remark}

\begin{remark}
It is conceivable that Corollary~\ref{cor:coxeter-elt-inf-type}
could be extended to an arbitrary infinite Coxeter group~$W$
(not necessarily crystallographic) with a bipartite Coxeter graph.
We note that in this generality, the Coxeter elements $t_+t_-$ were
studied in~\cite[Section~5]{mcmullen}.
\end{remark}

\pagebreak[3]

\section{Denominators and $\gg$-vectors on the bipartite belt}
\label{sec:denom-bipartite}

We stick to the notation of the preceding two sections.
In this section we compute, for any cluster variable
$x_{i;m}$ lying on the bipartite belt,
the denominator vector and the $\gg$-vector of $x_{i;m}$
with respect to the initial cluster~$\xx_0$.

We begin by obtaining an expression for the denominator vector
 of $x_{i;m}$
(see \eqref{eq:Laurent-normal-form}--\eqref{eq:denominator-vector})
that is similar to the one in
Definition~\ref{def:alpha-im}.

Let $\Phireal_{\geq -1}$ denote the union of the set of
real positive roots and the set of negative simple roots.
Following \cite{yga}, we introduce the involutive permutations
$\tau_+$ and $\tau_-$ of $\Phireal_{\geq -1}$ by setting
\begin{equation}
\label{eq:tau-action}
\tau_\varepsilon (\alpha) =
\begin{cases}
\alpha & \text{if $\alpha = -\alpha_j$ with $\varepsilon(j) = -\varepsilon$;} \\[.05in]
t_\varepsilon (\alpha)
 & \text{otherwise}
\end{cases}
\end{equation}
(see \eqref{eq:Tpm}).

\begin{example}
\label{example:A2-tau-tropical}
In type $A_2$ (cf.\
Example~\ref{example:A2-t+t-}), the actions of $\tau_+$
and $\tau_-$ on $\Phireal_{\geq -1}$ are given by
\begin{equation}
\label{eq:A2-tau-tropical}
\begin{array}{ccc}
-\alpha_1 & \stackrel{\textstyle\tau_+}{\verylongleftrightarrow}
~\alpha_1~ \stackrel{\textstyle\tau_-}{\verylongleftrightarrow}
~\alpha_1\!+\!\alpha_2~
\stackrel{\textstyle\tau_+}{\verylongleftrightarrow} ~\alpha_2~
\stackrel{\textstyle\tau_-}{\verylongleftrightarrow}
& -\alpha_2\,. \\
\circlearrowright & & \circlearrowright \\ \tau_- & & \tau_+
\end{array}
\end{equation}
\end{example}

\begin{definition}[cf.~Definition~\ref{def:alpha-im}]
\label{def:dd-im}
We define the vectors $\dd(i;m) \in Q$,
for all $i \in [1,n]$ and all $m \in \ZZ$ such that $\varepsilon(i) = (-1)^m$,
by setting, for all $r \geq 0$:
\begin{align}
\label{eq:dems-m-positive-tau}
 \dd(i;r) &=
\underbrace{\tau_- \tau_+ \cdots \tau_{\varepsilon(i)}}_{r
  \text{~factors}}(-\alpha_i) \qquad\, \text{for $\varepsilon(i) = (-1)^r$;}\\
\label{eq:dems-m-negative-tau}
 \dd(j;-r-1) &=
\underbrace{\tau_+ \tau_- \cdots \tau_{\varepsilon(j)}}_{r
\text{~factors}}(-\alpha_j) \qquad \text{for $\varepsilon(j) = (-1)^{r-1}$.}
\end{align}
In particular,
$\dd(i;m) = - \alpha_i$ for $m \in \{0, -1\}$.
\end{definition}

\begin{theorem}
\label{th:denominators-bipartite}
The denominator vector for a cluster variable
$x_{i;m}$ with respect~to the initial cluster~$\xx_0$
is equal to~$\dd(i;m)$,
for any $i \!\in\! [1,n]$ and $m\! \in\! \ZZ$ with $\varepsilon(i)\! =\! (-1)^m$.
\end{theorem}

\begin{proof}
It suffices to treat the case when the Cartan matrix~$A$ is indecomposable.
We start by clarifying the relationship between the vectors
$\dd(i;m)$ and $\alpha(i;m)$.
First suppose that~$A$ is of finite type.
Comparing Definitions~\ref{def:alpha-im} and~\ref{def:dd-im} and
using Proposition~\ref{pr:positive-roots-1}, we conclude that in
this case
\begin{equation}
\label{eq:dd-alpha-finite}
\dd(i;m) = \alpha(i;m) \quad (m \in [-h-2, h+1]),
\end{equation}
and
\begin{equation}
\label{eq:dd-finite-periodic}
\dd(i;m + 2(h + 2)) = \dd(i;m) \quad (m \in \ZZ).
\end{equation}
If $A$ is of infinite type then by Theorem~\ref{th:positive-roots-2}
\begin{equation}
\label{eq:dd-alpha-infinite}
\dd(i;m) = \alpha(i;m) \quad (m \in \ZZ).
\end{equation}

Turning to the proof of Theorem~\ref{th:denominators-bipartite},
note that it holds for $m \in \{-1, 0\}$ since in
this case, for $\varepsilon(i) = (-1)^m$,
the denominator vector of $x_{i;m}$ is equal to $\dd(i;m) = - \alpha_i$.
Applying \eqref{eq:exchange-denominator} to the exchange relations
\eqref{eq:xm-exchange}, we see that to finish the proof it
suffices to show that the vectors $\dd(i;m)$ with $\varepsilon(i) = (-1)^m$
satisfy the recurrence relation
\begin{equation}
\label{eq:denominator-recurrence}
\dd(j;m-1) + \dd(j;m+1) =
\Bigl[- \sum_{\varepsilon(i) = - \varepsilon(j)} a_{ij} \dd(i;m)\Bigr]_+
\end{equation}
(recall that the operation $a \mapsto [a]_+$ on vectors is
understood component-wise).

First suppose that $A$ is of infinite type, and so
$\dd(i;m) = \alpha(i;m)$ for all $m \in \ZZ$.
The equality \eqref{eq:denominator-recurrence} holds for $m = 0$
since $\dd(j;1) = \alpha(j;1) = \alpha_j$ for $\varepsilon(j) = -1$.
For $m > 0$, all roots $\alpha(i;m)$ are positive by
Theorem~\ref{th:positive-roots-2}, hence \eqref{eq:denominator-recurrence}
takes the form
\begin{equation}
\label{eq:alpha-recurrence}
\alpha(j;m-1) + \alpha(j;m+1) =
- \sum_{\varepsilon(i) = - \varepsilon(j)} a_{ij} \alpha(i;m) \ .
\end{equation}
The equality \eqref{eq:alpha-recurrence} is a direct consequence
of \eqref{eq:dems-m-positive}:
\begin{align*}
- \alpha(j;m-1) -
\sum_{\varepsilon(i) = (-1)^m} a_{ij} \,\alpha(i;m)
&= \underbrace{t_- t_+ t_- \cdots\,}_{m-1 \text{~factors}}\,
   (\alpha_j - \sum_{\varepsilon(i)= (-1)^m} a_{ij} \,\alpha_i)
\\
&= \underbrace{t_- t_+ t_- \cdots\,}_{m \text{~factors}}\,(\alpha_j)\\
&= \alpha(j;m+1).
\end{align*}
The case $m < 0$ is treated in the same way.

If $A$ is of finite type, we use \eqref{eq:dd-alpha-finite}, and
then the same argument as above proves \eqref{eq:denominator-recurrence}
for $m \in [-h-1, h]$.
The cases $m = -h-2$ and $m = h+1$ follow from
\eqref{eq:roots-periodicity-1} and \eqref{eq:roots-periodicity-2};
for $m$ outside the interval $[-h-2, h+1]$, use the periodicity
property~\eqref{eq:dd-finite-periodic}.
\end{proof}

\begin{example}
To illustrate Theorem~\ref{th:denominators-bipartite},
the denominator vectors for the elements $x_{i;m}$ shown in
Table~\ref{table:bipartite-A2} appear in the last column of that
table.
In type~$A_2$, we have $h=3$, so for all the values shown in the last
column, \eqref{eq:dd-alpha-finite} applies.
\end{example}

As a first corollary of Theorem~\ref{th:denominators-bipartite},
we show that Conjecture~\ref{con:sign-coherence-denom} holds for
the clusters on the bipartite belt.

\begin{corollary}
\label{cor:denominator-signs}
For any $m \in \ZZ$, the denominator vectors
of the cluster variables
\[
\{x_{j;m-1} \,:\, \varepsilon(j) = (-1)^{m-1}\}
\cup
\{x_{i;m}   \,:\, \varepsilon(i) = (-1)^{m}  \}.
\]
are sign-coherent (see Definition~\ref{def:coherent-vectors}).
\end{corollary}

\begin{proof}
Combine formulas \eqref{eq:dd-alpha-finite}--\eqref{eq:dd-alpha-infinite}
with Proposition~\ref{pr:positive-roots-1} and
Theorem~\ref{th:positive-roots-2}.
\end{proof}

The following result generalizes its finite-type counterpart obtained
in \cite[Theorem~1.9]{ca2}.

\begin{corollary}
\label{cor:nonzero-const-term}
Each cluster variable $x_{i;m}$ can be written as
\[
x_{i;m}=\frac{P_{i;m}(\xx_0)}{\xx_0^{\dd(i;m)}}\,,
\]
where $P_{i;m}$ is a polynomial with nonzero constant term.
\end{corollary}

\begin{proof}
The fact that $x_{i;m}$ is a Laurent polynomial whose numerator has
nonzero constant term is proved by induction on~$m$ employing
\cite[Lemma~5.2]{ca2}.
The denominator is then given by Theorem~\ref{th:denominators-bipartite}.
\end{proof}

As another application of Theorem~\ref{th:denominators-bipartite},
we compute the tropical evaluation of the coefficients $y_{j;m}$.
As in \eqref{eq:Y-solution}, we use
$\yy_0=(y_1,\dots,y_n)$ as the initial data.
We will use the shorthand
\[
\yy_0^\dd = \prod_{k \in [1,n]} y_{k}^{d_k}
\qquad \text{for}\qquad
\dd = \column{d_1}{d_n} \in \ZZ^n
\,.
\]

\begin{proposition}
\label{pr:exchange-bipartite-principal}
For all $m \in \ZZ$ and $j\in[1,n]$ with $\varepsilon(j) = (-1)^{m-1}$,
the elements $y_{j;m}$ evaluated in the tropical semifield
$\PP = \Trop(y_{1}, \dots, y_{n})$ are given by
\begin{equation}
\label{eq:ym-tropical}
y_{j;m}|_{\Trop(y_{1}, \dots, y_{n})} = \yy_0^{-\dd(j;m-1)}
\end{equation}
(see Definition~\ref{def:dd-im}).
\end{proposition}

\begin{proof}
We first verify \eqref{eq:ym-tropical}
by inspection for $m \in \{-1, 0\}$.
Assuming \eqref{eq:ym-tropical} for some $m \in \ZZ$,
we easily see from the definition of the tropical semifield
and the sign coherence property of Corollary~\ref{cor:denominator-signs}
that the right-hand side of \eqref{eq:ym-exchange} takes the form
$$\prod_{\varepsilon(j) = (-1)^{m-1}}
(y_{j;m} \oplus 1)^{-a_{ji}}
=
\Bigl(\prod_{\varepsilon(j) = (-1)^{m-1}} y_{j;m}^{-a_{ji}}\Bigr) \oplus 1 \ .$$
As a result, the $Y$-system dynamics
\eqref{eq:ym-exchange} becomes identical to the cluster dynamics
\eqref{eq:xm-exchange} (with trivial coefficients) after
identifying $y_{j;m}$ with $x_{j;m-1}\,$.
Thus \eqref{eq:ym-tropical} can be seen as a special case of
Theorem~\ref{th:denominators-bipartite}.
\end{proof}

\begin{example}
In type~$A_2$, the monomials
$y_{j;m}|_{\Trop(y_{1},y_{2})}$
already appeared in Table~\ref{table:seeds-geom-A2}.
Here we are only interested in the values $y_{j;m}$ with
$j\equiv m\bmod 2$. 
See the second column of Table~\ref{table:bipartite-A2-d-g}.
The corresponding vectors of exponents match the entries
located one row higher
in the fourth column of Table~\ref{table:bipartite-A2-d-g}
(or, equivalently, the last column of Table~\ref{table:bipartite-A2}).
\end{example}

\begin{table}[ht]
\begin{equation*}
\begin{array}{|c|cc|cc|cc|cc|}
\hline
&&&&&&&&\\[-4mm]
m & y_{1;m} & y_{2;m} & x_{1;m} & x_{2;m} & \dd(1;m) & \dd(2;m) &
\gg(1;m) & \gg(2;m) \\[1mm]
\hline
&&&&&&&&\\[-3mm]
-1 &
\frac{1}{y_1}
& & & x_2 & & -\alpha_2 & & \alpha_2
\\[2.5mm]
\hline
&&&&&&&&\\[-3mm]
0 &
&y_2&x_1& &-\alpha_1 & &\alpha_1 &
\\[2.5mm]
\hline
&&&&&&&&\\[-3mm]
1 &
y_1 &
& & \frac{x_1y_2+1}{x_2}
& & \alpha_2 && -\alpha_2
\\[2.5mm]
\hline
&&&&&&&&\\[-3mm]
2 &
&
\frac{1}{y_2}&
\frac{x_1y_1y_2+y_1+x_2}{x_1x_2} &
& \alpha_1+\alpha_2 & &-\alpha_1  &
\\[2.5mm]
\hline
&&&&&&&&\\[-3mm]
3 &
\frac{1}{y_1 y_2}&
&
& \frac{y_1+x_2}{x_1}
& & \alpha_1 && -\alpha_1+\alpha_2
\\[2.5mm]
\hline
&&&&&&&&\\[-3mm]
4 &
& \frac{1}{y_1} &
x_2 &
&-\alpha_2 & & \alpha_2 &
\\[2.5mm]
\hline
&&&&&&&&\\[-3mm]
5 & y_2 &
&
&x_1 & & -\alpha_1 & & \alpha_1
\\[2.5mm]
\hline
\end{array}
\end{equation*}
\smallskip

\caption{
Bipartite dynamics with principal coefficients in type~$A_2$}
\label{table:bipartite-A2-d-g}
\end{table}

Proposition~\ref{pr:exchange-bipartite-principal} enables us
to explicitly write the exchange relations~\eqref{eq:xm-exchange}
in the cluster algebra $\Aprin(B)$ with principal coefficients
at the initial bipartite seed $(\xx_0,\yy_0,B)$
(see Definitions~\ref{def:principal-coeffs}
and~\ref{def:bipartite-seed})
in terms of the denominators~$\dd(i;m)$.

\begin{corollary}
\label{cor:exchanges-on-big-circle}
In~$\Aprin(B)$, the exchange relations \eqref{eq:xm-exchange} take the form
\begin{equation}
\label{eq:xm-exchange-bipartite-principal}
x_{j;m-1}\, x_{j;m+1} = \yy_0^{[-\dd(j;m-1)]_+} \prod_{\varepsilon(i) = -\varepsilon(j)}
x_{i;m}^{-a_{ij}} \ + \ \yy_0^{[\dd(j;m-1)]_+} \,,
\end{equation}
for $m \in \ZZ$ and $\varepsilon(j) = (-1)^{m-1}$.
\end{corollary}

\begin{corollary}
\label{cor:conjectures-bipartite}
Conjectures~\ref{con:conjectures-on-denoms},
\ref{con:conjectures-on-denoms-2}, and
\ref{con:denominators-thru-F}
(hence Conjecture~\ref{con:g-thru-same-F})
hold in the case when clusters and cluster
variables in question belong to the bipartite belt.
Thus, a polynomial~$F^{B^0;t_0}_{\ell;t_m}\,$
has constant term~$1$, and among its monomials
there is one, namely
$\yy_0^{[\dd(\ell;m)]_+}$,
which is divisible by all others.
\end{corollary}

\begin{proof}
Conjectures~\ref{con:conjectures-on-denoms}
and \ref{con:conjectures-on-denoms-2}
follow from the observation that
the formula~\eqref{eq:xm-exchange-bipartite-principal}
satisfies condition~{\rm (iii)} of
Proposition~\ref{pr:equivalent-conjectures-on-denoms}.
It remains to show \eqref{eq:denom-thru-F}.
This is also an easy consequence of
\eqref{eq:xm-exchange-bipartite-principal}, which implies that
the denominator vectors $\dd(\ell;m)$ and the tropical evaluations in
 \eqref{eq:denom-thru-F} satisfy the same recurrence relations
 (we leave the details to the reader).
\end{proof}

\begin{remark}
\label{rem:Y-denominators}
Tropical evaluations of the elements $y_{j;m}$ in
Proposition~\ref{pr:exchange-bipartite-principal}
can be viewed as $Y$-system analogues of the denominator vectors of cluster
variables.
Indeed, by Theorem~\ref{th:Laurent-Y-system},
every element $y_{j;m}$ with $\varepsilon(j) = (-1)^{m-1}$
is an integer Laurent polynomial in the initial data $y_1, \dots, y_n$.
In analogy with \eqref{eq:Laurent-normal-form},
Proposition~\ref{pr:exchange-bipartite-principal} can be rephrased as
saying that this Laurent polynomial has the form
\[
y_{j;m} = \dfrac{N(y_{1}\,,\dots,y_{n})}{\yy_0^{\dd(j;m-1)}}\,,
\]
where~$N$ is an integer polynomial not divisible by any of its $n$ variables.
Thus, Proposition~\ref{pr:exchange-bipartite-principal} provides a ``denominator vector"
parameterization of the elements $y_{j;m}$ in a $Y$-system by the roots in $\Phireal_{\geq-1}$
contained in the $\langle\tau_+,\tau_-\rangle$-orbit of the negative
simple roots.
This parameterization is close to but \emph{different} from that
of the elements in a $Y$-system (by the same set of roots)
used in~\cite{yga} in the finite type situation.
The relationship between these two parameterizations will be discussed in
the next section.
\end{remark}

\begin{proof}[Proof of Theorem~\ref{th:belt-periodic}]
First suppose that an indecomposable Cartan matrix $A$ is of infinite type.
Combining \eqref{eq:dd-alpha-infinite} with Theorem~\ref{th:positive-roots-2},
we conclude that the vectors $\dd(i;m)$ for $m \in \ZZ$ and $\varepsilon(i) = (-1)^m$
are all distinct.
Therefore, the cluster variables $x_{i;m}$ (resp., the Laurent
polynomials $y_{i;m+1}$) are all distinct by Theorem~\ref{th:denominators-bipartite}
(resp., by Proposition~\ref{pr:exchange-bipartite-principal}),
proving Part (2) of Theorem~\ref{th:belt-periodic}.

Now suppose $A$ is of finite type. In view of \eqref{eq:y-period},
to prove Part (1) of Theorem~\ref{th:belt-periodic}, one only
needs to show that $x_{i;m} = x_{i;m + 2(h+2)}$ for all $i$ and
$m$ with $\varepsilon(i) = (-1)^m$. This follows from
Theorem~\ref{th:denominators-bipartite}, the periodicity property
\eqref{eq:dd-finite-periodic}, and the following result, which is
a consequence of \cite[Theorem~1.9]{ca2} (see also
Proposition~\ref{pr:fintype-belt} below): in a cluster algebra of
finite type, a cluster variable is uniquely determined by its
denominator vector with respect to a bipartite initial cluster
$\xx_0$.
\end{proof}

\begin{proof}[Proof of Proposition~\ref{pr:finite-eq-finite-mutation}]
The ``only if" part is immediate from the fact (established in \cite{ca2}
and already mentioned in Section~\ref{sec:exchange-graphs}) that the property of a seed $(\xx, \yy, B)$ to
define a cluster algebra of finite type depends only on the matrix~$B$.
To prove the ``if" part, we need to show the following:
\begin{enumerate}
\item[($*$)]
If $B$ gives rise to a cluster algebra of infinite type, then
there exists a matrix~$\tilde B$ with principal part~$B$
such that $\tilde B$ is of infinite mutation type.
\end{enumerate}
As shown in \cite{ca2}, if $B$ gives rise to a cluster algebra of
infinite type, then the same is true for some principal $2 \times 2$
submatrix $B'$ of a matrix obtained from~$B$ by a sequence of mutations.
It is then enough to show that ($*$) holds for~$B'$.
Since~$B'$ is bipartite, we can take the corresponding $\tilde B'$
as a $4 \times 2$ matrix giving rise to the cluster algebra with
principal coefficients; the fact that $\tilde B'$ is of infinite
mutation type is immediate from
Proposition~\ref{pr:exchange-bipartite-principal} and
Theorem~\ref{th:positive-roots-2}.
\end{proof}

Our next goal is to compute the $\gg$-vectors of the cluster
variables~$x_{i;m}$.
To state the result, we denote by $E$ the linear automorphism of
the root lattice~$Q$ acting on the basis of simple roots by
\begin{equation}
\label{eq:E(alpha)}
E(\alpha_i) = - \varepsilon(i)\, \alpha_i \qquad \text{for $i \in
  [1,n]$.}
\end{equation}

\begin{theorem}
\label{th:gg-im}
For every $i \in [1,n]$ and $m \in \ZZ$ with $\varepsilon(i) = (-1)^m$,
the $\gg$-vector $\gg(i;m)$ of the cluster variable
$x_{i;m}$ with respect to the initial seed
is given by
\begin{equation}
\label{eq:gim-dim}
\gg(i;m) = E \tau_- (\dd(i;m))
\end{equation}
(see \eqref{eq:E(alpha)}, \eqref{eq:tau-action} and
Definition~\ref{def:dd-im}).
\end{theorem}

\begin{proof}
The desired equality \eqref{eq:gim-dim} is easily
checked for $m \in \{-1, 0\}$.
Indeed, for $\varepsilon(j) = -1$, we have
\[
E \tau_- (\dd(j;-1))
= E \tau_- (-\alpha_j)
= E \alpha_j
= \alpha_j
= \gg(j;-1);
\]
and for $\varepsilon(i) = +1$, we have
\[
E \tau_- (\dd(i;0))
= E \tau_- (-\alpha_i)
= E(-\alpha_i)
= \alpha_i
= \gg(i;0).
\]

To prove \eqref{eq:gim-dim} for the rest of the values of~$m$,
recall \eqref{eq:g-cluster-variable} and~\eqref{eq:degree-vector}.
The vector
$\gg(i;m)$ is equal to the degree of $x_{i;m}$
in $\Aprin(B)$.
In view of Corollary~\ref{cor:exchanges-on-big-circle},
for any $m \in \ZZ$ and $\varepsilon(j) = (-1)^{m-1}$, we have
\begin{equation}
\label{eq:degree-bipartite-inductive-1}
\gg(j;m-1) +  \gg(j;m+1) =
{\rm deg}(\yy_0^{[\dd(j;m-1)]_+}) \ .
\end{equation}
Remembering the definition \eqref{eq:degrees-xy},
we rewrite the right-hand side of
\eqref{eq:degree-bipartite-inductive-1} as follows:
\begin{align}
\label{eq:gim-1}
\deg(\yy_0^{[\dd(j;m-1)]_+}) &=
\sum_k [\dd(j;m-1)_k]_+ \ {\rm deg}(y_{k})\\
\nonumber
& =
\sum_k [\dd(j;m-1)_k]_+ \sum_{i \neq k} \varepsilon(i) a_{ik}
\alpha_i \,,
\end{align}
where $\dd(j;m-1)_k$ is the $k$th component of the vector
$\dd(j;m-1)$.

To simplify \eqref{eq:gim-1} further, we follow \cite[(1.7)]{yga}
and extend the permutations $\tau_+$ and $\tau_-$ of
$\Phireal_{\geq -1}$ to piecewise-linear involutions of the
root lattice~$Q$ by setting
\begin{equation}
\label{eq:tau-pm-tropical-1}
\tau_\varepsilon (c_1 \alpha_1 + \cdots + c_n \alpha_n) =
c'_1 \alpha_1 + \cdots + c'_n \alpha_n,
\end{equation}
where
\begin{equation}
\label{eq:tau-pm-tropical-2}
c'_i =
\begin{cases}
- c_i - \sum_{j \neq i}  a_{ij} [c_j]_+
&
\text{if $\varepsilon (i) = \varepsilon$;} \\[.1in]
c_i & \text{otherwise.}
\end{cases}
\end{equation}
As an immediate consequence of
\eqref{eq:tau-pm-tropical-1}--\eqref{eq:tau-pm-tropical-2}, we
have, for any
$\alpha = c_1 \alpha_1 + \cdots + c_n \alpha_n \in Q$:
\begin{equation}
\label{eq:tau++tau-}
\tau_+ \alpha  + \tau_- \alpha =
-\sum_{i=1}^n \sum_{j\neq i} a_{ij}\,[c_j]_+\, \alpha_i \,.
\end{equation}

Using \eqref{eq:tau++tau-}, we rewrite
\eqref{eq:gim-1} as follows:
$${\rm deg}(\yy_0^{[\dd(j;m-1)]_+}) =
E(\tau_- \dd(j;m-1) + \tau_+ \dd(j;m-1)) \ .$$
In view of \eqref{eq:dems-m-positive-tau}--\eqref{eq:dems-m-negative-tau},
we have $\dd(j;m-1) = \tau_+ \tau_- \dd(j;m+1)$.
Hence \eqref{eq:degree-bipartite-inductive-1} becomes
\begin{equation*}
\label{eq:degree-bipartite-inductive-2}
\gg(j;m-1) +  \gg(j;m+1) =
E(\tau_- \dd(j;m-1) + \tau_- \dd(j;m+1)) \ .
\end{equation*}
It follows that each of the equalities
$\gg(j;m-1) = E \tau_- (\dd(j;m-1))$ and
$\gg(j;m+1) = E \tau_- (\dd(j;m+1))$ implies the other.
Since we already checked these equalities for two consecutive
values of~$m$, they are true for all $m \in \ZZ$, as desired.
\end{proof}

\begin{example}
In type~$A_2$ (with the conventions~\eqref{eq:A2-epsilon}),
the maps $\tau_\pm$ are given by
\begin{align*}
\tau_+(c_1\alpha_1+c_2\alpha_2)
&= (-c_1+[c_2]_+)\,\alpha_1 + c_2\alpha_2 \,,\\
\tau_-(c_1\alpha_1+c_2\alpha_2)
&= c_1\alpha_1 + (-c_2+[c_1]_+)\,\alpha_2 \,,
\end{align*}
in agreement with~\eqref{eq:A2-tau-tropical}.
The linear map~$E$ acts by
\[
E(\alpha_1)=-\alpha_1\,,\quad
E(\alpha_2)=\alpha_2\,.
\]
For the corresponding cluster algebra $\Aprin(B)$,
the elements $x_{i;m}$ and~$y_{j;m}$ are shown in
Table~\ref{table:bipartite-A2-d-g}.
These entries have been imported from Table~\ref{table:seeds-geom-A2}
(or set $u_1=y_1^{-1}$ and $u_2=y_2$ in Table~\ref{table:bipartite-A2}
and calculate in $\Trop(y_1,y_2)$).
The denominators $\dd(i;m)$ and degrees~$\gg(i;m)$ are copied from
Table~\ref{table:degrees-recursive}.
(Recall that, by convention,
$\alpha_1=\left[\begin{smallmatrix}1\\0\end{smallmatrix}\right]$,
$\alpha_2=\left[\begin{smallmatrix}0\\1\end{smallmatrix}\right]$.)
\end{example}

\begin{remark}
\label{rem:signs-gi-bipartite}
Theorem~\ref{th:gg-im} and Corollary~\ref{cor:denominator-signs}
imply that Conjecture~\ref{con:signs-gi} holds
in the case when the clusters and cluster
variables in question lie on the bipartite belt.
\end{remark}

By Theorems~\ref{th:denominators-bipartite} and
\ref{th:gg-im}, the denominator vector $\dd$ and
the $\gg$-vector $\gg$ of any cluster variable $x_{i;m}$ with respect to the initial seed
are related by $\gg = E \tau_- (\dd)$.
We now show that this relationship can be extended to all cluster
monomials on the bipartite belt.

\begin{corollary}
\label{cor:degrees-bipartite}
For every cluster monomial of the form
$x_{1;m}^{a_1} \cdots x_{n;m}^{a_n}$, its $\gg$-vector
$\gg$ and its denominator vector $\dd$ with respect to the initial seed
are related by
\begin{equation}
\label{eq:degree-dem-bipartite}
\gg = E \tau_- (\dd) \,.
\end{equation}
(Here $\tau_-$ is the piecewise linear involution
of~$Q$ given by
\eqref{eq:tau-pm-tropical-1}--\eqref{eq:tau-pm-tropical-2}.)
\end{corollary}

\begin{proof}
In view of \eqref{eq:delta-multiplicative} and
\eqref{eq:g-multiplicative}, it is enough to show
that
\[
\tau_- \Bigl(\sum_i a_i \ \dd(i;m)\Bigr) = \sum_i a_i \ \tau_-(\dd(i;m))
\]
for any nonnegative integers $a_1, \dots, a_n$.
This follows from the
definition
\eqref{eq:tau-pm-tropical-1}--\eqref{eq:tau-pm-tropical-2},
once we remember Corollary~\ref{cor:denominator-signs}.
\end{proof}

The last result in this section provides supporting evidence
for Conjectures~\ref{con:hk-gk} and~\ref{con:g-transition} in a bipartite setting.

\begin{proposition}
\label{pr:hk-gk-bipartite}
Let $\Sigma_0$ and
$\Sigma_1=\mu_\varepsilon(\Sigma_0)$ be
adjacent seeds on  the bipartite belt.
Then the $\gg$-vectors
\[
\gg(\ell;m)= \column{g_1}{g_n} \quad \text{and}\quad
\gg'(\ell;m)= \column{g'_1}{g'_n}
\]
of a cluster variable $x_{\ell;m}$ (with $\varepsilon(\ell)=(-1)^m$)
computed with respect to the seeds $\Sigma_0$ and~$\Sigma_1$,
respectively, are related by
\begin{equation}
\label{eq:g-transition-bipartite}
g_i =
\begin{cases}
-g'_i  & \text{if $\varepsilon(i) = \varepsilon$};\\[.05in]
g'_i + \displaystyle\sum_{\varepsilon(k)=\varepsilon}
([-b^0_{ik}]_+ g'_k - b^0_{ik} [-g_k']_+),
 & \text{if $\varepsilon(i) = -\varepsilon$}.
\end{cases}
\end{equation}
Here, as before, $(b^0_{ij})$ is the exchange matrix
at~$\Sigma_0\,$.
\end{proposition}

\begin{proof}
It follows from the definitions combined with Theorem~\ref{th:gg-im}
that
\begin{align}
\label{eq:gg=E...}
\gg(\ell;m) &= E\tau_- \dd(\ell;m), \\
\label{eq:gg'=-E}
\gg'(\ell;m) &= -E\tau_+ \dd(\ell;-m).
\end{align}
(If $B$ and $B'$ are the exchange
matrices at $\Sigma_0$ and $\Sigma_1$,
then $B' = -B$, hence the corresponding sign function
$\varepsilon': [1,n] \to \{1, -1\}$ is equal to $-\varepsilon$; this
is why in passing from \eqref{eq:gg=E...} to~\eqref{eq:gg'=-E},
$\tau_-$ gets replaced by $\tau_+$, and $E$ gets replaced by~$-E$.)
Definition~\ref{def:dd-im} implies that
$\dd(\ell;m)=\tau_-\dd(\ell;-m)$.
Consequently,
\[
\gg(\ell;m)
=E \dd(\ell;-m)
=E \tau_+ E (-\gg'(\ell;m)).
\]
Using \eqref{eq:E(alpha)}, \eqref{eq:tau-pm-tropical-2},
and~\eqref{eq:bipartite-B}, this can be checked to
match~\eqref{eq:g-transition-bipartite}.
\end{proof}

\pagebreak[3]

\section{Finite type}
\label{sec:finite-type}

In this section, we focus on the cluster algebras of finite type.
By the classification result of~\cite{ca2},
any cluster algebra of finite type has a bipartite seed,
so all the results in
Sections~\ref{sec:bipartite}--\ref{sec:denom-bipartite} apply.
In fact, a cluster algebra $\Acal=\Acal(\xx_0,\yy_0,B)$
with a bipartite initial matrix~$B$ is of finite type
if and only if the Cartan counterpart~$A=A(B)$ (see~\eqref{eq:assoc-cartan})
is a Cartan matrix of finite type.
(The ``only if'' direction follows from \cite[Propositions~9.3 and~9.7]{ca2}.)

The results in \cite{ca2, yga} imply that things are particularly nice
in finite type:
\pagebreak[3]

\begin{proposition}
\label{pr:fintype-belt}
In a cluster algebra of finite type,
\begin{enumerate}
\item[{\rm (1)}]
every cluster variable belongs to
a seed lying on the bipartite belt;
\item[{\rm (2)}]
the denominator vectors establish a bijection between cluster
variables and ``almost positive roots" (the elements of~$\Phi_{\geq
  -1}$),
and between cluster monomials and the root lattice.
\end{enumerate}
\end{proposition}

\begin{proof}
The first statement is a consequence of
Proposition~\ref{pr:positive-roots-1}
and \cite[Theorems 1.9 and~1.13]{ca2}.
To obtain the second statement, combine these theorems
with \cite[Theorems~1.8 and~1.10]{yga}.
\end{proof}

We now establish Conjecture~\ref{con:cl-mon-lin-indep}
for the cluster algebras of finite type.

\begin{theorem}
\label{th:cluster-monomials-linear-indep}
In a cluster algebra $\Acal$ of finite type,
the cluster monomials are linearly independent
over the ground ring~$\ZZ\PP$.
\end{theorem}

\begin{proof}
By Corollary~\ref{cor:nonzero-const-term}, each cluster monomial has
a unique ``leading'' term of smallest total degree with respect to the cluster
variables from the initial cluster.
Then Proposition~\ref{pr:fintype-belt} implies that
different cluster monomials have different leading
terms, hence are linearly independent.
%
%
\end{proof}

In \cite[Conjecture~4.28]{cdm}, we conjectured that
Theorem~\ref{th:cluster-monomials-linear-indep} can be strengthened as
follows:
in a cluster algebra of finite type, the cluster monomials  form an
additive basis.
For cluster algebras of infinite type, we expect this property to
fail.

\medskip

We next verify several conjectures from the previous sections
in the special case where $\Acal$ is a cluster algebra of finite type
and the initial seed $\Sigma_0=\Sigma_{t_0}$ is bipartite.

Conjecture~\ref{con:denominator-1}
(hence Conjecture~\ref{con:sign-coherence-denom},
already established in more general setting---see
Corollary~\ref{cor:denominator-signs})
is immediate from \cite[Theorems 1.9 and 1.13]{ca2}.

As an application of Proposition~\ref{pr:fintype-belt}, we extend
Corollary~\ref{cor:degrees-bipartite} to \emph{all} cluster
monomials (not necessarily associated with clusters lying on the bipartite belt).

\begin{proposition}
\label{pr:g-d-cluster-monomials-finite}
In a cluster algebra of finite type,
the $\gg$-vector $\gg$ and the denominator vector $\dd$
of any cluster monomial with respect to a bipartite initial seed
are related by \eqref{eq:degree-dem-bipartite},
where the piecewise linear involution $\tau_-$ and the linear
transformation~$E$ are given by
\eqref{eq:tau-pm-tropical-1}--\eqref{eq:tau-pm-tropical-2})
and~\eqref{eq:E(alpha)}, respectively.
\end{proposition}

\begin{proof}
In view of Proposition~~\ref{pr:fintype-belt} and
Theorem~\ref{th:gg-im}, the equality \eqref{eq:degree-dem-bipartite}
holds for any cluster variable.
It then extends to all cluster monomials by the same reasoning as
in Corollary~\ref{cor:degrees-bipartite}, using
Corollary~\ref{cor:denominator-signs}.
\end{proof}

We now turn to Conjectures~~\ref{con:denominator-2}
and~\ref{con:g-vector-2}.
In finite type, they can be strengthened as follows.

\begin{corollary}
\label{cor:fock-goncharov-finite-type-bipartite}
For a bipartite seed $\Sigma_{t_0}$ in a cluster algebra of
finite type, each of the maps
$x_{\aa;t} \mapsto \dd_{\aa;t}^{B;t_0}$ and
$x_{\aa;t} \mapsto \gg_{\aa;t}^{B;t_0}$
is a bijection between the set of all cluster monomials~and~$\ZZ^n$.
\end{corollary}

\begin{proof}
The statement about denominator vectors
follows from \cite[Theorems~1.9 and~1.13]{ca2}
and \cite[Theorems~1.8 and~1.10]{yga}.
It implies the statement about the $\gg$-vectors in view of
Proposition~\ref{pr:g-d-cluster-monomials-finite}.
\end{proof}




We next turn our attention to the
polynomials~$F^{B;t_0}_{\ell;t}\,$,
under the assumption that $B$ is of finite type,
and both seeds $t_0$ and~$t$ lie on the bipartite belt.
We will use the shorthand
\[
F_{\ell;m}=F^{B;t_0}_{\ell;t_m},
\]
where $t_m$ ($m\in\ZZ$) is a vertex on the bipartite belt with
an attached seed $\Sigma_m=\Sigma_{t_m}$
(see \eqref{eq:big-circle-seeds+}--\eqref{eq:big-circle-seeds-}).

We establish a link between the polynomials $F_{\ell;m}$
and the ``Fibonacci polynomials'' introduced and studied in~\cite{yga}.

\begin{definition}[\emph{Fibonacci polynomials}]
\label{def:fib-poly}
For $\varepsilon\in\{+,-\}$, let $E_\varepsilon$ denote the
automorphism of the field $\QQ(y_1,\dots,y_n)$ defined by
\begin{equation}
E_\varepsilon(y_i)=
\begin{cases}
y_i^{-1} & \text{if $\varepsilon(i)=\varepsilon$;}\\
y_i      & \text{otherwise.}
\end{cases}
\end{equation}
For $\ell\in[1,n]$ and $m\in\ZZ$, with $\varepsilon(\ell)=(-1)^m$, let
\begin{equation}
\label{eq:flm-def}
f_{\ell;m}\stackrel{\rm def}{=}
E_+(F_{\ell;m}(y_1,\dots,y_n))
\prod_{\varepsilon(j)=1} y_j^{[\dd(\ell;m)_j]_+},
\end{equation}
where $\dd(\ell;m)_j$ stands for the $j$th component of the
denominator vector $\dd(\ell;m)$.
By Corollary~\ref{cor:conjectures-bipartite}, $f_{\ell;m}$ is a polynomial in
$y_1,\dots,y_n\,$, which we call a \emph{Fibonacci polynomial}.
\end{definition}

In plain terms, $f_{\ell;m}$ is obtained from the polynomial
$F_{\ell;m}(y_1,\dots,y_n)$ by replacing each variable $y_j$
with $\varepsilon(j)=1$ by its inverse and clearing denominators.
It is then clear that, conversely,
\begin{equation}
\label{eq:Flm-via-flm}
F_{\ell;m}\stackrel{\rm def}{=}
E_+(f_{\ell;m}(y_1,\dots,y_n))
\prod_{\varepsilon(j)=1} y_j^{[\dd(\ell;m)_j]_+}.
\end{equation}

For the special case of type~$A_2$,
the polynomials $F_{\ell;m}$ and $f_{\ell;m}$ are shown in
Table~\ref{table:Fibonacci-A2}.

\begin{table}[ht]
\begin{equation*}
\begin{array}{|c|cc|cc|cc|}
\hline
&&&&&&\\[-4mm]
m & F_{1;m} & F_{2;m} & \dd(1;m) & \dd(2;m) &
f_{1;m} & f_{2;m} \\[1mm]
\hline
&&&&&&\\[-3mm]
0 & 1 & &-\alpha_1 & & 1 &
\\[2.5mm]
\hline
&&&&&&\\[-3mm]
1 & & y_2 + 1
& & \alpha_2 && y_2+1
\\[2.5mm]
\hline
&&&&&&\\[-3mm]
2 & y_1y_2+y_1+1 &
& \alpha_1+\alpha_2 & & y_1+y_2+1  &
\\[2.5mm]
\hline
&&&&&&\\[-3mm]
3 & & y_1+1
& & \alpha_1 && y_1+1
\\[2.5mm]
\hline
&&&&&&\\[-3mm]
4 & 1 &
&-\alpha_2 & & 1 &
\\[2.5mm]
\hline
&&&&&&\\[-3mm]
5 & & 1 & & -\alpha_1 & & 1
\\[2.5mm]
\hline
\end{array}
\end{equation*}
\smallskip

\caption{Fibonacci polynomials in type~$A_2$}
\label{table:Fibonacci-A2}
\end{table}

In \cite[Theorem~1.6]{yga}, we introduced a family of polynomials
\[
F[\alpha]\in\ZZ[y_1,\dots,y_n]
\]
labeled by the roots $\alpha\in\Phi_{\geq-1}\,$.
Combining \cite[Theorem~1.5, Theorem~1.7, and Remark~2.8]{yga}
and remembering Remark~\ref{rem:A-transpose} and
Definition~\ref{def:dd-im} of the current paper,
we see that the polynomials $F[\alpha]$ can be defined
as the unique family of polynomials with the following properties.
First, these polynomials satisfy
the initial conditions $F[-\alpha_i]=1$ (for all~$i$).
Second, for any $\alpha\!=\!\dd(j;-m-1)\in\Phi_{\geq -1}$
with $\varepsilon(j)\!=\!(-1)^{m-1}$, they satisfy the recurrence relation
\begin{equation}
\label{eq:F[alpha]-recur}
F[\tau_+\alpha]\, F[\tau_-\alpha]
=\yy_0^{[-\alpha]_+} \prod_{i\neq j} F[\dd(i;m)]^{-a_{ij}}
+\yy_0^{[ \alpha]_+}\,.
\end{equation}
(We note that $\tau_+\alpha=\dd(j;m-1)$ and
$\tau_-\alpha=\dd(j;m+1)$.)

\begin{theorem}
\label{th:fib=Fib}
In a cluster algebra of finite type,
$f_{\ell;m}=F[\alpha]$,
where $\alpha=\dd(\ell;m)$.
\end{theorem}

Thus, the terminology of Definition~\ref{def:fib-poly}
is consistent with that of~\cite{yga}:
the Fibonacci polynomials $f_{\ell;m}$ and $F[\alpha]$ are indeed the
same.

\begin{proof}
We need to verify that the polynomials $f_{\ell;m}$ satisfy the
recurrence~\eqref{eq:F[alpha]-recur}
(the initial conditions are immediate).
As before, we use the notation $\yy_0=(y_1,\dots,y_n)$.
Assuming that $\varepsilon(j)\!=\!(-1)^{m-1}$, we denote
\begin{align*}
\dd(j;m-1)&=(d_1,\dots,d_n),\\
\dd(j;m+1)&=(\tilde d_1,\dots,\tilde d_n).
\end{align*}
It then follows from
\eqref{eq:tau-pm-tropical-1}--\eqref{eq:tau-pm-tropical-2}
that
\[
\dd(j;-m-1)_\ell
=\begin{cases}
\tilde d_\ell & \text{if $\varepsilon(\ell)=1$;}\\
       d_\ell & \text{if $\varepsilon(\ell)=-1$.}
\end{cases}
\]
Thus, we need to show that
\begin{equation}
\label{eq:ff=...}
f_{j;m-1}\, f_{j;m+1}
=
\prod_{\varepsilon(\ell)=1} y_\ell^{[-\tilde d_\ell]_+}
\prod_{\varepsilon(k)=-1} y_k^{[-d_k]_+}
\prod_{i\neq j} f_{i;m}^{-a_{ij}}
+
\prod_{\varepsilon(\ell)=1} y_\ell^{[\tilde d_\ell]_+}
\prod_{\varepsilon(k)=-1} y_k^{[d_k]_+}
\,.
\end{equation}
This can be deduced from Corollary~\ref{cor:exchanges-on-big-circle}
as follows.
Combining the latter with~\eqref{eq:F-def}, we get
\[
F_{j;m-1}\, F_{j;m+1} = \prod_\ell y_\ell^{[-d_\ell]_+}
\prod_{\varepsilon(i) = -\varepsilon(j)}
F_{i;m}^{-a_{ij}} \ + \ \prod_\ell y_\ell^{[d_\ell]_+} \,.
\]
Applying $E_+$ and using~\eqref{eq:flm-def}, we obtain, after
straightforward simplifications:
\begin{align}
\label{eq:ff=pp...}
f_{j;m-1}\, f_{j;m+1} =
&\prod_{\varepsilon(\ell)=1} y_\ell^{d_\ell+[\tilde d_\ell]_+}
\prod_{\varepsilon(k)=-1} y_k^{[-d_k]_+}
\prod_{i\neq j} \Biggl(
f_{i;m}^{-a_{ij}} \prod_{\varepsilon(\ell)=1}  y_\ell^{a_{ij}
  [\dd(i;m)_\ell]_+} \Biggr) \\[.2in]
\nonumber
& +
\prod_{\varepsilon(\ell)=1} y_\ell^{[\tilde d_\ell]_+}
\prod_{\varepsilon(k)=-1} y_k^{[d_k]_+}
\,.
\end{align}
We next use~\eqref{eq:denominator-recurrence}, which
in view of Corollary~\ref{cor:denominator-signs} can be rewritten as
\begin{equation}
\label{denom-rec-bip}
d_\ell+\tilde d_\ell = -\sum_{i\neq j} a_{ij} [\dd(i;m)_\ell]_+\,.
\end{equation}
Putting \eqref{eq:ff=pp...} and \eqref{denom-rec-bip} together
yields~\eqref{eq:ff=...}.
\end{proof}

\begin{corollary}
Each polynomial $F_{\ell;m}(y_1, \dots, y_n)$ has
positive coefficients.
\end{corollary}

\begin{proof}
By \cite[Theorem~1.6]{yga},
each polynomial $F[\alpha]$ has positive coefficients.
The claim follows by Theorem~\ref{th:fib=Fib} and~\eqref{eq:Flm-via-flm}.
\end{proof}

In \cite[Section~2.4]{yga}, we gave explicit formulas for the
Fibonacci polynomials 
$F[\alpha]$ associated with the root systems of
classical types~$ABCD$. The Fibonacci polynomials of exceptional types
can be calculated on a computer using the defining
recurrence~\eqref{eq:flm-def}.
Combining these results with~\eqref{eq:Flm-via-flm},
one can compute all polynomials~$F_{\ell;m}\,$.
Corollary~\ref{cor:xjt=F/F} then yields explicit formulas
expressing an arbitrary cluster variable in an arbitrary cluster
algebra of finite type in terms of the 
bipartite initial seed.
These formulas are too bulky to be included in this paper:
as mentioned in \cite[p.~999]{yga}, the monomial expansion of a
Fibonacci polynomial of type~$E_8$ can have up to 26908 terms.

\section{Universal coefficients}
\label{sec:universal-coeffs}

In this section, we show that among all cluster algebras of a
given finite (Cartan-Killing) type, there is one with the
``universal" system of coefficients:
any other such algebra is obtained from this one
by a coefficient specialization.
To make this precise, we will need the following definition.

\begin{definition}
\label{def:coeff-specialization}
Let $\Acal$ and $\overline \Acal$ be cluster algebras
of the same rank~$n$ over the coefficient semifields $\PP$
and~$\overline \PP$, respectively,
with the respective families of cluster variables $(x_{i;t})_{i \in
  [1,n], t \in \TT_n}$ and
$(\overline x_{i;t})_{i \in [1,n], t \in \TT_n}$.
We say that $\overline \Acal$ is obtained from $\Acal$ by a
\emph{coefficient specialization} if:
\begin{enumerate}
\item $\Acal$ and $\overline \Acal$ have the same exchange
matrices $B_t = \overline B_t$ at each vertex $t \in \TT_n$;

\item there is a homomorphism of multiplicative
groups $\varphi: \PP \to \overline \PP$ that extends to a (unique) ring
homomorphism $\varphi: \Acal \to \overline \Acal$ such that
$\varphi(x_{i;t}) = \overline x_{i;t}$ for all~$i$ and~$t$.
\end{enumerate}
With some abuse of notation, we denote both maps $\varphi: \PP \to \overline \PP$
and $\varphi: \Acal \to \overline \Acal$ by the same symbol, and
refer to each of them as a \emph{coefficient specialization}.
\end{definition}

Recall from Section~\ref{sec:cluster-prelims}
that the cluster algebra structure is completely determined
by the underlying $Y$-pattern $t \mapsto (\yy_t,B_t)$, 
which is in turn completely determined by a $Y$-seed
$(\yy_{t_0},B_{t_0})$ attached to an arbitrary ``initial'' vertex
$t_0\in\TT_n$.
In the language of $Y$-patterns, Definition~\ref{def:coeff-specialization} can
be rephrased as follows.

\begin{proposition}
\label{pr:Y-base-change}
Let $\Acal$ and $\overline \Acal$ be cluster algebras
over the coefficient semifields $\PP$ and~$\overline \PP$,
respectively, sharing the same exchange matrices~$B_t\,$.
Let  the underlying $Y$-patterns for $\Acal$ and~$\overline \Acal$
be $t \mapsto (\yy_t,B_t)$ and $t \mapsto (\overline \yy_t,B_t)$,
respectively.
%
A homomorphism of multiplicative groups $\varphi: \PP \to \overline \PP$
is a coefficient specialization (that is, satisfies condition~{\rm
  (2)}
in Definition~\ref{def:coeff-specialization}) if and only if
\begin{equation}
\label{eq:phi-y-ybar}
\varphi(y_{j;t}) = \overline y_{j;t}, \quad
\varphi(y_{j;t} \oplus 1) = \overline y_{j;t} \oplus 1
\end{equation}
for all $j \in [1,n]$ and $t \in \TT_n$.
\end{proposition}

\begin{proof}
We start by restating the condition \eqref{eq:phi-y-ybar}
in terms of the coefficients appearing in the exchange relations.
Following \eqref{eq:p-through-y}, we denote by
\begin{equation}
\label{eq:pkt-through-ykt}
p_{j;t}^+ = \frac{y_{j;t}}{y_{j;t} \oplus 1}, \quad
p_{j;t}^- = \frac{1}{y_{j;t} \oplus 1}
\end{equation}
the two coefficients in the exchange relation for a cluster variable $x_{j;t}$
at a vertex~$t$ (cf.~\eqref{eq:exch-rel-usual});
also let $\overline p^\pm_{j;t}$ denote the counterparts of the
coefficients $p^\pm_{j;t}$ in the cluster algebra~$\overline \Acal$.
In view of \eqref{eq:pkt-through-ykt}, for a homomorphism of
multiplicative groups $\varphi : \PP \to \overline \PP$,
condition \eqref{eq:phi-y-ybar} is equivalent to requiring
\begin{equation}
\label{eq:phi-p-pbar}
\varphi(p^\pm_{j;t}) = \overline p^\pm_{j;t}
\end{equation}
for all $j \in [1,n]$ and $t \in \TT_n$.

Let $\xx_0$ and $\overline \xx_0$ be the clusters of
$\Acal$ and $\overline \Acal$ at some initial vertex $t_0$.
Since the elements of $\xx_0$ are algebraically independent over
the ground ring $\ZZ \PP$, any homomorphism of multiplicative groups
$\varphi: \PP \to \overline \PP$ extends to a (unique) ring
homomorphism $\varphi: \ZZ \PP[\xx_0^{\pm 1}] \to
\ZZ \overline \PP[\overline \xx_0^{\pm 1}]$ by setting
$\varphi(x_{i;t_0}) = \overline x_{i;t_0}$ for all~$i$.
Recall that by Theorem~\ref{th:Laurent-phenom},
$\Acal \subset  \ZZ \PP[\xx_0^{\pm 1}]$ and
$\overline \Acal \subset  \ZZ \overline \PP[\overline \xx_0^{\pm 1}]$.
In order to express a cluster variable at any vertex $t_1$ as a Laurent
polynomial in the elements of the initial cluster at~$t_0$, one
only needs to iterate the exchange relations of the form
\eqref{eq:exchange-rel-xx} along the path connecting $t_0$ and~$t_1$.
Since the elements of every cluster are algebraically independent,
the condition that $\varphi(x_{i;t}) = \overline x_{i;t}$ for all~$i$ and~$t$
is equivalent to the property that~$\varphi$ sends every relation
\eqref{eq:exchange-rel-xx} in~$\Acal$ to the
corresponding relation in $\overline \Acal$,
that is, satisfies the equivalent conditions \eqref{eq:phi-y-ybar}
and \eqref{eq:phi-p-pbar}.
\end{proof}

One way to get a coefficient specialization is to take
any semifield homomorphism $\varphi:\PP \to \overline \PP$.
For example, if $\Acal$ is a cluster algebra of geometric type
over a semifield $\Trop(u_j:j\in J)$, then for any subset $I\subset
J$, there is a semifield homomorphism $\varphi:\Trop(u_j:j\in
J)\to\Trop(u_j:j\in I)$ defined by setting $\varphi(u_j)=u_j$ for
$j\in I$ and $\varphi(u_j)=1$ for $j\notin I$.

\begin{definition}
\label{def:universal-coeffs}
We say that a cluster algebra~$\Acal$
(and the corresponding $Y$-pattern $t \mapsto (\yy_t,B_t)$)
has \emph{universal coefficients}
if every cluster algebra with the same family of exchange matrices~$(B_t)$
is obtained from~$\Acal$ by a unique coefficient specialization.
\end{definition}

By the standard universality argument, for a
given family of exchange matrices~$(B_t)$, a cluster algebra with
universal coefficients is \emph{unique} in the following sense:
for any two such algebras~$\Acal$ and $\overline \Acal$, there is a
canonical isomorphism of
multiplicative groups of their coefficient semifields that
extends to a canonical ring isomorphism $\Acal \to \overline \Acal$
sending all cluster variables and all coefficients in the exchange relations
for~$\Acal$ to their respective counterparts for $\overline \Acal$.
(A subtle point: the definition of universal
coefficients does \emph{not} imply the uniqueness of auxiliary
addition in the coefficient semifield.)

Note that the existence of a cluster algebra with universal coefficients
is by no means clear.
Below we prove it for the finite type case.
To state the result,
let us again recall the main theorem of~\cite{ca2}:
any cluster algebra~$\Acal$ of finite type has an
initial seed with a bipartite exchange matrix~$B$
such that the Cartan counterpart $A=A(B)$ is a Cartan matrix of finite
type.
That is, $A$ is associated with a \emph{finite} root system~$\Phi$,
uniquely defined by~$\Acal$ (up to isomorphism).

\begin{theorem}
\label{th:universal}
For any bipartite 
matrix $B$ such that $A=A(B)$ is a Cartan matrix of finite type,
there exists a universal cluster algebra~$\Auniv(B)=\Auniv(\Phi)$
having~$B$ as the initial exchange matrix.
Specifically, $\Auniv(B)$ can be realized as follows:
\begin{itemize}
\item
As a coefficient semifield for $\Auniv(B)$, take
\begin{equation}
\label{eq:semifield-universal-coeff}
\PP^{\rm univ} = \Trop (p[\alpha^\vee]: \alpha^\vee \in \Phi^\vee_{\geq -1}),
\end{equation}
the tropical semifield whose generators $p[\alpha^\vee]$ are labeled
by the set $\Phi^\vee_{\geq -1}$ of ``almost positive coroots"
(i.e., positive roots and negative simple roots in the
root system dual to~$\Phi$).
\item
Define the initial $Y$-seed $(\yy_0^{\rm univ},B)$, with
$\yy_0^{\rm univ} = (y_{1;0},\dots,y_{n;0})$, by
\begin{equation}
\label{eq:y-universal-initial}
y_{j;0} = \prod_{\alpha^\vee \in \Phi^\vee_{\geq -1}}
p[\alpha^\vee]^{\varepsilon(j) [\alpha^\vee : \alpha_j^\vee]} ,
\end{equation}
where $\varepsilon$ is the sign function from~\eqref{eq:bipartite-B},
and $[\alpha^\vee : \alpha_j^\vee]$ stands
for the coefficient of $\alpha_j^\vee$ in the expansion of
the coroot $\alpha^\vee$ in the basis of simple coroots.
\end{itemize}
\end{theorem}

The rank~$2$ special case of Theorem~\ref{th:universal} is implicitly
contained in \cite[Proposition~6.4 and Remark~6.5]{ca1}.

Before providing a proof of the theorem, we illustrate
it by an example.

\pagebreak[3]

\begin{example}[\emph{Universal cluster algebra of type~$A_2$}]
Let $B=\left[\begin{smallmatrix}
0 & 1\\
-1 & 0
\end{smallmatrix}
\right]$
as in Example~\ref{example:A2-xy-bipartite}.
Thus, $\varepsilon(1)=1$ and $\varepsilon(2)=-1$.
Since the root system $\Phi$ of type~$A_2$ is simply-laced,
we can identify the coroots with the roots.
By \eqref{eq:semifield-universal-coeff}, the universal coefficient
semifield~$\PP^{\rm univ}$ is the tropical semifield with $5$~generators:
\begin{equation}
\label{eq:PP-univ-A2}
\PP = \Trop
(p[-\alpha_1],p[-\alpha_2],p[\alpha_1],p[\alpha_1+\alpha_2],p[\alpha_2]).
\end{equation}
Formula \eqref{eq:y-universal-initial} gives
\begin{equation}
\label{eq:y0-univ-A2}
y_{1;0}=\frac{p[\alpha_1]\,p[\alpha_1+\alpha_2]}{p[-\alpha_1]}\,,\qquad
y_{2;0}=\frac{p[-\alpha_2]}{p[\alpha_2]\,p[\alpha_1+\alpha_2]}\,.
\end{equation}
Thus, $\Auniv(B)$ is the cluster algebra of geometric type defined by
the extended exchange matrix (see \eqref{eq:tilde-B})
\[
\tilde B=
\begin{bmatrix}
0  &  1\\
-1 &  0\\
-1 &  0\\
 0 &  1\\
 1 &  0\\
 1 & -1\\
 0 & -1
\end{bmatrix}\,.
\]
(The rows 3--7 of the matrix $\tilde B$ correspond to the variables
$p[\alpha]$ ($\alpha\in\Phi_{\geq -1}$) in the order listed
in~\eqref{eq:PP-univ-A2}.)
We then calculate the rest of the coefficients $y_{j;t}$ by substituting the
above expressions for $y_1$ and $y_2$ into the third column of
Table~\ref{table:seeds-A2}, and computing the resulting expressions
in~$\PP^{\rm univ}$. (It is useful to note the following: if $Y$ and $Y'$ are
two disjoint monomials in the generators of~$\PP^{\rm univ}$, then
$\frac{Y}{Y'}\oplus 1=\frac{1}{Y'}$.)
The results are shown in the third column of
Table~\ref{table:universal-coeff-A2}.

\begin{table}[ht]
\begin{equation*}
\begin{array}{|c|c|cc|cc|}
\hline
&&&&&\\[-4mm]
t & B_t & y_{1;t} & y_{2;t}  & x_{1;t} & x_{2;t}  \\[1mm]
\hline
&&&&&\\[-3mm]
0 & \begin{bmatrix}0&1\\-1&0\end{bmatrix} &
\dfrac{p[\alpha_1]\,p[\alpha_1+\alpha_2]}{p[-\alpha_1]}&
\dfrac{p[-\alpha_2]}{p[\alpha_2]\,p[\alpha_1+\alpha_2]}
&x[-\alpha_1] & x[-\alpha_2]
\\[4.5mm]
\hline
&&&&&\\[-3mm]
1 & \begin{bmatrix}0&-1\\1&0\end{bmatrix} &
\dfrac{p[\alpha_1]}{p[-\alpha_1]\,p[\alpha_2]} &
\dfrac{p[\alpha_2]\,p[\alpha_1+\alpha_2]}{p[-\alpha_2]}
& x[-\alpha_1]
& x[\alpha_2]
\\[4.5mm]
\hline
&&&&&\\[-3mm]
2 & \begin{bmatrix}0&1\\-1&0\end{bmatrix} &
\dfrac{p[-\alpha_1]\,p[\alpha_2]}{p[\alpha_1]}&
\dfrac{p[\alpha_1+\alpha_2]}{p[-\alpha_1]\,p[-\alpha_2]}&
x[\alpha_1+\alpha_2] &
x[\alpha_2]
\\[4.5mm]
\hline
&&&&&\\[-3mm]
3 & \begin{bmatrix}0&-1\\1&0\end{bmatrix} &
\dfrac{p[\alpha_2]}{p[-\alpha_2]\,p[\alpha_1]}&
\dfrac{p[-\alpha_1]\,p[-\alpha_2]}{p[\alpha_1+\alpha_2]}&
x[\alpha_1+\alpha_2]  &
x[\alpha_1]
\\[4.5mm]
\hline
&&&&&\\[-3mm]
4 & \begin{bmatrix}0&1\\-1&0\end{bmatrix} &
\dfrac{p[-\alpha_2]\,p[\alpha_1]}{p[\alpha_2]} &
\dfrac{p[-\alpha_1]}{p[\alpha_1]\,p[\alpha_1+\alpha_2]} &
x[-\alpha_2] &
x[\alpha_1]
\\[4.5mm]
\hline
&&&&&\\[-3mm]
5 & \begin{bmatrix}0&-1\\1&0\end{bmatrix} &
\dfrac{p[-\alpha_2]}{p[\alpha_2]\,p[\alpha_1+\alpha_2]} &
\dfrac{p[\alpha_1]\,p[\alpha_1+\alpha_2]}{p[-\alpha_1]} &
x[-\alpha_2] &
x[-\alpha_1]
\\[4.5mm]
\hline
\end{array}
\end{equation*}
\caption{Universal coefficients in type~$A_2$}
\label{table:universal-coeff-A2}
\end{table}

We conclude this example by writing down the exchange relations in the
universal cluster algebra of type~$A_2\,$.
The coefficients in each of these relations are obtained by splitting
the appropriate Laurent monomial in the third column of
Table~\ref{table:universal-coeff-A2} into its numerator and
its denominator.
To have a consistent notation for the $5$ cluster variables
in~$\Auniv(B)$, we use the denominator parametrization:
let $x[\alpha]$ denote the cluster variable
whose denominator vector is~$\alpha$.
See the fourth column of Table~\ref{table:universal-coeff-A2}.

With the notation agreed upon above, the exchange relations in~$\Auniv(B)$
are:
\begin{align}
\label{eq:exch-univ-A2-1}
x[-\alpha_2]\,x[\alpha_2]
&= p[-\alpha_2]\,x[-\alpha_1]+p[\alpha_2]\,p[\alpha_1+\alpha_2]\,,\\[.1in]
x[-\alpha_1]\,x[\alpha_1+\alpha_2]
&= p[\alpha_1]\,x[\alpha_2]+ p[-\alpha_1]\,p[\alpha_2]\,,\\[.1in]
x[\alpha_2]\,x[\alpha_1]
&= p[\alpha_1+\alpha_2]\,x[\alpha_1+\alpha_2]+ p[-\alpha_1]\,p[-\alpha_2]\,,\\[.1in]
x[\alpha_1+\alpha_2]\,x[-\alpha_2]
&= p[\alpha_2]\,x[\alpha_1]+ p[-\alpha_2]\,p[\alpha_1]\,,\\[.1in]
\label{eq:exch-univ-A2-5}
x[\alpha_1]\,x[-\alpha_1]
&= p[-\alpha_1]\,x[-\alpha_2]+p[\alpha_1]\,p[\alpha_1+\alpha_2]\,.
\end{align}
\end{example}
These exchange relations (thus, implicitly, the algebra~$\Auniv(B)$)
have already appeared in \cite[Example~7.6]{ca1}\footnote{
\emph{Correction:} in \cite[Figure~3]{ca1},
the exchange relation between $y_2$ and $y_5$ contained a typo.
}.

\pagebreak[3]

\begin{proof}[Proof of Theorem~\ref{th:universal}]
In the following argument, we fix the bipartite exchange matrix~$B$ of finite type
associated with a root system~$\Phi$, and work with the class of cluster algebras~$\Acal$
having~$B$ attached to an initial seed.
We start by taking an inventory of the exchange relations in~$\Acal$.
By Proposition~\ref{pr:fintype-belt},
the denominator vectors with respect to the initial cluster
provide a bijection
$\alpha \mapsto x[\alpha]$ between the set $\Phi_{\geq -1}$ of ``almost positive" roots
and the set of all cluster variables in~$\Acal$;
furthermore, this bijection extends to a bijection
$\gamma \mapsto x[\gamma]$ between the root lattice $Q$ and the set of all
cluster monomials in~$\Acal$.
Using these bijective parameterizations,
we can write every exchange relation in~$\Acal$ in the form
\begin{equation}
\label{eq:exchange-general-2}
x[\beta]\, x[\beta'] = p_1 x[\gamma_1] + p_2 x[\gamma_2],
\end{equation}
where  $\beta, \beta' \in \Phi_{\geq -1}$, $\gamma_1, \gamma_2 \in Q$, and
$p_1$ and $p_2$ belong to the coefficient semifield~$\PP$.
These exchange relations have the following properties:
\begin{itemize}
\item
The pairs $\{\beta, \beta'\}$ appearing in the left-hand sides
depend only on the root system~$\Phi$, that is, are the same for all choices of
coefficients; see \cite[Proposition~3.5 and Corollary~4.4]{ca2},
where such pairs are called \emph{exchangeable}.
\item
The pair $\{\gamma_1, \gamma_2\}$ on the right
is uniquely determined by $\{\beta, \beta'\}$; so it also
does not depend on the choice of the coefficient system
in~$\Acal$. (See \cite[(5.1)]{ca2}, where $\gamma_1$ and
$\gamma_2$ are explicitly expressed in terms of $\beta$ and $\beta'$.)
\end{itemize}
In view of these properties, we will denote the coefficients
in \eqref{eq:exchange-general-2} as follows:
\begin{equation}
\label{eq:p(beta,beta:gamma)}
p_i = p(\{\beta, \beta'\}; \gamma_i) \qquad (i = 1,2).
\end{equation}
This notation becomes slightly ambiguous if $\gamma_1 = \gamma_2$.
But this can only happen if $\beta = \alpha_i$ and $\beta' = -\alpha_i$ are
the two roots of an irreducible component of type~$A_1$ in~$\Phi$,
in which case $\gamma_1 = \gamma_2 = 0$; all the arguments below
will be  valid in this case as well, even if both coefficients will
carry the same notation.

Among the coefficients $p(\{\beta, \beta'\}; \gamma)$,
the following ones will be of special importance.

\begin{definition}
\label{def:primitive-relations}
A coefficient $p(\{\beta, \beta'\}; \gamma_1)$ in
\eqref{eq:exchange-general-2}--\eqref{eq:p(beta,beta:gamma)}
is called \emph{primitive} if $\gamma_2 = 0$,
i.e., the ``opposite" monomial $x[\gamma_2]$ is equal to~$1$.
\end{definition}

The following lemma provides several equivalent
descriptions of exchange relations with primitive coefficients.

\begin{lemma}
\label{lem:primitive-relations}
For an exchange relation \eqref{eq:exchange-general-2}
the following conditions are equivalent:
\begin{enumerate}
\item One of the coefficients $p(\{\beta, \beta'\}; \gamma_1)$
and $p(\{\beta, \beta'\}; \gamma_2)$ is primitive.
\item The relation in question 
appears on the
bipartite belt, that is,
$\{\beta, \beta'\} = \{\dd(j;m-1), \dd(j;m+1)\}$ for
some $j \in [1,n]$ and~$m \in \ZZ$ (see
Theorem~\ref{th:denominators-bipartite}).
\item $\{\beta, \beta'\} = \{\tau_+ \alpha, \tau_- \alpha\}$
for some $\alpha \in \Phi_{\geq -1}$.
\end{enumerate}
\end{lemma}

Recall that $\tau_+$ and $\tau_-$ appearing in condition (3)
are the involutive permutations of $\Phi_{\geq -1}$
given by \eqref{eq:tau-action}; they extend to piecewise-linear involutions of the
root lattice~$Q$ via
\eqref{eq:tau-pm-tropical-1}--\eqref{eq:tau-pm-tropical-2}.

\begin{proof}
The equivalence of (2) and (3) is clear from the formulas
\eqref{eq:dems-m-positive-tau}--\eqref{eq:dems-m-negative-tau}.
The implication (2)$\Rightarrow$(1) follows from the
definition of the bipartite belt.
To prove (1)$\Rightarrow$(2), we use the formula given in
\cite[(5.1)]{ca2},
which gives explicit expressions for $\gamma_1$ and $\gamma_2$ in
terms of $\{\beta, \beta'\}$.
As a special case of \emph{loc.~cit.}, condition (1) is equivalent
to the following: $\sigma (\beta) + \sigma (\beta') = 0$ for some transformation~$\sigma$
belonging to the group generated by~$\tau_+$ and~$\tau_-$.
Since both $\sigma (\beta)$ and $\sigma (\beta')$ belong to $\Phi_{\geq -1}$,
the only way they can add up to~$0$ is if say $\sigma (\beta) = - \alpha_j$
and $\sigma (\beta') = \alpha_j$ for some~$j$.
Again remembering
\eqref{eq:dems-m-positive-tau}--\eqref{eq:dems-m-negative-tau}, it is
easy to see that in that case, $\{\beta, \beta'\} = \{\dd(j;m-1),
\dd(j;m+1)\}$ for some~$m \in \ZZ$, and we are done.
\end{proof}

\begin{lemma}
\label{lem:primitive-coeffs-universal}
In the algebra $\Auniv(B)$ given by
\eqref{eq:semifield-universal-coeff} and
\eqref{eq:y-universal-initial},
the primitive coefficients are precisely
the generators $p[\alpha^\vee]$ of $\PP^{\rm univ}$.
More specifically, if $\{\beta, \beta'\} = \{\tau_+(\alpha), \tau_-(\alpha)\}$
for some $\alpha \in \Phi_{\geq -1}$ (see
Lemma~\ref{lem:primitive-relations}(3)), then the corresponding
primitive coefficient in \eqref{eq:exchange-general-2} is equal to
$p[\alpha^\vee]$, where $\alpha^\vee$ is the coroot associated with a root $\alpha$.
\end{lemma}

\begin{proof}
Let us explicitly solve the $Y$-system
\eqref{eq:ym-exchange} in the semifield~$\PP^{\rm univ}$ defined
by~\eqref{eq:semifield-universal-coeff} with the initial
conditions
\begin{equation}
\label{eq:y-universal-initial-2}
\begin{array}{rcl}
y_{i;-1}\!\! &=& \!\!\displaystyle\prod_{\alpha^\vee \in \Phi^\vee_{\geq -1}}
p[\alpha^\vee]^{- [\alpha^\vee : \alpha_i^\vee]}
\qquad (\varepsilon(i) = +1),\\[.3in]
y_{j;0}\!\! &=& \!\!\displaystyle\prod_{\alpha^\vee \in \Phi^\vee_{\geq -1}}
p[\alpha^\vee]^{- [\alpha^\vee : \alpha_j^\vee]}
\qquad (\varepsilon(j) = -1);
\end{array}
\end{equation}
by \eqref{eq:y-parity}, these conditions are equivalent
to~\eqref{eq:y-universal-initial}.
We will need the analogue of \eqref{eq:tau-pm-tropical-1}--\eqref{eq:tau-pm-tropical-2}
for the dual root system, which can be written as
\begin{equation}
\label{eq:tau-pm-tropical-dual}
[\tau_\varepsilon  \alpha^\vee: \alpha^\vee_j] =
\begin{cases}
- [\alpha^\vee: \alpha^\vee_j] -
\sum_{i \neq j}  a_{ij} [[\alpha^\vee: \alpha^\vee_i]]_+
&
\text{if $\varepsilon (j) = \varepsilon$;} \\[.1in]
[\alpha^\vee: \alpha^\vee_j] & \text{otherwise.}
\end{cases}
\end{equation}
Using \eqref{eq:tau-pm-tropical-dual} and proceeding by induction,
it is straightforward to verify that the solution of
\eqref{eq:ym-exchange} in~$\PP^{\rm univ}$ with the initial conditions
\eqref{eq:y-universal-initial-2} is given by
\begin{equation}
\begin{array}{rcl}
\label{eq:y-universal-solution}
y_{i;-r-1}\!\! &=&\!\! \displaystyle\prod_{\alpha^\vee \in \Phi^\vee_{\geq -1}}
p[\alpha^\vee]^{- [\tau_{\varepsilon(i)}^{(r)} \alpha^\vee : \alpha_i^\vee]}
\qquad (\varepsilon(i) = (-1)^r),\\
y_{j;r}\!\! &=&\!\! \displaystyle\prod_{\alpha^\vee \in \Phi^\vee_{\geq -1}}
p[\alpha^\vee]^{- [\tau_{\varepsilon(j)}^{(r)}\alpha^\vee : \alpha_j^\vee]}
\qquad (\varepsilon(j) = (-1)^{r-1})
\end{array}
\end{equation}
for all $r \geq 0$, where we abbreviate
$\tau_\varepsilon^{(r)} =
\underbrace{\tau_\varepsilon \tau_{-\varepsilon} \tau_\varepsilon
  \cdots}_{r \text{~factors}}$.

Substituting the expressions given by
\eqref{eq:y-universal-solution} into \eqref{eq:xm-exchange}, we
obtain explicit formulas for exchange relations in $\Auniv$ that
belong to the bipartite belt.
By comparing \eqref{eq:y-universal-solution} with the expressions
\eqref{eq:dems-m-positive-tau} and \eqref{eq:dems-m-negative-tau}
for the denominator vectors, it is easy to see that if
the exchange relation in question has
$\{\beta, \beta'\} = \{\tau_+(\alpha), \tau_-(\alpha)\}$ then the
corresponding primitive coefficient $p(\{\beta, \beta'\}; \gamma_1)$ is
equal to $p[\alpha^\vee]$, completing the proof.
(There is also a nice formula for the other coefficient in the same relation:
\begin{equation}
\label{eq:constant-term-universal}
p(\{\beta, \beta'\}; \gamma_2) = \prod_{\beta^\vee} p[\beta^\vee]^{(\beta^\vee \|
\alpha^\vee)},
\end{equation}
where $\beta^\vee$ runs over $\Phi_{\geq -1}^\vee$,
and $(\beta^\vee \| \alpha^\vee)$ is the compatibility degree
introduced in \cite[Section~3.1]{yga}. We do not need this
formula for the proof.)
\end{proof}

We continue the proof of Theorem~\ref{th:universal}.
By Lemma~\ref{lem:primitive-coeffs-universal},
every generator $p[\alpha^\vee]$ of $\PP^{\rm univ}$ appears
as a coefficient in an exchange relation on the
bipartite belt.
It follows that the coefficients in the exchange
relations generate $\PP^{\rm univ}$ as a multiplicative group.
In view of \eqref{eq:phi-p-pbar}, we arrive at the following conclusion:

\begin{itemize}
\item
If a cluster algebra $\Acal$ is obtained from $\Auniv$ by a
coefficient specialization~$\varphi$, then~$\varphi$ is \emph{unique}.
\end{itemize}

It remains to prove the \emph{existence} of a
coefficient specialization $\varphi: \Auniv \to \Acal$. \linebreak[3]
To do this, we introduce the following terminology.

\begin{definition}
A \emph{multiplicative coefficient identity} (MCI for
short) is a multiplicative relation among the coefficients
$p(\{\beta, \beta'\}; \gamma)$ that holds in \emph{all}
cluster algebras of a given type, independently of the choice of
coefficients.
\end{definition}

The remaining (existence) part of Theorem~\ref{th:universal} is an immediate consequence of
Lemma~\ref{lem:primitive-coeffs-universal} in conjunction with the following lemma.

\begin{lemma}
\label{lem:key}
For every non-primitive coefficient $p(\{\beta, \beta'\}; \gamma)$, there is an MCI
expressing it as a product of primitive ones (not necessarily distinct).
\end{lemma}

\begin{proof}
We start by deducing the desired statement from the following formally weaker
assertion:
\begin{align}
\label{eq:non-primitive-factors}
&\text{For every non-primitive coefficient $p(\{\beta, \beta'\};
\gamma)$, there is an MCI that}\\
\nonumber
&\text{expresses it as a product of two or more
other coefficients.}
\end{align}
Indeed, let us assume that \eqref{eq:non-primitive-factors}
holds, and fix an MCI as in \eqref{eq:non-primitive-factors} for
each non-primitive coefficient $p(\{\beta, \beta'\}; \gamma)$.
To prove that every non-primitive coefficient~$p$ factors into a
product of primitive ones, we write a chosen MCI for~$p$ as
$p = p_1 \cdots p_s$, then substitute for each~$p_i$ the product
given by the corresponding MCI, and continue in the same way.
To show that this process terminates, consider the following directed graph~$\Gamma$:
the vertices of~$\Gamma$ correspond to all the coefficients $p(\{\beta, \beta'\}; \gamma)$
in the exchange relations, and the edges are of the form
$p \to p'$ whenever~$p'$ is one of the factors in the chosen MCI for~$p$.
Since the sinks of~$\Gamma$ are precisely the primitive
coefficients, it suffices to show that~$\Gamma$ has
no oriented cycles.
Suppose on the contrary that $p_1 \to p_2 \to \cdots \to p_s \to
p_1$ is an oriented cycle in~$\Gamma$.
Multiplying the chosen MCIs along this cycle
and canceling $p_1\cdots p_s$, we obtain an MCI that
expresses~$1$ as a product of several coefficients:
$1=q_1 \cdots q_r$.
Being an MCI, this identity must be satisfied in \emph{any}
cluster algebra of the given type.
Returning for a moment to the nomenclature of coefficients used in
the preceding sections, we can write $q_1$ in the form
$q_1 = p^+_{j;t}$ for some $j \in [1,n]$ and $t \in \TT_n$.
Now let $\Aprin$ be the algebra (of the same finite type)
with principal coefficients at~$t$.
Its coefficient semifield is $\Trop(y_1, \dots, y_n)$, with
$q_1 = y_j$, while the rest of the factors $q_2, \dots, q_r$ are
monomials in $y_1, \dots, y_n$ (with nonnegative exponents).
Hence the equality $1=q_1 \cdots q_r$ does not hold in
$\Aprin$ and so cannot be an MCI.
This contradiction proves that Lemma~\ref{lem:key} is indeed a consequence
of \eqref{eq:non-primitive-factors}.

The proof of \eqref{eq:non-primitive-factors} follows from the results
in \cite{ca1}: applying the reduction procedure in
\cite[Section~2]{ca1}, we reduce the statement to the rank~$2$
case, and then the desired MCI can be obtained as one of the identities
\cite[(6.11)--(6.13)]{ca1}.
For the convenience of the reader, we provide a little more details here.
Let us start with a relation of the form \eqref{eq:exchange-general-2}
with the coefficient~$p_1$ non-primitive.
This means that the cluster monomial $x[\gamma_2]$ is non-trivial,
hence contains some cluster variable $x[\beta'']$.
Let $\xx$ be a cluster containing $x[\beta]$ and such that
$\xx - \{x[\beta]\} \cup \{x[\beta']\}$ is another cluster.
Let $\zz \subset \xx$ be the $(n-2)$-element subset
$\xx - \{x[\beta], x[\beta'']\}$.
Following \cite[Section~2]{ca1}, we can consider the rank~$2$
cluster algebra setup with the initial cluster
$\{x[\beta], x[\beta'']\}$ obtained by ``freezing" the cluster
variables from $\zz$, i.e., viewing them as ``constants" or ``new coefficients."
This algebra is of one of the finite types $A_2, B_2$ or~$G_2$;
the type $A_1 \times A_1$ cannot occur because then all the
cluster variables occurring in $x[\gamma_1]$ or $x[\gamma_2]$
would belong to~$\zz$.
Its exchange graph is an $(h+2)$-cycle
(cf.\ Theorem~\ref{th:belt-periodic}), that is, a $5$-cycle
(resp., $6$-cycle, $8$-cycle) for the type~$A_2$ (resp.,
$B_2$,~$G_2$).
In determining the MCIs along this cycle, we can ignore the terms
from~$\zz$, which reduces the proof of \eqref{eq:non-primitive-factors}
to the rank~$2$ case, and brings us into the setup studied in
\cite[Section~6]{ca1}.

Let us rewrite the results in \emph{loc.\ cit.} in our current notation.
Assume that $\Phi$ is of one of the types $A_2$, $B_2$ or~$G_2$.
Let $\Phi_{\geq -1} = \{\beta_1, \dots, \beta_{h+2}\}$, where the
roots are ordered cyclically around the origin; for instance, in
type~$A_2$ we use the same ordering as in~\eqref{eq:PP-univ-A2}.
We use the convention that the
subscripts are considered modulo $(h+2)$.
Then the exchange relations and the corresponding MCIs can be
described as follows.

\smallskip

\noindent {\bf Type $A_2$.} The exchange relations have the form
\[
x[\beta_{m-1}]\, x[\beta_{m+1}] = q_m\, x[\beta_{m}] + r_m\,.
\]
For every~$m$, there is an MCI expressing the non-primitive
coefficient~$r_m$ as a product $r_m = q_{m+2} q_{m+3}\,$.

\smallskip

\noindent {\bf Type $B_2$.}
In this case, in the sequence
$\beta_1, \beta_2,\dots$, the long and short roots alternate.
The exchange relations have the form
\[
x[\beta_{m-1}]\, x[\beta_{m+1}] = q_m\, x[\beta_{m}]^{b_m} + r_m,
\]
where $b_m = 1$ (resp., $2$) if $\beta_{m-1}$ and $\beta_{m+1}$ are
short (resp., long).
The corresponding MCIs look as follows: $r_m = q_{m+2}\, q_{m+3}^{b_m}\, q_{m+4}$.

\smallskip

\noindent {\bf Type $G_2$.} Again, in the sequence
$\beta_1, \beta_2, \dots$, the long and short roots alternate.
The exchange relations have the form
\[
x[\beta_{m-1}]\, x[\beta_{m+1}] = q_m\, x[\beta_{m}]^{b_m} + r_m,
\]
where $b_m = 1$ (resp., $3$) if $\beta_{m-1}$ and $\beta_{m+1}$ are
short (resp., long).
The corresponding MCIs look as follows:
$$
r_m =
\begin{cases}
q_{m+2}\, q_{m+3}\, q_{m+4}^2 \,q_{m+5}\, q_{m+6}  & \text{if $\beta_{m-1}$ and $\beta_{m+1}$ are short};\\[.05in]
q_{m+2}\, q_{m+3}^3\, q_{m+4}^2\, q_{m+5}^3\, q_{m+6}
 &  \text{if $\beta_{m-1}$ and $\beta_{m+1}$ are long}.
\end{cases}
$$
This completes the proof of Lemma~\ref{lem:key}.
\end{proof}

Theorem~\ref{th:universal} is proved.
\end{proof}

\begin{remark}
Once it is proved that the cluster monomials form a $\ZZ \PP$-basis
in any cluster algebra~$\Acal$ of finite type
(see \cite[Conjecture~4.28]{cdm}),
Theorem~\ref{th:universal}
will imply that any such~$\Acal$ can be obtained by a ``base change" from
the universal cluster algebra $\Auniv$ of the same type: that is,
$\Acal$
is canonically isomorphic to $\Auniv \otimes_{\ZZ \PP^{\rm univ}}
\ZZ \PP$.
\end{remark}

\section{Index of conjectures}

Table~\ref{table:index-of-conjectures} gives pointers to
partial results which provide supporting evidence towards
conjectures proposed
in
this paper:
either for bipartite initial seeds, or under the additional assumption
of finite type.

\medskip

\begin{table}[ht]
\begin{tabular}{|c|c|c|}
\hline
                    & \textbf{bipartite case;} & \textbf{finite type;} \\
\textbf{conjecture} & \textbf{seeds/variables on} &\textbf{bipartite}\\
                     & \textbf{the bipartite belt} & \textbf{initial seed}\\
\hline
\hline
&&\\[-.15in]
Conjecture~\ref{con:exchange-graph-independence} & &
 \cite[Theorem~1.13]{ca2}\\
\cite[Conjecture~4.14(1)]{cdm} & &
\\[.05in] \hline
&&\\[-.15in]
Conjectures~\ref{con:conjectures-on-denoms}--\ref{con:conjectures-on-denoms-2}
 & Corollary~\ref{cor:conjectures-bipartite} &
\\[.05in] \hline
&&\\[-.15in]
Conjecture~\ref{con:hk-gk} & Proposition~\ref{pr:hk-gk-bipartite} &
\\[.05in] \hline
&&\\[-.15in]
Conjecture~\ref{con:g-thru-same-F} & Corollary~\ref{cor:conjectures-bipartite} &
\\[.05in] \hline
&&\\[-.15in]
Conjecture~\ref{con:signs-gi} & Remark~\ref{rem:signs-gi-bipartite}&
\\[.05in] \hline
&&\\[-.15in]
Conjecture~\ref{con:cl-mon-lin-indep} & & Theorem~\ref{th:cluster-monomials-linear-indep}\\
\cite[Conjecture~4.16]{cdm}&&
\\[.05in] \hline
&&\\[-.15in]
Conjecture~\ref{con:denominator-1} && \cite[Theorems 1.9, 1.13]{ca2}
\\[.05in] \hline
&&\\[-.15in]
Conjecture~\ref{con:sign-coherence-denom} &
 Corollary~\ref{cor:denominator-signs} &
\\[.05in] \hline
&&\\[-.15in]
Conjecture~\ref{con:denominator-2} & & Proposition~\ref{pr:fintype-belt} \\
\cite[Conjecture~4.17]{cdm} &&
\\[.05in] \hline
&&\\[-.15in]
Conjecture~\ref{con:g-vector-2} & & Corollary~\ref{cor:fock-goncharov-finite-type-bipartite}
\\[.05in] \hline
&&\\[-.15in]
Conjecture~\ref{con:g-transition} &
Proposition~\ref{pr:hk-gk-bipartite}
&
\\[.05in] \hline
&&\\[-.15in]
Conjecture~\ref{con:denominators-thru-F} & Corollary~\ref{cor:conjectures-bipartite} &
\\[.05in] \hline
\end{tabular}

\vspace{.1in}

\caption{Index of conjectures}
\label{table:index-of-conjectures}
\end{table}

\section*{Acknowledgments}

The authors thank the following institutions and colleagues whose
hospitality and support stimulated their work on this research project:
University of Bielefeld \linebreak[3]
(Claus Michael Ringel);
Frederick W.\ and Lois B.\ Gehring Fund and University of Michigan;
Institut Fourier (Michel Brion);
and Institut Mittag-Leffler (Anders Bj\"orner and Richard Stanley).
We also thank Vladimir Fock, Alexander Goncharov, and
Alek Vainshtein for stimulating discussions.

\end{document}